\documentclass[11pt]{article}
\usepackage[a4paper, margin = 2.2cm]{geometry}
\usepackage[english]{babel}
\usepackage{latexsym,amsmath,enumerate,graphics,enumerate,amsthm,tikz,hyperref,float}
\usepackage[mathscr]{euscript}
\usepackage[affil-it]{authblk}
\usepackage{enumerate,amsthm,dsfont,pstricks}
\usepackage{latexsym,amsmath,amssymb,amscd,wrapfig,graphicx}
\usepackage{wrapfig,graphicx,curves}

\newtheorem{theorem}{Theorem}[section]
\newtheorem{lemma}[theorem]{Lemma}
\newtheorem{proposition}[theorem]{Proposition}
\newtheorem{corollary}[theorem]{Corollary}
\newtheorem{problem}[theorem]{Problem}

\theoremstyle{definition}

\newtheorem{remark}[theorem]{Remark}
\theoremstyle{remark}

\let\phi=\varphi
\def\epsilon{\varepsilon}

\def\0{\mathbf{0}}

\newcommand{\comment}[1]{}

\numberwithin{equation}{section}

\let\epsilon=\varepsilon

\makeatletter
\def\@maketitle{%
  \newpage
  \null
  \vskip 1em%
  \begin{center}%
  \let \footnote \thanks
    {\Large\bfseries \@title \par}%
    \vskip 1.5em%
    {\normalsize
      \lineskip .5em%
      \begin{tabular}[t]{c}%
        \@author
      \end{tabular}\par}%
    \vskip 1em%
    {\normalsize \@date}%
  \end{center}%
  \par
  \vskip 1.5em}
\makeatother

\begin{document}

\title{\sc  \huge Horofunction compactifications and duality}

\author{Bas Lemmens and Kieran Power%
\thanks{Email: \texttt{B.Lemmens@kent.ac.uk} and \texttt{kcp6@kent.ac.uk}, Bas Lemmens gratefully acknowledges the support of the EPSRC (grant EP/R044228/1) and Kieran Power  gratefully acknowledges the support of the EPSRC (grant EP/V520093/1)}}
\affil{School of Mathematics, Statistics \& Actuarial Science,
University of Kent, \\Canterbury, CT2 7NX, United Kingdom}

\maketitle
\date{}

\begin{abstract}
We study the global topology and geometry of the horofunction compactification of certain simply connected smooth manifolds with a Finsler distance. The main goal is to show, for various classes of these spaces, that the  horofunction compactification is naturally homeomorphic to the closed unit ball of the dual norm of the norm in the tangent space (at the basepoint) that generates the Finsler distance.  We construct explicit homeomorphisms for a variety of spaces in three settings: bounded convex domains in $\mathbb{C}^n$ with the Kobayashi distance,  Hilbert geometries, and finite dimensional normed spaces.  For the spaces under consideration,  the horofunction boundary has an intrinsic partition into so called parts.  The natural connection with the dual norm arises through the fact that the homeomorphism maps each part in the horofunction boundary  onto the relative interior of a boundary face of the dual unit ball.  For normed spaces the connection between the global topology of the horofunction boundary and the dual norm was suggested by Kapovich and Leeb.  We confirm  this connection for Euclidean Jordan algebras equipped with the spectral norm.
\end{abstract}


{\small {\bf Subject Classification:} Primary 53C60; Secondary 32Q45, 46B20}
{\footnotesize\tableofcontents}

\section{Introduction}
A well known result in the theory of manifolds of nonpositive curvature says that if $M$ is a complete simply connected Riemannian manifold of nonpositive sectional curvature, then the horofunction compactification of $M$ is homeomorphic to the closed unit ball of the Hilbert norm in the tangent space at the basepoint  in $M$, e.g., \cite[Proposition 1.7.6]{Eb} or \cite{EO}. 
The main goal of this paper is to establish analogues of this result for various classes of simply connected smooth manifolds with a Finsler distance. 

Recall that a {\em Finsler distance} $d_F$  on a smooth manifold $M$ has an infinitesimal form  $F\colon TM\to \mathbb{R}$ on the tangent bundle $TM$, such that $d_F(x,y)$  is the infimum of {\em lengths},
\[
L(\gamma) = \int_0^1 F(\gamma(t),\gamma'(t))\,\mathrm{d}t,
\] 
 over piecewise $C^1$-smooth paths $\gamma\colon [0,1]\to M$  with $\gamma(0)=x$ and $\gamma(1)=y$. 

More explicitly, we analyse the following general question.
\begin{problem}\label{keyproblem} Suppose $M$ is a smooth simply connected manifold with a Finsler distance, such that  the restriction of $F$ to the tangent space $T_bM$ at $b$ is a norm. When does there exist a homeomorphism from the horofunction compactification of $M$, with basepoint $b$, onto the dual unit ball $B_1^*$ of the norm in $T_bM$ such that the homeomorphism maps each part in the horofunction  boundary onto the relative interior of a boundary face of $B_1^*$?
\end{problem}
It should be noted that the answer to the question may depend on the basepoint $b\in M$, as the norm in $T_bM$ may have a different facial structure for different basepoints.  However, the spaces we consider here are homogeneous in the sense that the facial structure of the compact convex set $\{v\in T_bM\colon F(b,v)\leq 1\}$ is the same for all $b\in M$. 

We confirm the existence of such a homeomorphism for a variety of manifolds in three settings:  bounded convex domains in $\mathbb{C}^n$ with the Kobayashi distance, finite dimensional normed spaces, and Hilbert geometries. For finite dimensional normed spaces the connection between the horofunction compactification and dual unit ball was suggested by Kapovich and Leeb \cite[Question 6.18]{KL}, who asked if for finite dimensional normed spaces the horofunction compactification is homeomorphic to the closed unit ball of the dual normed space.  This was confirmed by Ji and Schilling \cite{JS1,JS} for polyhedral normed spaces. 

For the Kobayashi distance on bounded convex domains, we consider product domains $B=B_1\times \cdots \times B_r$ in $\mathbb{C}^n$,  where each $B_i$ is the open unit ball of a norm with a strongly convex $C^3$-boundary. Prime examples are polydiscs. The Finsler structure, i.e., the infinitesimal  Kobayashi metric,  in the tangent space at $0$ is given by the norm $\|\cdot\|_B$ whose (open) unit ball is $B$, see \cite[Proposition 2.3.34]{Ab}. We show that the horofunction compactification is naturally homeomorphic to the closed ball of the dual norm of $\|\cdot\|_B$. 
For domains $D\subset \mathbb{C}^n$ with the Kobayashi distance, various conditions are known that imply that the identity map on $D$ extends as a homeomorphism from the horofunction compactification  of $D$ onto the norm closure $\mathrm{cl}\, D$, see \cite[Theorem 1.2]{AFGG} and \cite{BB,BGZ,Z}. These conditions typically involve strong convexity and smoothness properties of the domain. In our setting, however, the domains are not smooth, and the identity does not extend as a homeomorphism, as different geodesics converging to the same point in the norm boundary of the domain can yield different horofunctions. 

For finite dimensional normed spaces we will focus on  the finite dimensional Euclidean Jordan algebras equipped with the spectral norm, which are precisely the finite dimensional formally real JB-algebras \cite{AS1}. A prime example is the real vector space $\mathrm{Herm}_n(\mathbb{C})$ consisting of all $n\times n$ Hermitian matrices equipped with the norm,  $\|A\|=\max\{|\lambda|\colon \lambda \in\sigma(A)\}$.
The Jordan algebra structure allows us to give a complete characterisation of the horofunctions of these normed spaces.  We use this characterisation to provide a natural homeomorphism of the horofunction compactification onto the closed unit ball of the dual space.  To prove the results we  do not rely on the characterisation of the Busemann points in arbitrary normed spaces obtained by Walsh \cite{Wa2,Wa4}, but instead we exploit the Jordan algebra structure.  

For Hilbert geometries $(\Omega,d_H)$ we consider domains $\Omega$ that are obtained by intersecting a symmetric cone with hyperplane. A prime example is the space of  strictly-positive definite $n\times n$ Hermitian matrices with trace $n$. These Hilbert geometries are homogeneous in the sense that $\mathrm{Isom}(\Omega)$ acts transitively on $\Omega$, which ensures that the unit balls in the tangent spaces all have the same facial structure. We show, for these Hilbert geometries, that the horofunction compactification is naturally homeomorphic to the closed dual unit ball of the norm in the tangent space at the unit.  We use the cone  version of the Hilbert distance, see \cite{LNBook}, which provides a convenient way to analyse its Finsler structure \cite{Nu} and the dual of its norm. The horofunction compactification of these Hilbert geometries was determined in \cite[Theorem 5.6]{LLNW} and is naturally described in terms of the Euclidean Jordan algebra associated to the symmetric cone, which will be exploited in the analysis.

The origins of the horofunction compactification go back to Gromov \cite{Ba,Gr} who associated a boundary at infinity to any locally compact geodesic metric space. It has found  numerous applications in  diverse areas of mathematics including, geometric group theory \cite{BH}, noncommutative geometry \cite{Ri}, complex analysis \cite{Ab,AFGG, BB,BG,BGZ,Z}, Teichm\"uller theory \cite{DF,GJ,Ka2,Mi,Wa3}, dynamical systems and ergodic theory \cite{Be,GK,Ka1,LLNW} and in the study of Satake compactifications of noncompact type symmetric spaces \cite{HSWW,KL,Sc}.  A more general set up was discussed by Rieffel \cite{Ri}, who recasted the horofunction compactification of a locally compact geodesic metric space as a maximal ideal space of  a commutative $C^*$-algebra. Rieffel's set up works for any metric space, but if the metric space is not proper, then the embedding into its horofunction compactification need not be a homeomorphism.

The horofunction compactification is a particularly powerful tool to study isometry groups of metric spaces and isometric embeddings between metric spaces, see \cite{L1,LW,Wa1, Wa4}. Especially useful  in this context are the so called Busemann points in the horofunction compactification, which are limits of almost geodesics. They were introduced by Rieffel \cite{Ri}, who asked whether every horofunnction is a Busemann point in a finite dimensional normed space. Walsh \cite{Wa2} gave a complete solution to this problem and found necessary and sufficient conditions for a finite dimensional normed to have the property that all horofunctions are Busemann points.  

In the metric spaces under consideration in this paper, all horofunctions are Busemann points. On the set of Busemann points one can define a metric known as the detour distance \cite{AGW,LW}, which partitions the set of Busemann points into parts consisting of Busemann points that have finite detour distance to each other. So, for the metric spaces $M$ in this paper, the horofunction compactification is the disjoint union of $M$ and the parts in the horofunction boundary.  If two Busemann points have finite detour distance, it means that the corresponding almost geodesics are in some sense asymptotic. Moreover, any isometry on the  metric space $M$ induces an isometry on the set of Busemann points under the detour distance.  In each of our settings we will give an explicit homeomorphism that maps  $M$ onto the interior of the closed dual unit ball, and each part in the horofunction boundary of $M$ onto the relative interior of a boundary face of the dual unit ball.  It is this property of the homeomorphism that {\em naturally} connects the global topology of the horofunction compactification to the closed unit {\em dual} ball in each of our spaces.

In general it is hard to determine the horofunction compactification explicitly, and only  in relatively few spaces has this been accomplished, even in the context of normed spaces. We give an incomplete list of results in this direction.  For $\mathrm{CAT}(0)$ spaces the horofunction compactification is well understood, see \cite[Chapter II.8]{BH} and coincides with the visual boundary.  At present the horofunction compactification has been determined explicitly for a variety of normed spaces. Guti\`errez \cite{Gu1,Gu2,Gu3} computed the horofunction  compactification of several classes of $L_p$-spaces. It has also been identified for finite dimensional polyhedral normed spaces, see \cite{CKS,HSWW,JS,KMN}. In that case, the horofunction compactification is homeomorphic to closed unit ball of the dual space \cite{JS1} and closely related to projective toric varieties \cite{JS}.  For arbitrary (possibly infinite dimensional) normed spaces the Busemann points have been characterised by Walsh \cite{Wa4}. For Hilbert geometries there exists a characterisation of the Busemann points \cite{Wa3}, and for the Hilbert distance on a symmetric cone in a Euclidean Jordan algebra, the horofunction compactification was obtained in \cite{LLNW}, for the cone in a (possibly infinite dimensional) spin factor  in \cite{Cl}, and for the $p$-metrics, with $1\leq p<\infty$, on the symmetric cone in $\mathrm{Herm}_n(\mathbb{C})$ in \cite{FF}.

\section{Metric geometry preliminaries} We start by recalling the construction of the horofunction compactification and the detour distance. 

Let $(M,d)$ be a metric space and let $\mathbb{R}^M$ be the space of all real functions on $M$ equipped with the  topology of pointwise convergence. Fix  a $b\in M$, which is called the {\em basepoint}, and let $\mathrm{Lip}^c_b(M)$ denote the set of all functions $h\in\mathbb{R}^M$ such that $h(b)=0$ and $h$ is $c$-Lipschitz, i.e., $|h(x)-h(y)|\leq cd(x,y)$ for all $x,y\in M$. 

Then $\mathrm{Lip}^c_b(M)$ is a compact subset of $\mathbb{R}^M$. Indeed, the complement of $\mathrm{Lip}^c_b(M)$ is open, so $\mathrm{Lip}^c_b(M)$ is closed subset of $\mathbb{R}^M$. Moreover, as $|h(x)|= |h(x)-h(b)|\leq cd(x,b)$  for all $h\in \mathrm{Lip}^c_b(M)$ and $x\in M$, we get that $\mathrm{Lip}_b^c(M)\subseteq [-cd(x,b),cd(x,b)]^M$, which is compact by Tychonoff's theorem. 

For $y\in M$ define the real valued function, 
\begin{equation}\label{internalpoint}
h_{y}(z) = d(z,y)-d(b,y)\mbox{\quad with $z\in M$.}
\end{equation}
Then $h_y(b)=0$ and $|h_y(z)-h_y(w)| = |d(z,y)-d(w,y)|\leq d(z,w)$. Thus, $h_y\in  \mathrm{Lip}_b^1(M)$ for all $y\in M$.  Using the previous observation one  now defines the {\em horofunction compactification} of  $(M,d)$ to be the closure of $\{h_y\colon y\in M\}$ in $\mathbb{R}^M$, which is a compact subset of $ \mathrm{Lip}_b^1(M)$ and is denoted by $\overline{M}^h$. Its elements are called {\em metric functionals}, and the boundary $\partial \overline{M}^h= \overline{M}^h\setminus \{h_y\colon y\in M\}$ is called the {\em horofunction boundary}. The metric functionals in $\partial \overline{M}^h$ are called {\em horofunctions}, and all other metric functionals are said to be {\em internal points}. 

The topology of pointwise convergence on $ \mathrm{Lip}_b^1(M)$ coincides with the topology of uniform convergence on compact sets, see \cite[Section 46]{Mun}.  In general the topology of pointwise convergence on  $ \mathrm{Lip}_b^1(M)$ is not metrizable, and hence  horofunctions are limits of nets rather than sequences. If, however, the metric space is separable, then the pointwise convergence topology on $\mathrm{Lip}_b^1(M)$ is metrizable and each horofunction is the limit of a sequence.   It should be noted that the embedding $\iota\colon M \to   \mathrm{Lip}_b^1(M)$, where $\iota(y) = h_y$, may not have a continuous inverse on $\iota(M)$, and hence the metric compactification is not always a compactification in the strict topological sense.  If, however, $(M,d)$ is proper (i.e. closed balls are compact) and geodesic, then $\iota$ is a homeomorphism from $M$ onto $\iota(M)$.  Recall that a map $\gamma$ from a (possibly unbounded)  interval $I\subseteq \mathbb{R}$ into a metric space $(M,d)$ is called a {\em geodesic path} if 
\[
d(\gamma(s),\gamma(t)) = |s-t|\mbox{\quad for all }s,t\in I.
\]
The image, $\gamma(I)$, is called a {\em geodesic}, and a metric space $(M,d)$ is said to be {\em geodesic} if for each $x,y\in M$ there exists a geodesic path $\gamma\colon [a,b]\to M$ connecting $x$ and $y$, i.e, $\gamma(a)=x$ and $\gamma(b) = y$.  We call a geodesic $\gamma([0,\infty))$ a {\em geodesic ray}. 

The following fact, which is slightly weaker than \cite[Theorem 4.7]{Ri}, will be useful in the sequel. 
\begin{lemma}\label{Rieffel} 
If $(M,d)$ is a proper geodesic metric space, then $h\in  \partial \overline{M}^h$ if and only if there exists a sequence $(x^n)$ in $M$ with $d(b,x^n)\to \infty$ such that $(h_{x^n})$ converges to $h\in \overline{M}^h$ as $n\to \infty$.
\end{lemma}

A net $(x^\alpha)$  in $(M,d)$ is called an {\em almost geodesic net} if there exists $w\in M$ and  for all $\epsilon>0$ there exists a $\beta$ such that 
\[
d(x^\alpha,x^{\alpha'}) +d(x^{\alpha'},w) - d(x^\alpha,w) <\epsilon\mbox{\quad for all }\alpha\geq \alpha'\geq \beta.
\]
The notion of an almost geodesic sequence  goes back to  Rieffel \cite{Ri} and was further developed by Walsh and co-workers in \cite{AGW,L1,LW,Wa4}. 
In particular, every unbounded almost geodesic net yields a horofunction for a complete metric space \cite{Wa4}. 
\begin{lemma} Let $(M,d)$ be a complete metric space. If $(x^\alpha)$ is an unbounded almost geodesic net in $M$, then 
\[
h(z) = \lim_{\alpha} d(z,x^\alpha)-d(b,x^\alpha)
\]
exists for all $z\in M$ and $h\in\partial \overline{M}^h$.
\end{lemma}
Given a  complete metric space $(M,d)$, a horofunction $h\in\overline{M}^h$ is called a {\em Busemann point} if there exists an almost geodesic  net $(x^\alpha)$ in $M$ such that $h(z) = \lim_{\alpha} d(z,x^\alpha)- d(b,x^\alpha)$ for all $z\in M$. We denote the collection of all Busemann points by $\mathcal{B}_M$. 

Suppose that $(M,d)$ is a complete metric space and $h,h'\in \partial \overline{M}^h$ be horofunctions. Let $W_h$ be the collection of neighbourhoods of $h$ in $\overline{M}^h$. The {\em detour cost} is given by  
\[
H(h,h') = \sup_{W\in W_h}\left(\inf_{x\colon \iota(x)\in W} d(b,x)  +h'(x)\right).
\]
The {\em detour distance} is given by 
\[
\delta(h,h') = H(h,h')+H(h',h).
\]

It is known \cite{Wa4} that  if $(x^\alpha)$ is an almost geodesic net converging to a horofunction $h$, then 
\begin{equation}\label{detourcost}
H(h,h') = \lim_\alpha d(b,x^\alpha) +h'(x^\alpha).
\end{equation}
for all horofunctions $h'$. Moreover, on the set of Busemann points $\mathcal{B}_M$ the detour distance  is a metric where points can be at infinite distance from each other, see \cite{Wa4}.   The detour distance yields a partition of $\mathcal{B}_M$ into equivalence classes,  called {\em parts}, where $h$ and $h'$ are equivalent if  $\delta(h,h') <\infty$. The equivalence class of $h$ is denoted by $\mathcal{P}_h$. So $(\mathcal{P}_h,\delta)$ is a metric space and $\mathcal{B}_M$ is the disjoint union of metric spaces under the detour distance. 
Unlike in the setting of $\mathrm{CAT}(0)$ spaces, where each part is a singleton,  the parts  in the spaces under consideration in this paper are nontrivial.

\section{Complex manifolds}
In this section we investigate Problem \ref{keyproblem} for certain bounded convex domains in $\mathbb{C}^n$ with the Kobayashi distance. We will start by recalling some basic concepts. 
\subsection{Product domains and Kobayashi distance} 
On a convex domain $D\subseteq \mathbb{C}^n$ the {\em Kobayashi distance} is given by  
\[
k_D(z,w) =\inf\{ \rho(\zeta,\eta)\colon \mbox{ $\exists f\colon \Delta \to D$ holomorphic with $f(\zeta)=z$ and $f(\eta)=w$}\}.
\]
for all $z,w\in D$., where 
\[
\rho(z,w) = \log \frac{ 1+\left|\frac{w-z}{1-\bar{z}w}\right| }{1-\left|\frac{w-z}{1-\bar{z}w}\right| }=2\tanh^{-1}\left( 1 -\frac{(1-|w|^2)(1-|z|^2)}{|1-w\bar{z}|^2}\right)^{1/2}\]
is  the {\em hyperbolic distance}  on the open disc, $\Delta:=\{z\in\mathbb{C}\colon |z|<1\}$. 

It is known, see \cite[Proposition 2.3.10]{Ab},  that if $D\subset \mathbb{C}^n$ is bounded convex domain, then $(D,k_D)$ is a proper metric space, whose topology coincides with the usual topology on $\mathbb{C}^n$. Moreover, $(D,k_D)$ is a geodesic metric space containing geodesics rays, see \cite[Theorem 2.6.19]{Ab} or \cite[Theorem 4.8.6]{Ko}.

For the Euclidean ball $B^n =\{(z_1,\ldots,z_n)\in\mathbb{C}^n\colon \|z\|^2<1\}$, where $\|z\|^2 = \sum_i |z_i|^2$, the Kobayashi distance satisfies
\[
k_{B^n}(z,w) = 2\tanh^{-1}\left( 1 -\frac{(1-\|w\|^2)(1-\|z\|^2)}{|1-\langle z,w\rangle |^2}\right)^{1/2}
\]
for all $z,w\in B^n$, see \cite[Chapters 2.2 and 2.3]{Ab}. 

In our setting we will consider product domains $B=\prod_{i-1}^r B_i$, where each $B_i$ is a open unit ball of a norm in $\mathbb{C}^{n_i}$, and we will  use the product property of $k_B$, which says that 
\[k_B(z,w) = \max_{i=1,\ldots,r} k_i(z_i,w_i),\] 
where $k_i$ is the Kobayashi distance on $B_i$, see \cite[Theorem 3.1.9]{Ko}.  So for the polydisc $\Delta^r=\{(z_1,\ldots,z_r)\in\mathbb{C}^r\colon \max_i |z_i|<1\}$, the Kobayashi distance satisfies
\[
k_{\Delta^r}(z,w) =\max_i \rho(z_i,w_i)\mbox{\quad for all $w=(w_1,\ldots,w_r), z=(z_1,\ldots,z_r)\in\Delta^r$.}
\]

For the Euclidean ball, $B^n$, it is well known that the  horofunctions of $(B^n,k_{B^n})$, with basepoint $b=0$, are given by 
\begin{equation}\label{hor}
h_\xi(z) =  \log\frac{ |1 -\langle z,\xi\rangle |^2}{1-\|z\|^2}\mbox{\quad for all $z\in B^n$,}
\end{equation}
 where $\xi\in\partial B^n$.  Moreover,  each horofunction $h_\xi$ is a Busemann point, as it is the limit induced by the geodesic ray $t\mapsto \frac{e^t-1}{e^t+1}\xi$, for $0\leq t<\infty$. 
 
Moreover, if $B$ is a product of Euclidean balls, then the horofunctions are known, see  \cite[Proposition 2.4.12]{Ab} and \cite[Corollary 3.2]{L1}. Indeed, for a product of Euclidean balls $B^{n_1}\times\cdots\times B^{n_r}$ the Kobayashi distance horofunctions  with basepoint $b=0$ are precisely the functions of the form, 
\[
h(z) = \max_{j\in J} \left( h_{\xi_j}(z_j)-\alpha_j\right),
\]
where $J\subseteq \{1,\ldots,r\}$  nonempty,  $\xi_j\in\partial B^{n_j}$ for $j\in J$, and $\min_{j\in J} \alpha_{j}=0$.  
Moreover, each horofunction is a Busemann point. 

The form of the horofunctions of the product of Euclidean balls is essentially due to the product property of the Kobayashi distance and the smoothness and convexity properties of the balls.  Indeed, more generally, the following result holds, see \cite[Section 2 and Lemma 3.3]{L1}.

\begin{theorem}\label{thm:cnhoro} If $D_i\subset \mathbb{C}^{n_i}$ is a bounded strongly convex domain with $C^3$-boundary,  then for each $\xi_i\in\partial D_i$ there exists a unique horofunction $h_{\xi_i}$ which is the limit of a geodesic $\gamma$ from the basepoint $b_i\in D_i$ to $\xi_i$. Moreover, these are all horofunctions. 
If $D = \prod_{i=1}^r D_i$, where each $D_i$  is a bounded strongly convex domain with $C^3$-boundary, then each horofunction $h$ of $(D,k_D)$ (with respect to the basepoint $b=(b_1,\ldots,b_r)$ is of the form, 
\begin{equation}\label{eq:4.1}
h(z) = \max_{j\in J} \left( h_{\xi_j}(z_j)-\alpha_j\right),
\end{equation}
where $J\subseteq \{1,\ldots,r\}$  nonempty,  $\xi_j\in\partial D_j$ for $j\in J$, and $\min_{j\in J} \alpha_{j}=0$. 
Furthermore, each horofunction is a Busemann point, and  the part of $h$ consists of those horofunctions $h'$ with
\[
h'(z) = \max_{j\in J} \left( h_{\xi_j}(z_j)-\beta_j\right),
\]
with $\min_{j\in J}\beta_j =0$. 
\end{theorem}

Now let $D = \prod_{i=1}^r D_i$, where each $D_i$  is a bounded strongly convex domain with $C^3$-boundary. Given  $J\subseteq \{1,\ldots,r\}$  nonempty,  $\xi_j\in\partial D_j$ for $j\in J$, and $\alpha_j\geq 0$ for $j\in J$ with $\min_{j\in J} \alpha_{j}=0$, we can find geodesic paths $\gamma_j\colon [0,\infty)\to D_j$ from $b_j$ to $\xi_j$, and form the path $\gamma\colon [0,\infty)\to D$, where 
\begin{equation}\label{gamma} 
\gamma(t)_j = \left[\begin{array}{ll} \gamma_j(t-\alpha_j) &  \mbox{ for all $j\in J$ and $t\geq \alpha_j$}\\ b_j &\mbox{ otherwise.}\end{array}\right.
\end{equation}  
\begin{lemma}\label{lem:geo} The path $\gamma\colon [0,\infty)\to D$  in (\ref{gamma}) is a geodesic path, and $h_{\gamma(t)}\to h$ where $h$ is given by (\ref{eq:4.1}).
\end{lemma}
\begin{proof} Let $k_i$ denote the Kobayashi distance on $D_i$. 
By the product property we have that 
\[
k_D(\gamma(s),\gamma(t)) =\max_i k_i(\gamma(s)_i,\gamma(t)_i)
\]
for all $s\geq t\geq 0$. By construction $k_i(\gamma(s)_i,\gamma(t)_i)\leq k_i(\gamma_i(s),\gamma_i(t))=s-t$ for all $i$ and $s\geq t\geq 0$. For $j\in J$ with $\alpha_j=0$ we have that $k_j(\gamma(s)_j,\gamma(t)_j) = k_j(\gamma_j(s),\gamma_j(t)) =s-t$ for all $ s\geq t\geq 0$, and hence 
\[
k_D(\gamma(s),\gamma(t)) =\max_i k_i(\gamma(s)_i,\gamma(t)_i) = s-t
\]
for all $s\geq t\geq 0$.

Note that for $z\in D$ we have
\begin{eqnarray*}
\lim_{t\to\infty}h_{\gamma(t)}(z) &= &\lim_{t\to\infty} k_D(z,\gamma(t)) - k_D(\gamma(t),b) \\
  & = & \lim_{t\to\infty} \max_i (k_i(z_i,\gamma(t)_i) - t)\\
  & = & \lim_{t\to\infty} \max_{j\in J} (k_j(z_j,\gamma(t)_j) - t)\\
  & = & \lim_{t\to\infty} \max_{j\in J} (k_j(z_j,\gamma_j(t-\alpha_j)) - k_j(\gamma_j(t-\alpha_j),b_j)-\alpha_j)\\
  & = & \max_{j\in J} \left( h_{\xi_j}(z_j)-\alpha_j\right),
\end{eqnarray*}
which shows that $h_{\gamma(t)}\to h$. 
\end{proof}

Consider $B=\prod_{i=1}^r B_i\subseteq \mathbb{C}^n$, where each $B_i$ is an open unit ball of a norm in $\mathbb{C}^{n_i}$. Then $B$ is the open unit ball of the norm $\|\cdot\|_B$ on $\mathbb{C}^n$. In fact, 
\[
\|w\|_B =\max_{i=1,\ldots,r} \|w_i\|_{B_i}.
\]
Its dual norm satisfies $\|z\|_B^* = \sum_{i=1}^r \|z_i\|_{B_i}^*$ 
and has closed unit ball, 
\[
B_1^* =\{ z\in\mathbb{C}^n\colon \mathrm{Re} \langle w,z\rangle \leq 1\mbox{ for all }w\in\mathrm{cl}{B} \}.  
\]

Now suppose that each $B_i$ is strictly convex and smooth. Then the closed ball $B_1^*$ has extreme points $p(\xi^*_i)= (0,\ldots,0,\xi_i^*,0,\ldots,0)$, where $\xi_i^*\in\mathbb{C}^{n_i}$ is the unique supporting functional at $\xi_i\in\partial B_i$, i.e.,  $\mathrm{Re} \langle \xi_i,\xi_i^*\rangle = 1$ and $\mathrm{Re} \langle w_i,\xi_i^*\rangle < 1$ for $w_i\in\mathrm{cl}\, B_i$ with $w_i\neq \xi_i$. 

The relatively open faces of $B_1^*$ are the sets of the form:
\[
F(\{\xi_j\in\partial B_j\colon j\in J\}) = \left\{\sum_{j\in J}\lambda_jp(\xi_j^*)\colon \sum_{j\in J}\lambda_j =1\mbox{ and }\lambda_j>0\mbox{ for all }j\in J\right\},
\]
where $J\subseteq \{1,\ldots,r\}$ is nonempty and $\xi_j\in\partial B_j$ are fixed. 

On $B$ the Kobayashi distance has a Finsler structure in terms of the infinitesimal Kobayashi metric, see e.g., \cite[Chapter 2.3]{Ab}. Indeed, we have that  
\[
k_B(z,w) =\inf_\gamma L(\gamma),
\]
where the infimum is taken over all piecewise $C^1$-smooth paths $\gamma\colon [0,1]\to B$ with $\gamma(0)=z$ and $\gamma(1) =w$, and 
\[
L(\gamma) = \int_0^1 \kappa_B(\gamma(t),\gamma'(t))\mathrm{d}t,
\]
with
\[
\kappa_B(u,v) = \inf\{|\xi|\colon \exists \phi\in\mathrm{Hol}(\Delta,B)\mbox{ such that }\phi(0) = u\mbox{ and }(\mathrm{D}\phi)_0(\xi) = v\}. 
\]
\begin{proposition}{\cite[Proposition 2.3.24]{Ab}} If $B$ is the open unit ball of a norm on $\mathbb{C}^n$, then 
\[
\kappa_B(0,v) = \|v\|_B\mbox{\quad for all }v\in\mathbb{C}^n.
\]
\end{proposition}

For $z\in B$ and $i=1,\ldots,r$, if $z_i\neq 0$, then we let $z_i' = \|z_i\|_{B_i}^{-1}z_i\in\partial B_i$ and we write $p(z_i^*) = (0,\ldots,0,z_i^*,0,\ldots,0)$, where $z_i^*$ is the unique supporting functional at $z_i'\in\partial B_i$. If $z_i=0$, we set $p(z_i^*)=0$. 

We will now define a map $\phi_B\colon\overline{B}^h\to B_1^*$ and show in the remainder of this section that it is a homeomorphism.  For $z\in B$ let 
\[
\phi_B(z) = \frac{1}{\sum_{i=1}^r e^{k_i(z_i,0)} + e^{-k_i(z_i,0)}} \left(\sum_{i=1}^r (e^{k_i(z_i,0)} - e^{-k_i(z_i,0)})p(z^*_i)\right).
\] 
For a horofunction $h$ given by (\ref{eq:4.1}) we define 
\[
\phi_B(h) = \frac{1}{\sum_{j\in J} e^{-\alpha_j}}\left( \sum_{j\in J} e^{-\alpha_j}p(\xi_j^*)\right).
\]

In fact, we will prove the following theorem.
\begin{theorem}\label{thm:4.3} 
If $B = \prod_{i=1}^r B_i$, where each $B_i$  is the open unit ball of a norm on $\mathbb{C}^{n_i}$ which is strongly convex and has a $C^3$-boundary, then $\phi_B\colon \overline{B}^h\to B_1^*$ is a homeomorphism, which maps each part of $ \overline{B}^h$ onto the relative interior of a boundary face of $B_1^*$. 
\end{theorem}

In view of this result the following version of Problem \ref{keyproblem} is of interest. 
\begin{problem} Suppose that $B$ is the open unit ball  of a norm on $\mathbb{C}^n$ and equipped with the Kobayashi distance. For which $B$ does there exist a homeomorphism from $\overline{B}^h$ onto  $B_1^*$ which maps each part of $ \overline{B}^h$ onto the relative interior of a boundary face of $B_1^*$?  
Of particular interest are bounded symmetric domains realised as the open unit ball in a $\mathrm{JB}^*$-triple, see \cite{Kaup}. 
\end{problem}

\subsection{The map $\phi_B$: injectivity and surjectivity} 
Throughout the remainder of this section we assume that $B = \prod_{i=1}^r B_i$ and each $B_i$  is the open unit ball of a norm on $\mathbb{C}^{n_i}$, which is strongly convex and has a $C^3$-boundary. So for each $\xi_i\in\partial B_i$ there exists a unique $\xi_i^*\in \mathbb{C}^{n_i}$ such that 
\[
\mathrm{Re}\langle \xi_i,\xi_i^*\rangle = 1\mbox{ and } \mathrm{Re}\langle w_i, \xi_i^*\rangle < 1\mbox{ for all }w_i\in \mathrm{cl}\,B_i\mbox{ with }w_i\neq \xi_i,
\]
as $\mathrm{cl}\,B_i$ is strictly convex and smooth.

We start with the following basic observation. 
\begin{lemma}\label{lem:4.4} 
For each $z\in B$ we have that $\phi_B(z) \in \mathrm{int}\, B_1^*$, and  $\phi_B(h) \in\partial B_1^*$ for all $h\in \partial\overline{B}^h$. 
\end{lemma}
\begin{proof}
Note that for $z\in B$ and $w\in\mathrm{cl}\,B$ we have that 
\begin{eqnarray*}
\mathrm{Re} \langle w,\phi_B(z)\rangle &= &\frac{1}{\sum_{i=1}^r e^{k_j(z_i,0)} + e^{-k_i(z_i,0)}} \left(\sum_{i=1}^r (e^{k_i(z_i,0)} - e^{-k_i(z_i,0)})\mathrm{Re}\langle w_i,z^*_i\rangle\right)\\
 &\leq &\frac{1}{\sum_{i=1}^r e^{k_i(z_i,0)} + e^{-k_i(z_i,0)}} \left(\sum_{i=1}^r e^{k_i(z_i,0)} - e^{-k_i(z_i,0)}\right)\\
 &< &1-\delta
\end{eqnarray*}
for some $0<\delta <1$, which is independent of $w$. Thus, $\sup_{w\in \mathrm{cl} B} \mathrm{Re} \langle w,\phi_B(z)\rangle<1-\delta<1$, and hence $\phi_B(z) \in\mathrm{int}\, B_1^*$. 

To see that $\phi_B(h)\in\partial B_1^*$, note that for $w =\sum_{j\in J}p(\xi_j)\in \mathrm{cl}\, B$, where $p(\xi_j) = (0,\ldots,0,\xi_j,0,\ldots,0)$,  we have that 
$\mathrm{Re}\langle w,\phi_B(h)\rangle = 1$.
\end{proof}
To  show that $\phi_B$ is injective on $B$, we need the following basic calculus fact, which can also be found in \cite{JS1}. For completeness we include the proof. 
\begin{lemma}\label{lem:calc} If $\mu\colon\mathbb{R}^r\to \mathbb{R}$ is given by $\mu(x_1,\ldots,x_r) = \sum_{i=1}^r e^{x_i} +e^{-x_i}$, then 
$x\mapsto \nabla \log \mu(x) $ is injective on $\mathbb{R}^r$.
\end{lemma} 
\begin{proof}
For $0<t<1$ we let $ p=1/t\geq 1$ and $q = 1/(1-t)\geq 1$. Then  by H\"older's  inequality we have that 
\begin{eqnarray*} 
\mu( tx +(1-t)y) & = & \sum_{i=1}^r e^{tx_i}e^{(1-t)y_i} + \sum_{i=1}^r e^{-tx_i}e^{-(1-t)y_i}\\
  & \leq & \left(\sum_{i=1}^r (e^{tx_i})^p + \sum_{i=1}^r (e^{-tx_i})^p\right)^{1/p}\left(\sum_{i=1}^r (e^{(1-t)y_i})^q + \sum_{i=1}^r (e^{-(1-t)y_i})^q\right)^{1/q}\\
  & \leq & \left(\sum_{i=1}^r e^{x_i}+ e^{-x_i}\right)^{t}\left(\sum_{i=1}^r e^{y_i} + e^{-y_i}\right)^{1-t},
\end{eqnarray*}
which implies that $\mu(tx+(1-t)y)\leq \mu(x)^t\mu(y)^{1-t}$. Moreover, equality holds if and only if 
\[ 
e^{\pm x_i} = (e^{\pm tx_i})^p = C(e^{\pm(1-t)y_i})^q = Ce^{\pm y_i}
\] 
for all $i$ and some fixed $C>0$. This is equivalent to $ \pm x_i =\pm y_i +\log C$ for all $i$, and hence we have equality if and only if $x=y$. 

Thus, $x\mapsto  \log \mu(x) $  is a strictly convex function on $\mathbb{R}^r$.  By strict convexity we have that 
\[
  \log \mu(x)  - \log \mu(y) > \frac{\log\mu(y+1/2(x-y)) -\log\mu(y)}{1/2} >\frac{\log\mu(y+1/4(x-y)) -\log\mu(y)}{1/4}>\ldots. 
\]
so that $ \log \mu(x)  - \log \mu(y) >\nabla\log\mu(y)\cdot(x-y)$. 
Likewise, $ \log \mu(y)  - \log \mu(x)> \nabla\log \mu(x)\cdot (y-x)$. Combining the inequalities, we see that $0>(\nabla\log \mu(y)-\nabla\log \mu(x))\cdot (x-y)$ for all $x\neq y$, and hence $x\mapsto \nabla \log \mu(x) $ is injective on $\mathbb{R}^r$.
\end{proof}
Note that 
\[
(\nabla\log \mu(x))_j = \frac{e^{x_j}-e^{-x_j}}{\sum_{i=1}^r e^{x_i}+e^{-x_i}}\mbox{\quad  for all } j.
\] 

\begin{lemma}\label{lem:4.5} The map $\phi_B$ is a continuous bijection from $B$ onto $\mathrm{int}\, B_1^*$.
\end{lemma}
\begin{proof}
Cleary $\phi_B$ is continuous on $B$ and $\phi_B(z) =0$ if and only if $z=0$.  Suppose that $z,w\in B\setminus\{0\}$ are such that $\phi_B(z) = \phi_B(w)$.  For simplicity write
\[
\alpha_j = \frac{e^{k_j(z_j,0)} - e^{-k_j(z_j,0)}}{\sum_{i=1}^r e^{k_i(z_i,0)} + e^{-k_i(z_i,0)}}\geq 0
\mbox{\quad  and \quad }
\beta_j = \frac{e^{k_j(w_j,0)} - e^{-k_j(w_j,0)}}{\sum_{i=1}^r e^{k_i(w_i,0)} + e^{-k_i(w_i,0)}}\geq 0.
\]
Note that $\alpha_jp(z_j^*)=0$ if and only if $z_j=0$, and, $\beta_jp(w_j^*)=0$ if and only if $w_j=0$.  Thus, $z_j=0$ if and only if $w_j=0$. 
Now suppose that $z_j\neq 0$, so $w_j\neq 0$. Then $\langle p(v_j),\phi_B(z) \rangle = \langle p(v_j),\phi_B(w) \rangle$ for each $v_j\in B_j$. This implies that
\[
\alpha_j\langle v_j,z^*_j \rangle = \beta_j\langle v_j,w^*_j \rangle\mbox{\quad for all }v_j\in B_j,
\]
and hence $\alpha_j z^*_j = \beta_j w^*_j$. It follows that $\alpha_j =\beta_j$ and $z^*_j=w^*_j$.  Thus $z_j = \mu_j w_j$ for some $\mu_j>0$. As $\alpha_i=\beta_i$ for all $i\in\{1,\ldots,r\}$, we know by  Lemma \ref{lem:calc} that $k_j(z_j,0) = k_j(w_j,0)$, and  hence  $z_j =w_j$ by \cite[Proposition 2.3.5]{Ab}.  So $z=w$, which shows that $\phi_B$ is injective.

As $\phi_B$ is injective and continuous on $B$, it follows from Brouwer's  domain invariance theorem that $\phi_B(B)$ is an open subset of $\mathrm{int}\,B^*_1$ by Lemma \ref{lem:4.4}. 
Suppose, by way of contradiction, that $\phi_B(B)\neq \mathrm{int}\, B^*_1$.  Then $ \partial \phi_B(B)\cap \mathrm{int}\, B_1^*$ is nonempty, as otherwise $\phi_B(B)$ is closed and open in $\mathrm{int}\, B^*_1$, which would imply that $\mathrm{int}\,B_1^*$ is the disjoint union of the nonempty open sets $\phi_B(B)$ and its complement contradicting the connectedness of $\mathrm{int}\, B^*_1$. So let $w\in\partial \phi_B(B)\cap \mathrm{int}\, B_1^*$ and  $(z^n)$ be a sequence in $B$ such that  $\phi_B(z^n)\to w$.  As $\phi_B$ is continuous on $B$, we have that $k_B(z^n,0)\to\infty$. 

Using the product property, $k_B(z^n,0) = \max_i k_i(z^n_i,0)$, we may assume after taking subsequences  that $\alpha_i^n =  k_B(z^n,0) - k_{i}(z^n_i,0)\to \alpha_i\in [0,\infty]$ and $z^n_i\to \zeta_i\in\mathrm{cl}\, B_i$ for all 
$i$. Let $I = \{i\colon \alpha_i<\infty\}$, and note that for each $i\in I$, $\zeta_i\in\partial B_i$, as $k_i(z_i^n,0)\to\infty$. Then 
\begin{eqnarray*}
 \phi_B(z^n) & = & \frac{1}{\sum_{i=1}^r e^{k_i(z^n_i,0)} + e^{-k_i(z^n_i,0)}} \left(\sum_{i=1}^r (e^{k_i(z^n_i,0)} - e^{-k_i(z^n_i,0)})p((z_i^n)^*)\right)\\
       & = & \frac{1}{\sum_{i=1}^r e^{-\alpha^n_i} + e^{-k_B(z^n,0) -k_i(z^n_i,0)}} \left(\sum_{i=1}^r (e^{-\alpha^n_i} - e^{-k_B(z^n,0)-k_i(z^n_i,0)})p((z^n_i)^*)\right).
\end{eqnarray*} 
Letting $n\to\infty$, the righthand side converges to 
\[
\frac{1}{\sum_{i\in I} e^{-\alpha_i}} \left(\sum_{i\in I} e^{-\alpha_i}p(\zeta_i^*)\right) =w.
\]
But this implies that $w\in\partial B^*_1$, as $\mathrm{Re}\langle \sum_{i\in I}p(\zeta_i), w\rangle = 1$ and $\sum_{i\in I} p(\zeta_i)\in \mathrm{cl}\, B$, where $p(\zeta_i) = (0,\ldots,0,\zeta_i,0,\ldots,0)$. This is impossible and hence $\phi_B(B)=\mathrm{int}\, B^*_1$. 
\end{proof}
We now analyse $\phi_B$ on $\partial \overline{B}^h$. 
\begin{lemma}\label{lem:4.6} The map $\phi_B$ maps $\partial \overline{B}^h$ bijectively onto $\partial B_1^*$. Moreover, the part $\mathcal{P}_h$, where $h$ is given by (\ref{eq:4.1}), is mapped onto the relative open boundary face 
\[
F(\{\xi_j\in\partial B_j\colon j\in J\}) = \left\{\sum_{j\in J}\lambda_jp(\xi_j^*)\colon \sum_{j\in J}\lambda_j =1\mbox{ and }\lambda_j>0\mbox{ for all }j\in J\right\}.
\] 
\end{lemma}
\begin{proof}
We know from Lemma \ref{lem:4.4} that $\phi_B$ maps $\partial \overline{B}^h$ into $\partial B_1^*$. To show that it is onto we let $w\in\partial B_1^*$. As $B^*_1$ is the disjoint union of its relative open faces (see \cite[Theorem 18.2]{Rock}), there exist $J\subseteq \{1,\ldots,r\}$, extreme points $p(\xi^*_j)$ of $B_1^*$, and $0<\lambda_j\leq 1$ for $j\in J$ with $\sum_{j\in J}\lambda_j=1$ such that $w=\sum_{j\in J}\lambda_j p(\xi^*_j)$.  Let $\mu_j = -\log \lambda_j$ and $\mu^* = \min_{j\in J} \mu_j$. 
Now set $\alpha_j = \mu_j -\mu^*$ for $j\in J$.  Then $\alpha_j \geq 0$ for $j\in J$ and $\min_{j\in J} \alpha_j =0$. 

Let $h\in\partial \overline{B}^h$ be given by $h(z) = \max_{j\in J}  (h_{\xi_j}(z_j)-\alpha_j)$. 
Then 
\[\phi_B(h) = \frac{\sum_{j\in J} e^{-\alpha_j} p(\xi_j^*)}{\sum_{j\in J} e^{-\alpha_j}} =  \frac{\sum_{j\in J} e^{-\mu_j} p(\xi_j^*)}{\sum_{j\in J} e^{-\mu_j} }
= \frac{\sum_{j\in J} \lambda_j p(\xi^*_j)}{\sum_{j\in J} \lambda_j}=w.
\]

To prove injectivity let $h,h'\in\partial \overline{B}^h$, where $h$ is as in (\ref{eq:4.1}) and 
\begin{equation}\label{eq:h'}
h'(z) = \max_{j\in J'} (h_{\eta_j}(z_j)-\beta_j)
\end{equation}
for $z\in B$. Suppose that $\phi_B(h)=\phi_B(h')$, so 
\[
\phi_B(h) = \frac{\sum_{j\in J}e^{-\alpha_j}p(\xi_j^*)}{\sum_{j\in J} e^{-\alpha_j}}= \frac{\sum_{j\in J'}e^{-\beta_j}p(\eta_j^*)}{\sum_{j\in J'} e^{-\beta_j} } =\phi_B(h').
\]  
We have that $J=J'$. Indeed, if $k\in J$ and $k\not\in J'$, then 
\[
0=\mathrm{Re}\langle p(\xi_k),\phi_B(h')\rangle  = \mathrm{Re}\langle p(\xi_k),\phi_B(h)\rangle >0,
\]
which is impossible. In the other case a contradiction can be derived in the same way. 

Now suppose there exists $k\in J$ such that $\xi_k\neq \eta_k$. If 
\[
\frac{e^{-\alpha_k}}{\sum_{j\in J} e^{-\alpha_j}}\leq \frac{e^{-\beta_k}}{\sum_{j\in J} e^{-\beta_j}},
\]
then 
\[
\mathrm{Re}\langle p(\eta_k) ,\phi_B(h)\rangle  = \frac{e^{-\alpha_k}}{\sum_{j\in J} e^{-\alpha_j}}\mathrm{Re}\langle \eta_k,\xi^*_k\rangle < \frac{e^{-\alpha_k}}{\sum_{j\in J} e^{-\alpha_j}}\leq \frac{e^{-\beta_k}}{\sum_{j\in J} e^{-\beta_j}} = \mathrm{Re}\langle p(\eta_k) ,\phi_B(h')\rangle, 
\]
as $\mathrm{cl}\, B_k$ is smooth and strictly convex. This is impossible. The other case goes in the same way.  Thus, $J=J'$ and $\xi_j=\eta_j$ for all $j\in J$. 

It follows that 
\[
\frac{e^{-\alpha_k}}{\sum_{j\in J} e^{-\alpha_j}} = \mathrm{Re}\langle p(\xi_k) ,\phi_B(h)\rangle  = \mathrm{Re}\langle p(\eta_k) ,\phi_B(h')\rangle = \frac{e^{-\beta_k}}{\sum_{j\in J} e^{-\beta_j}}
\]
for all $k\in J$. We now show that $\alpha_k=\beta_k$ for all $k\in J$ by using ideas similar to the ones used in the proof of Lemma \ref{lem:calc}. 

Let $\nu\colon \mathbb{R}^J\to\mathbb{R}$ be given by $\nu(x) = \sum_{j\in J}e^{-x_j}$. Then for $x,y\in\mathbb{R}^J$ and $0<t<1$ we have that 
\[
\nu(tx+(1-t)y) \leq \nu(x)^t\nu(y)^{1-t},
\]
and we have equality if and only if there exists a constant $c$ such that $x_k = y_k+c$ for all $k\in J$.  So, if $x\neq y +(c,\ldots,c)$ for all $c$, then  $-\nabla \log \nu(x) \neq -\nabla \log \nu (y)$.  

As $\min_{j\in J}\alpha_j =0=\min_{j\in J}\beta_j$, we can conclude that $\alpha_k=\beta_k$ for all $k\in J$. 
This shows that  $h = h'$ and hence $\phi_B$ is injective on $\partial \overline{B}^h$. 

To complete the proof note that  $\phi_B(h)$ is in the relative open boundary face $F(\{\xi_j\in\partial B_j\colon j\in J\})$ of $B_1^*$. Moreover, $h'$ given by  (\ref{eq:h'})  is in the same part as $h$ if, and only if, $J=J'$ and $\xi_j=\eta_j$ for all $j\in J$ by \cite[Propositions 2.8 and 2.9]{L1}.  So, $\phi_B(h')$  lies in $F(\{\xi_j\in\partial B_j\colon j\in J\})$ if and only if $h'$ lies in the same part as $h$.   
\end{proof}
\subsection{Continuity and the proof of Theorem \ref{thm:4.3}}
We now show that $\phi_B$ is continuous on $\overline{B}^h$. 
\begin{proposition}\label{prop:4.7}
The map $\phi_B\colon \overline{B}^h\to B_1^*$ is continuous.
\end{proposition}
\begin{proof}
Clearly $\phi_B$ is continuous on $B$.  Suppose that $(z^n)$ is sequence in $B$ converging to $h\in\partial \overline{B}^h$, where $h$ is given by (\ref{eq:4.1}). To show that $\phi_B(z^n)\to \phi_B(h)$ we show that every subsequence of $(\phi_B(z^n))$ has a subsequence converging to $\phi_B(h)$.  So, let $(\phi_B(z^{n_k}))$ be a  subsequence. Then we can take a further subsequence $(z^{n_{k,m}})$ such that 
\begin{enumerate}[(1)]
\item 
\[
\beta^m_j = k_B(z^{n_{k,m}},0) - k_j(z^{n_{k,m}}_j,0) \to \beta_j\in [0,\infty]\mbox{\quad for all }j\in\{1,\ldots,r\}.
\]
\item There exists $j_0$ such that $\beta^m_{j_0} =0$ for all $m\geq 1$. 
\item $(z^{n_{k,m}}_j)$ converges to $\eta_j\in \mathrm{cl}\, B_j$ and $h_{z^{n_{k,m}}}\to h_{\eta_j}$ for all $j\in\{1,\ldots,r\}$.
\end{enumerate}

Let $J'=\{j\colon \beta_j<\infty\}$. Then $h_{z^{n_{k,m}}}\to h'$, where $h'(z) = \max_{j\in J'} (h_{\eta_j}(z_j) - \beta_j)$ for $z\in B$, as 
\[
\lim_{m\to\infty} k_B(z,z^{n_{k,m}}) -  k_B(z^{n_{k,m}},0) = \lim_{m\to\infty} \max_j (k_j(z_j,z^{n_{k,m}}_j) -  k_j(z^{n_{k,m}}_j,0) -\beta^m_j) 
= \max_{j\in J'} (h_{\eta_j}(z_j) - \beta_j),
\]
by the product property of $k_B$. 

As $h=h'$, we know by \cite[Propositions 2.8 and  2.9]{L1} that $J=J'$,  $\xi_j=\eta_j$ and $\alpha_j =\beta_j$ for all $j\in J$.  We also know by Lemma \ref{Rieffel} that $k_B(z^{n_{k,m}},0)\to\infty$, as $h$ is a horofunction. So, 
\[
\phi_B(z^{n_{k,m}}) = \frac{\sum_{i=1}^r (e^{-\beta_i^m} - e^{-k_B(z^{n_{k,m}},0) -k_i(z^{n_{k,m}}_i,0)})p((z^{n_{k,m}})^*)} {\sum_{i=1}^r e^{-\beta_i^m} - e^{-k_B(z^{n_{k,m}},0) -k_i(z^{n_{k,m}}_i,0)}}\to \frac{\sum_{j\in J} e^{-\beta_j}p(\eta_j^*)}{\sum_{j\in J} e^{-\beta_j}} = \phi_B(h),
\]
which shows that $\phi_B(z^n)\to \phi_B(h)$. 

We know from Lemma \ref{lem:4.4} that $\phi_B(B)\subseteq \mathrm{int}\, B_1^*$ and $\phi_B(\partial \overline{B}^h)\subseteq \partial B_1^*$. So, to complete the proof it remains to show that if $(h_n)$ in $\partial\overline{B}^h$ converges to $h\in\partial\overline{B}^h$, where $h$ is as in (\ref{eq:4.1}), then $\phi_B(h_n)\to \phi_B(h)$. 
For $n\geq 1$ let $h_n$ be given by 
\[
h_n(z) = \max_{j\in J_n} (h_{\eta_j^n}(z_j) - \beta_j^n)
\]
for $z\in B$. We show that every subsequence of $(\phi_B(h_n))$  has a convergent subsequence with limit $\phi_B(h)$. 

So let $(\phi_B(h_{n_k}))$ be a subsequence. Then we can take a further subsequence $(\phi_B(h_{k_m}))$ to get that  \begin{enumerate}[(1)]
\item There exists $J_0\subseteq \{1,\ldots,r\}$ such that $J_{k_m} =J_0$ for all $m$. 
\item There exists $j_0\in J_0$ such that $\beta^{k_m}_{j_0}=0$ for all $m$.
\item $\beta_j^{k_m}\to\beta_j\in [0,\infty]$ for all $j\in J_0$.
\item $\eta_j^{k_m}\to\eta_j$ for all $j\in J_0$.
\end{enumerate} 
Note that for each $j\in J_0$ we have that $h_{\eta_j^{k_m}}\to h_{\eta_j}$ in $\overline{B}_j^h$, as the identity map on $\mathrm{cl} B_j$, that is $\xi_j\in\mathrm{cl} B_j\to h_{\xi_j}\in\overline{B}_j^h$, is a homeomorphism by \cite[Theorem 1.2]{AFGG}.  

Let $J'=\{j\in J_0\colon \beta_j<\infty\}$ and note that $j_0\in J'$. Then for each $z\in B$ we have that 
\[
\lim_{m\to\infty} h_{k_m}(z) = \lim_{m\to\infty} \max_{j\in J_0} (h_{\eta_j^{k_m}}(z_j) -\beta_j^{k_m}) = \lim_{m\to\infty} \max_{j\in J'} (h_{\eta_j^{k_m}}(z_j) -\beta_j^{k_m}) =\max_{j\in J'} (h_{\eta_j}(z_j) -\beta_j).
\]
So, if we let $h'(z) = \max_{j\in J'} (h_{\eta_j}(z_j) -\beta_j)$ for $z\in B$, then $h'$ is a horofunction by Theorem \ref{thm:cnhoro} and $h_{k_m}\to h'$ in $\overline{B}^h$. As $h_n\to h$, we conclude that $h'=h$. This implies that $J'=J$ and $\eta_j=\xi_j$ and $\beta_j=\alpha_j$ for all $j\in J$, as otherwise $\delta(h,h') \neq 0$ by \cite[Proposition 2.9 and Lemma 3.3]{L1}. This implies that $\beta_j^{k_m}\to \alpha_j$ and $\eta_j^{k_m}\to\xi_j$ for all $j\in J'$. Moreover, by definition  $\beta_j^{k_m}\to\infty$ for all $j\in J_0\setminus J'$. Thus, 
\[
\phi_B(h_{k_m}) = \frac{\sum_{j\in J_0} e^{-\beta_j^{k_m}} p((\eta_j^{k_m})^*)
}{\sum_{j\in J_0}e^{-\beta_j^{k_m}}} \to \frac{\sum_{j\in J} e^{-\alpha_j} p(\xi_j^*)
}{\sum_{j\in J}e^{-\alpha_j}} = \phi_B(h),\]
which completes the proof.
\end{proof}
The proof of Theorem \ref{thm:4.3} is now straightforward. 
\begin{proof}[Proof of Theorem \ref{thm:4.3}]
It follows from Lemmas \ref{lem:4.5} and \ref{lem:4.6} and  Proposition \ref{prop:4.7} that $\phi_B\colon \overline{B}^h\to B^*_1$ is a continuous bijection. As  
$ \overline{B}^h$ is compact and $B^*_1$ is Hausdorff, we conclude that $\phi_B$ is a homeomorphism.	
Moreover, $\phi_B$  maps each part of $\partial\overline{B}^h$ onto the relative interior of a boundary face of $B_1^*$ by Lemma \ref{lem:4.6}.
\end{proof}

\section{Finite dimensional normed spaces}
Every finite dimensional normed space $(V,\|\cdot\|)$ has a Finsler structure. Indeed, if we let 
\[
L(\gamma) = \int_0^1 \|\gamma'(t)\|\mathrm{d}t
\]
be the length of a piecewise $C^1$-smooth path $\gamma\colon [0,1]\to V$, then 
\[
\|x-y\| =\inf_\gamma L(\gamma),
\] 
where the infimum is taken over all $C^1$-smooth paths $\gamma\colon [0,1]\to V$ with $\gamma(0)=x$ and $\gamma(1)=y$. 
So for normed spaces $V$ the unit ball in the tangent space $T_bV$ is the same for all $b\in V$. 

We are interested in the  following more explicit version of Problem \ref{keyproblem}, which was posed by Kapovich and Leeb \cite[Question 6.18]{KL}. 
\begin{problem} 
For which finite dimensional normed spaces $(V,\|\cdot\|)$ does there exist a homeomorphism $\phi_V$ from the horofunction compactification of $(V,\|\cdot\|)$ onto the closed dual unit ball $B_1^*$ of $V$, which maps each part of the horofunction boundary onto the relative interior of a boundary face of $B_1^*$?
\end{problem}

We show that such a homeomorphism exists for Euclidean Jordan algebras equipped with the spectral norm. So we will consider the Euclidean Jordan algebras not as inner-product spaces, but as an order-unit space, which makes it a finite dimensional formally real JB-algebra, see \cite[Theorem 1.11]{AS1}.  We will give an explicit description of the horofunctions of these normed spaces and identify the parts and the detour distance.  In our analysis we make frequent use of the theory of Jordan algebras and order-unit spaces. For the reader's convenience we will recall some of the basic 
concepts. Throughout the paper we will follow the terminology used in \cite{AS0,AS1} and \cite{FK}.

\subsection{Preliminaries}
\paragraph{Order-unit spaces} A {\em cone} $V_+$ in a real vector space $V$ is a convex subset of $V$ with $\lambda V_+\subseteq V_+$ for all $\lambda\geq 0$ and $V_+\cap -V_+=\{0\}$. The cone $V_+$ induces a partial ordering $\leq $ on $V$ by $x\leq y$ if $y-x\in V_+$. We write $x<y$ if $x\leq y$ and $x\neq y$.  The cone $V_+$ is said to be {\em Archimedean} if for each $x\in V $ and $y\in V_+$ with $nx\leq y$ for all $n\geq 1$ we have that $x\leq 0$. An element $u$ of $V_+$ is called an {\em order-unit} if for each $x\in V$ there exists $\lambda\geq 0$ such that $-\lambda u\leq x\leq \lambda u$. The triple $(V,V_+,u)$, where $V_+$ is an Archimedean cone and $u$ is an order-unit, is called an {\em order-unit space}. An order-unit space admits a  norm, 
\[
\|x\|_u =\inf\{\lambda\geq 0\colon -\lambda u\leq x\leq \lambda u\},
\]
which is called the {\em order-unit norm}, and we have that $-\|x\|_uu \leq x\leq \|x\|_uu$ for all $x\in V$. The cone $V_+$ is closed under the order-unit norm and $u\in \mathrm{int}\, V_+$. 

A linear functional $\phi$ on an order-unit space is said to be {\em positive} if $\phi(x) \geq 0$ for all $x\in V_+$. It is called a {\em state} if it is positive and $\phi(u) =1$. The set of all states is denoted by $S(V)$ and is called the {\em state space}, which is convex set. In our case, the order-unit space is finite dimensional, and hence $S(V)$ is compact. The extreme points of $S(V)$ are called the {\em pure states}.  

The dual space $V^*$ of an order-unit space $V$ is a {\em base norm space}, see \cite[Theorem 1.19]{AS0}. More specifically, $V^*$ is an ordered normed vector space with cone $V^*_+=\{\phi\in V^*\colon \phi\mbox{ is positive}\}$, $V^*_+-V^*_+ = V^*$, and the unit ball of the norm of $V^*$ is given by 
\[
B_1^* = \mathrm{conv} (S(V)\cup -S(V)).
\]

\paragraph{Jordan algebras} Important examples of order-unit spaces come from Jordan algebras.  A {\em Jordan algebra} (over $\mathbb{R}$) is a real vector space $V$ equipped with a commutative bilinear product $\bullet$ that satisfies the identity
\[
x^2\bullet (y\bullet x) = (x^2\bullet y)\bullet x\mbox{\quad for all }x,y\in V. 
\]
A basic example is the space $\mathrm{Herm}_n(\mathbb{C})$ consisting of $n\times n$ Hermitian matrices with Jordan product $A\bullet B = (AB+BA)/2$. 

Throughout the paper we will assume that $V$ has a {\em unit}, denoted $u$.  For $x\in V$ we let $L_x$ be the linear map on $V$ given by $L_xy =x\bullet y$. A  finite dimensional Jordan algebra is said to be {\em Euclidean} if there exists an inner-product $(\cdot|\cdot)$ on $V$ such that 
\[
(L_xy| z) = (y| L_xz)\mbox{\quad for all }x,y,z\in V.
\]
A Euclidean Jordan algebra has a cone $V_+=\{x^2\colon x\in V\}$. The interior of $V_+$ is a {\em symmetric cone}, i.e., it is self-dual and $\mathrm{Aut}(V_+) =\{A\in\mathrm{GL}(V)\colon A(V_+) = V_+\}$ acts transitively on the interior  of $V_+$. In fact, the Euclidean Jordan algebras are in one-to-one correspondence with the symmetric cones by the Koecher-Vinberg theorem, see for example \cite{FK}.  

The algebraic unit $u$ of a Euclidean Jordan algebra is an order-unit for the cone $V_+$, so the triple $(V,V_+,u)$ is an order-unit space. We will consider the Euclidean Jordan algebras as an order-unit space equipped with the order-unit norm. These are precisely the finite dimensional formally real JB-algebras, see \cite[Theorem 1.11]{AS1}. In the analysis, however, the inner-product structure on $V$ will be exploited to identify $V^*$ with $V$.   
 
Throughout we will fix the rank of the Euclidean Jordan algebra $V$ to be $r$.  In a Euclidean Jordan algebra each $x$ can be written in a unique way as $x= x^+-x^-$, where $x^+$ and $x^-$ are orthogonal element $x^+$ and $x^-$ in $V_+$, see \cite[Proposition 1.28]{AS1}. This is called the {\em orthogonal decomposition of $x$}.  

Given  $x$ in a Euclidean Jordan algebra $V$, the  {\em spectrum} of $x$ is  given by 
$\sigma(x)=\{ \lambda\in\mathbb{R}\colon \lambda u -x\mbox{ is not invertible}\}$, and we have that $V_+ =\{x\in V\colon \sigma(x)\subset [0,\infty)\}$. We write $\Lambda(x) =\inf\{\lambda\colon x\leq\lambda u\}$ and note that $\Lambda(x)=\max\{\lambda\colon \lambda\in\sigma(x)\}$, so that  
 \[
 \|x\|_u= \max\{\Lambda(x),\Lambda(-x)\} = \max\{|\lambda|\colon \lambda \in\sigma(x)\}
 \]
 for all $x\in V$. We also note that 
 \[
 \Lambda(x+\mu u) = \Lambda(x)+\mu
 \]
 for all $x\in V$ and $\mu\in\mathbb{R}$. Moreover, if $x\leq y$, then $\Lambda(x)\leq\Lambda(y)$.

Recall that $p\in V$ is an {\em idempotent} if $p^2=p$. If, in addition, $p$ is non-zero and cannot be written as the sum of two non-zero idempotents, then it is said to be a {\em primitive} idempotent. The set of all primitive idempotent is denoted $\mathcal{J}_1(V)$ and is known to be a compact set \cite{Hi}. Two idempotents $p$ and $q$ are said to be orthogonal if $p\bullet q=0$, which is equivalent to $(p|q)=0$.  According to the spectral theorem \cite[Theorem III.1.2]{FK}, 
each $x$ has a {\em spectral decomposition}, $x = \sum_{i=1}^r \lambda_i p_i$, where each $p_i$ is a primitive idempotent, the $\lambda_i$'s are the eigenvalues of $x$ (including multiplicities), and $p_1,\ldots,p_r$ is a Jordan frame, i.e., the $p_i$'s are mutually orthogonal and $p_1+\cdots+p_r =u$. 

 Throughout the paper we will fix the inner-product on $V$ to be 
 \[
 (x|y) = \mathrm{tr}(x\bullet y), 
 \]
 where $\mathrm{tr}(x) = \sum_{i=1}^r \lambda_i$ and $x = \sum_{i=1}^r \lambda_ip_i$ is the spectral decomposition of $x$. 
 
For $x\in V$ we denote the {\em quadratic representation} by $U_x\colon V\to V$, which is the linear map, 
\[
U_x y = 2 x\bullet (x\bullet y) - x^2\bullet y = 2L_x(L_x y) - L_{x^2} y.
\]
In case of a Euclidean Jordan algebra $U_x$ is self-adjoint, $(U_xy| z) = (y| U_xz)$.

We identify $V$ with $V^*$ using the inner-product. So,  $S(V) =\{w\in V_+\colon (u|w)=1\}$ which is a  compact convex set,  as $V$ is finite dimensional. Moreover, the extreme points of $S(V)$ are the primitive idempotents, see \cite[Proposition IV.3.2]{FK}.  The dual space $(V,\|\cdot\|_u^*)$  is a base norm space with norm,
\[
\|z\|_u^* =\sup\{(x|z)\colon x\in V\mbox{ with } \|x\|_u=1\}.
\]
If $V$ is a Euclidean Jordan algebra, it is known that the (closed) boundary faces of the dual ball $B_1^*= \mathrm{conv} (S(V)\cup -S(V))$ are precisely the sets of the form, 
\begin{equation}\label{ER}
\mathrm{conv}\, ( (U_p(V)\cap S(V))\cup(U_q(V)\cap -S(V))),
\end{equation}
where $p$ and $q$ are orthogonal idempotents, see \cite[Theorem 4.4]{ER}.

\subsection{Summary  of results}
To conveniently describe the horofunction compactification $\overline{V}^h$ of   $(V,\|\cdot\|_u)$, where $V$ is a Euclidean Jordan algebra, we need some additional notation. 
Throughout this section  we will fix the basepoint $b\in V$ to be $0$.  

Let $p_1,\ldots, p_r$ be a Jordan frame in $V$. Given $I\subseteq\{1,\ldots, r\}$ nonempty, we write $p_I = \sum_{i\in I} p_i$ and we let $V(p_I) = U_{p_I}(V)$. For convenience we set $p_\emptyset =0$, so $V(p_\emptyset) =\{0\}$. 

Recall \cite[Theorem IV.1.1]{FK}  that $V(p_I)$ is the Peirce 1-space of the idempotent $p_I$:   
\[
V(p_I)=\{x\in V\colon p_I\bullet x =x\},
\]
which is a subalgebra. Given $z\in V(p_I)$, we write $\Lambda_{V(p_I)}(z)$ to denote the maximal eigenvalue of $z$ in the subalgebra $V(p_I)$. 

The following theorem characterises the horofunctions in $\overline{V}^h$. 
 \begin{theorem}  \label{main1}
 Let $p_1,\ldots,p_r$ be a Jordan frame, $I,J\subseteq\{1,\ldots,r\}$, with $I\cap J=\emptyset$ and $I\cup J$ nonempty, and 
$\alpha\in\mathbb{R}^{I\cup J}$ such that $\min\{\alpha_i \colon i\in I\cup J\}=0$. The function $h\colon V\to \mathbb{R}$ given by, 
\begin{equation}\label{horofunction1}
h(x)  = \max\left\{\Lambda_{V(p_I)}\left(-U_{p_I}x - \sum_{i\in I} \alpha_ip_i\right),  \Lambda_{V(p_J)}\left(U_{p_J}x - \sum_{j\in J} \alpha_j p_j\right)\right\} \mbox{\quad for $x\in V$,} 
\end{equation}
is a horofunction, where we  use the convention that if $I$ or $J$ is empty, the corresponding term is omitted from the maximum. 
Each horofunction in $\overline{V}^h$ is of the form (\ref{horofunction1}) and a Busemann point. 
\end{theorem}
 
To conveniently describe the parts and the detour distance we  introduce the following notation. Given orthogonal idempotents $p_I$ and $p_J$ we let $V(p_I,p_J)=V(p_I)+ V(p_J)$, which is a subalgebra of $V$ with unit $p_{IJ}=p_I+p_J$. The subspace $V(p_I,p_J)$ can be equipped with the {\em variation norm}, 
\[
\|x\|_\mathrm{var} = \Lambda_{V(p_I,p_J)}(x) + \Lambda_{V(p_I,p_J)}(-x)=\mathrm{diam}\, \sigma_{V(p_I,p_J)}(x), 
\]
which is a semi-norm on $V(p_I,p_J)$. The variation norm is, however, a norm on the quotient space $V(p_I,p_J)/\mathbb{R}p_{IJ}$. 

\begin{theorem}\label{main2}
Given horofunctions $h$ and $h'$, where 
\begin{equation}
h(x)  = \max\left\{\Lambda_{V(p_I)}\left(-U_{p_I}x - \sum_{i\in I} \alpha_ip_i\right),  \Lambda_{V(p_J)}\left(U_{p_J}x - \sum_{j\in J} \alpha_j p_j\right)\right\}
\end{equation}
and 
\begin{equation}\label{horofuncton2}
h'(x)  = \max\left\{\Lambda_{V(q_{I'})}\left(-U_{q_{I'}}x - \sum_{i\in I'} \beta_iq_i\right),  \Lambda_{V(q_{J'})}\left(U_{q_{J'}}x - \sum_{j\in J'} \beta_j q_j\right)\right\},
\end{equation}
we have that 
\begin{enumerate}[(i)]
\item $h$ and $h'$ are in the same part if and only if $p_I=q_{I'}$ and  $p_J = q_{J'}$. 
\item If $h$ and $h'$ are in the same part, then $\delta(h,h') = \|a -b\|_{\mathrm{var}}$, where $a=  \sum_{i\in I} \alpha_ip_i + \sum_{j\in J} \alpha_j p_j$ and $b =  \sum_{i\in I'} \beta_iq_i + \sum_{j\in J'} \beta_j q_j$ in $V(p_I,p_J)$. 
\item The part $(\mathcal{P}_h,\delta)$ is isometric to $(V(p_I,p_J)/\mathbb{R}p_{IJ},\|\cdot\|_{\mathrm{var}})$. 
\end{enumerate}
\end{theorem}
\begin{remark}
A basic example is $(\mathbb{R}^n,\|\cdot\|_\infty)$, where $\|z\|_\infty = \max_i |z_i|$, which is an associative Euclidean Jordan algebra. In that case  every horofunction is a Busemann points and of the form,
\[
h(x) = \max\{ \max_{i\in I} ( -x_i -\alpha_i),\max_{j\in J} (x_j -\alpha_i)\}
\]
where $I,J\subseteq \{1,\ldots,n\}$ are disjoint, $I\cup J$ is nonempty and $\alpha\in \mathbb{R}^{I\cup J}$ with $\min_{k\in I\cup J}\alpha_k =0$, see also \cite[Theorem 5.2]{Gu1} and \cite{L1}. Moreover, $(\mathcal{P}_h,\delta)$ is isometric to $(\mathbb{R}^{I\cup J}/\mathbb{R}\mathbf{1},\|\cdot\|_\mathrm{var})$, where $\mathbf{1}=(1,\ldots,1)\in \mathbb{R}^{I\cup J}$. 
\end{remark}

We will show that the following map is a homeomorphism from $\overline{V}^h$ onto $B^*_1$. Let $\phi\colon \overline{V}^h\to B^*_1$  be given by 
\begin{equation}\label{phi1}
\phi(x) = \frac{e^x - e^{-x}}{(e^x+e^{-x}|u)} = \frac{1}{\sum_{i=1}^r e^{\lambda_i} +e^{-\lambda_i}}\left(\sum_{i=1}^r (e^{\lambda_i} -e^{-\lambda_i})p_i\right)
\end{equation}
 for $x = \sum_{i=1}^r \lambda_ip_i\in V$, and 
\begin{equation}\label{phi2}
\phi(h) = \frac{1}{\sum_{i\in I} e^{-\alpha_i} +\sum_{j\in J} e^{-\alpha_j}}\left(\sum_{i\in I} e^{-\alpha_i}p_i -\sum_{j\in J} e^{-\alpha_j}p_j\right)
\end{equation}
for $h\in \partial \overline{V}^h$ given by (\ref{horofunction1}). 

We should note that $\phi$ is well defined. To verify this assume that the horofunction $h$ given by (\ref{horofunction1}) is represented as 
\[
h(x)  = \max\left\{\Lambda_{V(q_{I'})}\left(-U_{q_{I'}}x - \sum_{i\in I'} \beta_iq_i\right),  \Lambda_{V(q_{J'})}\left(U_{q_{J'}}x - \sum_{j\in J'} \beta_j q_j\right)\right\}
\]
for $x\in V$. Then it follows from Theorem \ref{main2} that $p_I=q_{I'}$ and  $p_J = q_{J'}$. Moreover, as $\delta(h,h) = 0$, we have that  $a= \sum_{i\in I} \alpha_ip_i + \sum_{j\in J} \alpha_j p_j=  \sum_{i\in I'} \beta_iq_i + \sum_{j\in J'} \beta_j q_j=b$, as $\min\{\alpha_i \colon I\cup J\} = 0 =\min\{\beta_i\colon I\cup J\}$.  
This implies that $U_{p_I} a = U_{q_{I'}}b$ and $U_{p_J} a = U_{q_{J'}}b$, so that 
\[
\sum_{i\in I} \alpha_ip_i = \sum_{i\in I'} \beta_iq_i \mbox{\quad and\quad }\sum_{j\in J} \alpha_j p_j = \sum_{j\in J'} \beta_j q_j. 
\]
Using the map $v\in V\mapsto e^{-v}$ we deduce that $\sum_{i\in I} e^{-\alpha_i}p_i +(u-p_I)  = \sum_{i\in I'} e^{-\beta_i}q_i +(u-q_{I'})$,  and hence  $\sum_{i\in I} e^{-\alpha_i}p_i  = \sum_{i\in I'}e^{-\beta_i}q_i$. Likewise $\sum_{j\in J} e^{-\alpha_j}p_j  = \sum_{j\in J'}e^{-\beta_j}q_j$.  We also find that 
\begin{eqnarray*}
\sum_{i\in I} e^{-\alpha_i}+\sum_{j\in J} e^{-\alpha_j} &= &(\sum_{i\in I} e^{-\alpha_i}p_i +\sum_{j\in J} e^{-\alpha_j}p_j |p_I +p_J) \\
 &= &( \sum_{i\in I'}e^{-\beta_i}q_i+\sum_{j\in J'}e^{-\beta_j}q_j|q_{I'}+q_{J'}) \\
 &= &\sum_{i\in I'}e^{-\beta_i}+ \sum_{j\in J'}e^{-\beta_j},
\end{eqnarray*}
and hence $\phi(h)$ is well defined.

We  will not only prove that $\phi\colon \overline{V}^h\to B^*_1$ is homeomorphism, but also show that $\phi$ maps each part of the horofunction boundary onto the relative interior of a boundary face of the dual unit ball. Indeed, we will establish the following theorem.

\begin{theorem}\label{hom} Given a Euclidean Jordan algebra $(V,\|\cdot\|_u)$, the map $\phi\colon \overline{V}^h\to B^*_1$ is a homeomorphism. Moreover the part $\mathcal{P}_h$, with $h$ given by (\ref{horofunction1}), is mapped onto the relative interior of the closed boundary face
\[
\mathrm{conv}\, ( U_{p_I}(V)\cap S(V))\cup(U_{p_J}(V)\cap -S(V))).
\]
\end{theorem}

\subsection{Horofunctions}
In this subsection we will prove Theorem \ref{main1}. We first make some preliminary observations. 
Note that $x\leq \lambda u$ if and only if $0\leq \lambda u -x$, which by the Hahn-Banach separation theorem is equivalent to $(\lambda u -x|w)\geq 0$ for all $w\in S(V)$.  As the state space is compact, we have for each $x\in V$ that 
\begin{equation}\label{states}
\Lambda(x) = \max_{w\in S(V)} (x|w).
 \end{equation}

As $\|\cdot\|_u$ is the JB-algebra norm,  $\|x\bullet y\|_u \leq \|x\|_u\|y\|_u$,  see \cite[Theorem 1.11]{AS1}.
It follows that if $x^n\to x$ and $y^n\to y$ in $(V,\|\cdot\|_u)$, then $x^n\bullet y^n\to x\bullet y$, since 
\[
\|x^n\bullet y^n - x\bullet y\|_u\leq \|x^n\bullet (y^n -y)\|_u+\|(x^n -x)\bullet y\|_u\leq \|x^n\|_u\|y^n-y\|_u+\|x^n -x\|_e\|y\|_u.
\]
Thus, we have the following lemma. 
\begin{lemma}\label{conv} If $x^n\to x$ and $y^n\to y$ in $(V,\|\cdot\|_u)$, then $U_{x^n} y^n \to U_x y$.
\end{lemma}
We will also use the following technical lemma several times. 
\begin{lemma}\label{lem:C} For $n\geq 1$, let $p^n_1,\ldots,p^n_r$ be a Jordan frame in $V$ and 
$I\subseteq \{1,\ldots,r\}$ nonempty. Suppose that 
\begin{enumerate}[(i)]
\item  $p_i^n\to p_i$ for all $i\in I$. 
\item $x^n \in V(p_I^n)$ with $x^n\to x\in V(p_I)$.
\item $\beta_i^n\geq 0$  with $\beta_i^n\to \beta_i\in [0,\infty]$ for all $i\in I$. 
\end{enumerate} 
If $I'=\{i\in I\colon \beta_i<\infty\}$ is nonempty, then 
\[
\lim_{n\to\infty}\Lambda_{V(p_I^n)}(x^n - \sum_{i\in I}\beta_i^np_i^n) =\Lambda_{V(p_{I'})}(U_{p_{I'}}x-\sum_{i\in I'}\beta_ip_i).
\] 
\end{lemma}
\begin{proof} We will show that every subsequence of $(\Lambda_{V(p_I^n)}(x^n - \sum_{i\in I}\beta_i^np_i^n))$ has a convergent subsequence with limit $\Lambda_{V(p_{I'})}(U_{p_{I'}}x-\sum_{i\in I'}\beta_ip_i)$. So let 
$(\Lambda_{V(p_I^{n_k})}(x^{n_k} - \sum_{i\in I}\beta_i^{n_k}p_i^{n_k}))$ be a subsequence. By (\ref{states}) there exists $d^{n_k}\in S(V(p_I^{n_k}))$ with 
\[
\Lambda_{V(p_I^{n_k})}(x^{n_k} - \sum_{i\in I}\beta_i^{n_k}p_i^{n_k}) = (x^{n_k} - \sum_{i\in I}\beta_i^{n_k}p_i^{n_k}|d^{n_k}).
\]
By taking subsequences we may assume that $d^{n_k}\to d\in S(V(p_I))$. 

Using the Peirce decomposition with respect to the Jordan frame $p^{n_k}_i$, $i\in I$, in $V(p_I^{n_k})$, we can write 
\[
d^{n_k} =\sum_{i\in I} \mu_i^{n_k} p^{n_k}_i +\sum_{i<j\in I} d_{ij}^{n_k}.
\] 
Note that as $d^{n_k}\geq 0$, we have that $\mu_i^{n_k} =(d^{n_k}|p^{n_k}_i) \geq 0$ for all $i\in I$. 

We claim that for each $i\in I\setminus I'$ we have that $\mu_i^{n_k}\to 0$. Indeed, as $I'$ is nonempty, there exist $l\in I'$ and a constant $C>0$ such that 
\[
(x^{n_k} -\sum_{i\in I} \beta_i^{n_k}p_i^{n_k}|d^{n_k})\geq (x^{n_k} -\sum_{i\in I} \beta_i^{n_k}p_i^{n_k}|p^{n_k}_l) = (x^{n_k}|p^{n_k}_l) -\beta_l^{n_k}\geq -\|x^{n_k}\|_u -\beta_l^{n_k}>-C 
\]
for all $k$, since $(x^{n_k}|p^{n_k}_l) \leq \|x^{n_k}\|_u$. Moreover,  
\[
(x^{n_k} -\sum_{i\in I} \beta_i^{n_k}p_i^{n_k}|d^{n_k}) = (x^{n_k}|d^{m_k})  -\sum_{i\in I} \beta_i^{n_k}\mu_i^{n_k}\leq \|x^{n_k}\|_u - \sum_{i\in I'} \beta_i^{n_k}\mu_i^{n_k}- \sum_{i\in I\setminus I'} \beta_i^{n_k}\mu_i^{n_k}.
\]
As $\beta_i^{n_k},\mu_i^{n_k}\geq 0$ for all $i\in I$ and $\beta_i^{n_k}\to\infty$ for all $i\in I\setminus I'$, we conclude from the previous two inequalities that $\mu_i^{n_k}\to 0$ for all $i\in I\setminus I'$. 

Using the Peirce decomposition with respect to the Jordan frame $p_i$, $i\in I$, we write 
\[
d = \sum_{i\in I}\mu_ip_i +\sum_{i<j\in I} d_{ij}.
\]
We now show that 
\begin{equation}\label{eq1c}
d=\sum_{i\in I'} \mu_{i} p_i +\sum_{i<j\in I'} d_{ij}, 
\end{equation}
and hence $d\in V(p_{I'})$.
Note that 
\[
\mu_i -\mu_i^{n_k}=(d|p_i)-(d^{n_k}|p^{n_k}_i)  = (d-d^{n_k}|p_i)+(d^{n_k}|p_i-p_i^{n_k})\to 0.
\]
We conclude that $\mu_i^{n_k} \to \mu_i$ for all $i\in I$, and hence $(d|p_j)=\mu_j=0$ for all $j\in I\setminus I'$. 
This implies by \cite[III, Exercise 3]{FK} that $d\bullet p_j=0$ for all $j\in I\setminus I'$. So,
\[
0 =d\bullet p_j = \frac{1}{2}\left(\sum_{l<j} d_{lj} + \sum_{j<m} d_{jm}\right), 
\]
which shows that $d_{lj}=0 =d_{jm}$ for all $l<j<m$, as they are all orthogonal. This implies (\ref{eq1c}).  

Next we show that 
$\lim_{k\to\infty} \Lambda_{V(p_I^{n_k})}(x^{n_k} - \sum_{i\in I}\beta_i^{n_k}p_i^{n_k}) = (U_{p_{I'}}x - \sum_{i\in I'}\beta_ip_i|d)$. First note that 
\begin{eqnarray*} \Lambda_{V(p_I^{n_k})}(x^{n_k} - \sum_{i\in I}\beta_i^{n_k}p_i^{n_k})
& = & 
(x^{n_k} -\sum_{i\in I'} \beta_i^{n_k}p_i^{n_k}|d^{n_k}) -\sum_{i\in I\setminus I'} (\beta_i^{n_k}p_i^{n_k}|d^{n_k})\\
& = & 
(x^{n_k} -\sum_{i\in I'} \beta_i^{n_k}p_i^{n_k}|d^{n_k}) -\sum_{i\in I\setminus I'} \beta_i^{n_k}\mu_i^{n_k}\\
 & \leq &
(x^{n_k} -\sum_{i\in I'} \beta_i^{n_k}p_i^{n_k}|d^{n_k})
\end{eqnarray*}
as $\beta_i^{n_k},\mu_i^{n_k}\geq 0$ for all $i$ and $k$.  This implies that 
\[
\limsup_{k\to\infty}  \Lambda_{V(p_I^{n_k})}(x^{n_k} - \sum_{i\in I}\beta_i^{n_k}p_i^{n_k})
 \leq \lim_{k\to\infty} (x^{n_k} -\sum_{i\in I'} \beta_i^{n_k}p_i^{n_k}|d^{n_k})
 =  (x -\sum_{i\in I'} \beta_ip_i|d)
\]
As $U_{p_{I'}}d =d$ and $U_{p_{I'}}$ is self-adjoint, we find that 
\[
(x -\sum_{i\in I'} \beta_ip_i|d) = (x -\sum_{i\in I'} \beta_i p_i|U_{p_{I'}}d) = (U_{p_{I'}}x -\sum_{i\in I'} \beta_ip_i|d),
\]
so that 
\begin{equation}\label{eq3c}
\limsup_{k\to\infty}  \Lambda_{V(p_I^{n_k})}(x^{n_k} - \sum_{i\in I}\beta_i^{n_k}p_i^{n_k}) \leq (U_{p_{I'}}x -\sum_{i\in I'} \beta_ip_i|d).
\end{equation}
Now let $p^{n_k}_{I'} = \sum_{i\in I'} p_i^{n_k}$.  As $p^{n_k}_{I'}\to p_{I'}$, it follows from Lemma \ref{conv}  that $U_{p^{n_k}_{I'}}d\to U_{p_{I'}}d = d$. This implies that 
 \[
(x^{n_k} -\sum_{i\in I} \beta_i^{n_k}p_i^{n_k}|U_{p_{I'}^{n_k}}d)(U_{p_{I'}^{n_k}}d|p_{I}^{n_k})^{-1}\leq 
 \Lambda_{V(p_I^{n_k})}(x^{n_k} - \sum_{i\in I}\beta_i^{n_k}p_i^{n_k})
 \]
 for all $k$ large, as $(U_{p^{n_k}_{I'}}d|p^{n_k}_{I})\to(U_{p_{I'}}d|p_{I})= (d|U_{p_{I'}}p_{I}) =(d|p_{I'}) =(d|p_{I})=1$.   Moreover,  
\begin{eqnarray*}
\lim_{k\to\infty} (x^{n_k} -\sum_{i\in I} \beta_i^{n_k}p_i^{n_k}|U_{p_{I'}^{n_k}}d)(U_{p_{I'}^{n_k}}d|p_{I}^{n_k})^{-1} & = & 
 \lim_{k\to\infty}(U_{p_{I'}^{n_k}}x^{n_k} -\sum_{i\in I'} \beta_i^{n_k}p_i^{n_k}|d)(U_{p_{I'}^{n_k}}d|p_{I}^{n_k})^{-1}\\
  & = & 
(U_{p_{I'}}x -\sum_{i\in I'} \beta_ip_i|d).
\end{eqnarray*}
 This shows that  
$(U_{p_{I'}}x -\sum_{i\in I'} \beta_ip_i|d)\leq \liminf_{k\to\infty}  \Lambda_{V(p_I^{n_k})}(x^{n_k} - \sum_{i\in I}\beta_i^{n_k}p_i^{n_k})$. From (\ref{eq3c}) we conclude that 
 \[
(U_{p_{I'}}x -\sum_{i\in I'} \beta_ip_i|d)= \lim_{k\to\infty}  \Lambda_{V(p_I^{n_k})}(x^{n_k} - \sum_{i\in I}\beta_i^{n_k}p_i^{n_k}).
 \]

To complete the proof we show that 
\begin{equation}\label{eq:5c}
(U_{p_{I'}}x -\sum_{i\in I'} \beta_ip_i|d) = \Lambda_{V(p_{I'})}(U_{p_{I'}}x -\sum_{i\in I'}\beta_ip_i).
\end{equation}
As $(d|p_{I'})=(d|p_{I})=1$, we know that $d\in S(V_{p_{I'}})$, and hence we get from by (\ref{states}) that  
\[
(U_{p_{I'}}x -\sum_{i\in I'} \beta_ip_i|d)\leq \sup_{z\in S(V(p_{I'}))} (U_{p_{I'}}x -\sum_{i\in I'} \beta_ip_i|z), 
\]
so that $(U_{p_{I'}}x -\sum_{i\in I'} \beta_ip_i|d) \leq \Lambda_{V(p_{I'})}(U_{p_{I'}}x -\sum_{i\in I'}\beta_ip_i)$.
On the other hand, if $w\in S(V(p_{I'}))$  is such that 
\[
(U_{p_{I'}}x -\sum_{i\in I'} \beta_ip_i|w) =   \sup_{z\in S(V(p_{I'}))} (U_{p_{I'}}x -\sum_{i\in I'} \beta_ip_i|z),
\]
then by definition of $d^{n_k}$ we get for all $k$ large that
\begin{eqnarray*}
(x -\sum_{i\in I} \beta_i^{n_k} p_i^{n_k}|d^{n_k}) & \geq & 
			(x -\sum_{i\in I} \beta_i^{n_k} p_i^{n_k}|U_{p^{n_k}_{I'}}w)(U_{p^{n_k}_{I'}}w|p^{n_k}_I)^{-1}\\
			  & =& 
			(U_{p^{n_k}_{I'}}x -\sum_{i\in I'} \beta_i^{n_k} p_i^{n_k}|w)(U_{p^{n_k}_{I'}}w|p^{n_k}_I)^{-1}, 
\end{eqnarray*}
as $(U_{p^{n_k}_{I'}}w|p^{n_k}_I)\to (U_{p_{I'}}w|p_I) = (w|p_{I'})=1$.

This implies that 
\[\lim_{k\to\infty}  \Lambda_{V(p_I^{n_k})}(x^{n_k} - \sum_{i\in I}\beta_i^{n_k}p_i^{n_k})
 \geq \lim_{k\to\infty}(U_{p^{n_k}_{I'}}x -\sum_{i\in I'} \beta_i^{n_k} p_i^{n_k}|w)(U_{p^{n_k}_{I'}}w|p^{n_k}_I)^{-1} = 
 (U_{p_{I'}}x -\sum_{i\in I'} \beta_i p_i|w),
 \]
 and hence  (\ref{eq:5c}) holds. 
\end{proof}

To prove that all horofunctions in $\overline{V}^h$ are of the form (\ref{horofunction1}), we first establish the following proposition by using the previous lemma.
\begin{proposition}\label{prop:3.1} Let $(y^n)$ be a sequence in $V$, with $y^n =\sum_{i=1}^r \lambda_i^n p^n_i$. Suppose that  $h_{y^n}\to h\in\partial \overline{V}^h$ and $(y^n)$ satisfies the following properties: 
\begin{enumerate}[(1)] 
\item There exists $1\leq s\leq r$ such  that $|\lambda^n_s|=r^n$ for all $n$, where $r^n =\|y^n\|_u$.  
\item $p_k^n\to p_k$ for all $1\leq k\leq r$. 
\item There exist $I,J\subseteq \{1,\ldots,r\}$ disjoint with $I\cup J$ nonempty, and $\alpha\in \mathbb{R}^{I\cup J}$  with $\min\{\alpha_i\colon i\in I\cup J\}=0$ such that $r^n - \lambda_i^n\to \alpha_i$ for all $i\in I$, $r^n+\lambda_j^n\to \alpha_j$ for all $j\in J$, and $r^n- |\lambda_k^n|\to \infty$ for all $k\not\in I\cup J$.
\end{enumerate}
Then  $h$ satisfies (\ref{horofunction1}). 
\end{proposition}
\begin{proof}
Take $x\in V$ fixed. Note that for all $n\geq 1$, 
\[
\|x-y^n\|_u -\|y^n\|_u = \max\{\Lambda(x -y^n),\Lambda(-x+y^n)\} -r^n = \max\{ \Lambda(x -y^n -r^nu),\Lambda(-x+y^n -r^nu)\}
\]
As $h$ is a horofunction, $\|y^n\|_u = r^n\to\infty$ by Lemma \ref{Rieffel}. Thus,  $\lambda^n_i\to\infty$ for all $i\in I$ and $\lambda_j^n\to -\infty$ for all $j\in J$. 
Now suppose that $J$ is nonempty. Then $r^n+\lambda_k^n\geq r^n-|\lambda_k^n|\to\infty$ for all $k\not\in I\cup J$. As 
\[
\Lambda(x -y^n -r^nu) = \Lambda(x - \sum_{j\in J} (r^n+\lambda_j^n)p^n_j - \sum_{k\not \in J} (r^n+\lambda_k^n)p^n_k),
\]
it follows that 
\[
\lim_{n\to\infty} \Lambda(x -y^n -r^nu) = \Lambda_{V(p_J)}(U_{p_J}x -\sum_{j\in J}\alpha_jp_j)
\]
by Lemma \ref{lem:C}. Likewise, if $I$ is nonempty, then  
\[
\lim_{n\to\infty} \Lambda(-x +y^n -r^nu) = \Lambda_{V(p_I)}(-U_{p_I}x -\sum_{i\in I}\alpha_ip_i)
\]
by Lemma \ref{lem:C}.
We conclude that if $I$ and $J$ are both nonempty, then 
\begin{eqnarray*}
h(x) = \lim_{n\to\infty} \|x-y^n\|_u -\|y^n\|_u&  = & \lim_{n\to\infty} \max\{\Lambda( -x +y^n-r^nu),\Lambda( x -y^n-r^nu) \} \\
 & = & \max\{\Lambda_{V(p_I)}(-U_{p_I}x -\sum_{i\in I}\alpha_ip_i) , \Lambda_{V(p_J)}(U_{p_J}x -\sum_{j\in J}\alpha_jp_j)\}.
\end{eqnarray*}

To complete the proof it remains to show that $\lim_{n\to\infty} \|x-y^n\|_u -\|y^n\|_u =  \lim_{n\to\infty} \Lambda( -x +y^n-r^nu)$ if $J$ is empty, and 
 $\lim_{n\to\infty} \|x-y^n\|_u -\|y^n\|_u =  \lim_{n\to\infty} \Lambda( x -y^n-r^nu)$ if $I$ is empty.
Suppose that $I$ empty, so $J$ is nonempty. For each $i\in\{1,\ldots,r\}$ we have that $r^n-\lambda_i^n\to\infty$. 
Note that 
\[
-x+y^n-r^nu= -x -\sum_i (r^n-\lambda_i^n)p_i^n\leq -x -\min_i (r^n-\lambda_i^n)u\leq (\|x\|_u-\min_i (r^n-\lambda_i^n))u.
\]
Thus, $\Lambda(-x+y^n-r^nu) \leq  \Lambda((\|x\|_u-\min_i (r^n-\lambda_i^n))u) =  \|x\|_u-\min_i (r^n-\lambda_i^n)$ for all $n$,  and hence 
$\Lambda(-x+y^n-r^nu)\to-\infty$. As 
\[
\max\{\Lambda(x-y^n-r^nu),\Lambda(-x+y^n-r^nu)\} = \|x-y^n\|_u -\|y^n\|_u\geq -\|x\|_u>-\infty,
\]
we conclude that $\|x-y^n\|_u -\|y^n\|_u = \Lambda(x-y^n-r^nu)$ for all $n$ sufficiently large, and hence 
\[
h(x) = \lim_{n\to\infty}  \Lambda(x-y^n-r^nu) = \Lambda_{V(p_J)}(U_{p_J}x -\sum_{j\in J}\alpha_jp_j).
\]
The argument for the case where $J$ is empty goes in the same way. 
\end{proof}
The following corollary shows that each horofunction is of the form (\ref{horofunction1}). 
 \begin{corollary}\label{cor1}  If $h$ is a horofunction in $\overline{V}^h$, then there exist a Jordan frame $p_1,\ldots,p_r$ in $V$, disjoint subsets $I,J\subseteq\{1,\ldots,r\}$, with $I\cup J$ nonempty, and 
$\alpha\in\mathbb{R}^{I\cup J}$ with $\min\{\alpha_i \colon i\in I\cup J\}=0$, such that $h\colon V\to \mathbb{R}$ satisfies (\ref{horofunction1}) 
for all $x\in V$. 
\end{corollary}
\begin{proof} Suppose that $(y^n)$ is a sequence in $V$  with $h_{y^n}\to h$ in $\overline{V}^h$. Then  for each $x\in V$ we have that  
\[
\lim_{n\to\infty} \|x-y^n\|_u -\|y^n\|_u = h(x)
\] 
and $\|y^n\|_u\to\infty$ by Lemma \ref{Rieffel}.  

To show that the limit is equal to (\ref{horofunction1}) it suffices to show that we can take  a subsequences of $(y^n)$ that satisfies the conditions in Proposition \ref{prop:3.1}. First we note that by the spectral theorem \cite[Theorem III.1.2]{FK}, there exist for each $n\geq 1$ a Jordan frame $p^n_1,\ldots,p^n_r$ in $V$ and   $\lambda^n_1,\ldots,\lambda^n_r\in\mathbb{R}$ such that 
\[
y^n = \lambda^n_1p^n_1+\cdots+\lambda^n_r p^n_r,
\]
where $r$ is the rank of $V$.  Denote $r^n=\|y^n\|_u = \max_i |\lambda^n_i|$. 

Now by taking subsequences we may assume that there exist $I_+\subseteq \{1,\ldots, r\}$ and $1\leq s\leq r$ such that for each $n\geq 1$  we have $r^n = |\lambda^n_s|$ and  
\[
\lambda_i^n>0\mbox{ for all }i\in I_+ \mbox{\quad and \quad}\lambda^n_i\leq 0 \mbox{ for all }i\not\in I_+. 
\]
Now for each $i\in\{1,\ldots,r\}$ and $n\geq 1$ define 
\[
\alpha_i^n = \left[\begin{array}{cl} r^n - \lambda_i^n & \mbox{ for } i\in I_+\\ r^n +\lambda^n_i & \mbox{ for } i\not\in I_+.\end{array}\right.
\]
Note that $\alpha^n_i \in [0,\infty)$ for all $i$. Again by taking subsequences we may assume that $\alpha_i^n\to \alpha_i\in [0,\infty]$, as $n\to \infty$, for all $i$. Recall that  $\alpha^n_s = 0$ for all $n$, so  $\alpha_s=0$. Furthermore we may assume that $p^n_i\to p_i$ in $\mathcal{J}_1(V)$ for all $i$, as it is a compact set \cite{Hi}. Note that $p_1,\ldots,p_r$ is a Jordan frame in $V$. 
  
Now let 
\[
I =\{i\colon \alpha_i<\infty\mbox{ and } i\in I_+\}\mbox{\quad and \quad} J=\{j\colon \alpha_j<\infty\mbox{ and }j\not\in I_+\}.
\]
So $I\cap J$ is empty, $s\in I\cup J$ and $\min\{\alpha_i\colon i\in I\cup J\} =\alpha_s=0$. 
Then the subsequence of $(y^n)$ satisfies the conditions in Proposition \ref{prop:3.1}, and hence $h$ is a horofunction of the form  (\ref{horofunction1}).
\end{proof}

The next proposition shows that each function of the form (\ref{horofunction1}) can be realised as a horofunction, and is a Busemann point. 
 \begin{proposition}  \label{prop:2}
Let $p_1,\ldots,p_r$ be a Jordan frame in $V$. Given $I,J\subseteq\{1,\ldots,r\}$, with $I\cap J=\emptyset$ and $I\cup J$ nonempty, and 
$\alpha\in\mathbb{R}^{I\cup J}$ with $\min\{\alpha_i \colon i\in I\cup J\}=0$, 
For $n\geq 1$ let  $y^n = \lambda_1^np_1+\cdots+\lambda_r^n p_r$, where 
\[
\lambda_i^n =\left[\begin{array}{rl}
n -\alpha_i & \mbox{ if } i\in I\\
-n+\alpha_i & \mbox{ if } i\in J\\
0& \mbox{ otherwise.}
\end{array}\right.
\]
Then $(y^n)$ is an almost geodesic sequence and $h_{y^n}\to h$ where  $h$ satisfies (\ref{horofunction1}) for all $x\in V$. In particular, $h$ is a Busemann point in $\overline{V}^h$. 
\end{proposition}
\begin{proof}
We will use Proposition \ref{prop:3.1}. Let $k\geq \max\{ \alpha_i\colon i\in I\cup J\}$ and note that  for $n\geq k$ we have that $r^n = \|y^n\|_u= n$, as  $ \min\{\alpha_i \colon i\in I\cup J\}=0$. The sequence $(y^n)$, where $n\geq k$, satisfies the conditions in Proposition \ref{prop:3.1}. Indeed, for $n\geq k$ we have that 
$r^n -\lambda_i^n = \alpha_i$ for all $i\in I$, $r^n+\lambda_i^n =\alpha_i$ for all $i\in J$, and $r^n -\lambda_i^n= n$ otherwise. Also for $1\leq s\leq r$ with $\alpha_s = \min\{\alpha_i \colon i\in I\cup J\}$, we have that 
$|\lambda_s^n| = n = \|y^n\|_u$. 

Finally to see that $(h_{y^n})$ converges, we note that  if we define $z = \sum_{i\in I} -\alpha_ip_i + \sum_{j\in J} \alpha_jp_j$ and $w= \sum_{i\in I} p_i -\sum_{j\in J}p_j$, then $y^n= nw+z$, which lies on the straight-line 
$t\mapsto tw+z$. Hence $(y^n)$ is an almost geodesic sequence, so 
\[
h(x) = \lim_{n\to\infty} \|x-y^n\|_u -\|y^n\|_u
\]
exists of all $x\in V$. Thus, we can apply Proposition \ref{prop:3.1} and conclude that $h$ satisfies (\ref{horofunction1}). Moreover, as  $(y^n)$ is an almost geodesic sequence,  $h$ is a Busemann point in the horofunction boundary.
\end{proof}
 Combining the results so far we now prove Theorem \ref{main1}. 
 \begin{proof}[Proof of Theorem \ref{main1}]
 Corollary \ref{cor1} shows that each horofunction in $\overline{V}^h$ is of the form (\ref{horofunction1}). It follows from Proposition \ref{prop:2} that any function of the form (\ref{horofunction1}) is a horofunction and  by the second part of that proposition each horofunction is a Busemann point. 
 \end{proof}

 \subsection{Parts and the detour metric}
 In this subsection we will identify the parts in the horofunction boundary of $\overline{V}^h$, derive a formula for the detour distance, and establish Theorem \ref{main2}. We begin by proving the following proposition. 
 
 \begin{proposition}  \label{prop:4.9}
If
\begin{equation}\label{h1}
h(x)  = \max\left\{\Lambda_{V(p_{I})}\left(-U_{p_{I}}x - \sum_{i\in I} \alpha_ip_i\right),  \Lambda_{V_(p_{J})}\left(U_{p_{J}}x - \sum_{j\in J} \alpha_j p_j\right)\right\},
\end{equation}
and 
\begin{equation}\label{h2}
h'(x)  = \max\left\{\Lambda_{V(q_{I'})}\left(-U_{q_{I'}}x - \sum_{i\in I'} \beta_iq_i\right),  \Lambda_{V(q_{J'})}\left(U_{q_{J'}}x - \sum_{j\in J'} \beta_j q_j\right)\right\},
\end{equation}
are horofunctions  with $p_{I} = q_{I'}$ and $p_{J}=q_{J'}$, then $h$ and $h'$ are in the same part and  
\[
\delta(h,h') = \|a-b \|_\mathrm{var} = \Lambda_{V(p_I,p_J)}(a-b) + \Lambda_{V(p_I,p_J)}(b-a), 
\]
where $a =\sum_{i\in I} \alpha_ip_i + \sum_{j\in J} \alpha_jp_j$ and $b =  \sum_{i\in I'} \beta_iq_i + \sum_{j\in J'} \beta_j q_j$ in $V(p_I,p_J)$.
\end{proposition}
\begin{proof} 
As in Proposition \ref{prop:2}, for $n\geq 1$ let  $y^n = \lambda_1^np_1+\cdots+\lambda_r^n p_r$, where 
\[
\lambda_i^n =\left[\begin{array}{rl}
n -\alpha_i & \mbox{ if } i\in I\\
-n+\alpha_i & \mbox{ if } i\in J\\
0& \mbox{ otherwise.}
\end{array}\right.
\]
and let  $w^n = \mu_1^nq_1+\cdots+\mu_r^n q_r$, where 
\[
\mu_i^n =\left[\begin{array}{rl}
n -\beta_i & \mbox{ if } i\in I'\\
-n+\beta_i & \mbox{ if } i\in J'\\
0& \mbox{ otherwise.}
\end{array}\right.
\]
By Proposition \ref{prop:2} we know that $(y^n)$ and $(w^n)$ are almost geodesic sequences with $h_{y^n}\to h$ and $h_{w^n}\to h'$. Note that 
\[
U_{p_I} w^m= U_{q_{I'}} w^m = \sum_{i\in I'}\mu_i^m U_{q_{I'}}q_i = \sum_{i\in I'}\mu_i^m q_i
\]
for all $m$, so that 
\begin{eqnarray*}
\Lambda_{V(p_I) }(-U_{p_I}w^m -\sum_{i\in I}\alpha_ip_i +\|w^m\|_u p_I) & =  &\Lambda_{V(p_I)}(-U_{q_{I'}}w^m -\sum_{i\in I} \alpha_ip_i +\|w^m\|_uq_{I'}) \\
  & = &  \Lambda_{V(p_I)} (\sum_{i\in I'} (\|w^m\|_u -\mu_i^m)q_i -\sum_{i\in I} \alpha_ip_i ).
\end{eqnarray*}
Thus, 
\begin{eqnarray*}
\lim_{m\to\infty} \Lambda_{V(p_I) }(-U_{p_I}w^m -\sum_{i\in I}\alpha_ip_i +\|w^m\|_u p_I) &= &\lim_{m\to\infty}  \Lambda_{V(p_I)} (\sum_{i\in I'} (\|w^m\|_u -\mu_i^m)q_i -\sum_{i\in I} \alpha_ip_i ) \\
&= &  \Lambda_{V(p_I)} (\sum_{i\in I'} \beta_iq_i -\sum_{i\in I} \alpha_ip_i )\\
 & = &  \Lambda_{V(p_I)}(b-a).
\end{eqnarray*}

In the same way it can be shown that 
\begin{equation*}
\lim_{m\to\infty} \Lambda_{V(p_J)}(U_{p_J}w^m -\sum_{j\in J}\alpha_j p_i +\|w^m\|_u p_J) =   \Lambda_{V(p_J)} (\sum_{j\in J'} \beta_jq_j -\sum_{j\in J} \alpha_jp_j) =  \Lambda_{V(p_J)}(b-a).
\end{equation*}

So, it follows from (\ref{detourcost}) that 
\begin{eqnarray*} 
H(h,h')    & = & \lim_{m\to\infty} \|w^m\|_u + \max\{\Lambda_{V(p_I)}(-U_{p_I}w^m - \sum_{i\in I} \alpha_ip_i),  \Lambda_{V(p_J)}(U_{p_J}w^m - \sum_{j\in J} \alpha_j p_j)\} \\
    & = & \lim_{m\to\infty} \max\{\Lambda_{V(p_I)}(-U_{p_I}w^m - \sum_{i\in I} \alpha_ip_i + \|w^m\|_u p_I),  \Lambda_{V(p_J)}(U_{p_J}w^m - \sum_{j\in J} 
   \alpha_j p_j  + \|w^m\|_u p_J)\} \\
 & = & \max\{  \Lambda_{V(p_I)} (\sum_{i\in I'} \beta_iq_i -\sum_{i\in I} \alpha_ip_i ),   \Lambda_{V(p_J)} (\sum_{j\in J'} \beta_jq_j -\sum_{j\in J} \alpha_jp_j)\}\\
 & = &  \Lambda_{V(p_I,p_J)}(b-a).
 \end{eqnarray*}

Interchanging the roles of $h$ and $h'$ gives 
\[
H(h',h) = \Lambda_{V(p_I,p_J)}(a-b),
\]
and hence $\delta(h,h') = \|a-b\|_\mathrm{var}$. 
\end{proof}
To show that $h$ and $h'$ are in different parts, if $p_I\neq q_{I'}$ or $p_j\neq q_{J'}$, we need the following lemma. 
\begin{lemma}\label{pneq}
If $p$ and $q$ are idempotents in $V$ with $p\nleq q$, then $U_pq<p$.
\end{lemma}
\begin{proof}
We have that $U_pq\leq U_p u = p$. In fact, $U_p q < p$. Indeed, if $U_p q= p$, then 
\[
p = U_p u = U_p (u-q) +U_p q = U_p( u-q) +p,
\]
and hence $U_p (u-q) =0$. This implies that $p +(u-q) \leq u$ by \cite[Lemma 4.2.2]{HOS}, so that $p\leq q$. This is impossible, as $p\nleq q$, and hence $U_p q < p$.
\end{proof}
\begin{proposition}  \label{prop:4.11}
If $h$ and $h'$ are horofunctions given by (\ref{h1}) and (\ref{h2}), respectively,  and $p_{I} \neq q_{I'}$ or $p_{J}\neq q_{J'}$, then 
\[
\delta(h,h') = \infty. 
\]
\end{proposition}
\begin{proof}
Suppose that $p_I\neq q_{I'}$. Then $p_I\nleq q_{I'}$ or $q_{I'}\nleq p_I$. Without loss of generality assume that $p_I\nleq q_{I'}$. 
Let $(y^n)$ in $V(p_I)$ and $(w^n)$ in $V(q_{I'})$ be as in Proposition \ref{prop:2}, so $h_{y^n}\to h$ and $h_{w^m}\to h'$. To prove the statement in this case, we use (\ref{detourcost}) and show that 
\begin{equation}\label{infty}
H(h',h) = \lim_{m\to\infty} \|w^m\|_u + h(w^m) = \infty. 
\end{equation}
Note that 
\[
\|w^m\|_u +h(w^m) \geq \|w^m\|_u + \Lambda_{V(p_I)}(-U_{p_I}w^m -\sum_{i\in I} \alpha_i p_i ) = \Lambda_{V(p_I)}(-U_{p_I}w^m -\sum_{i\in I} \alpha_i p_i +\|w^m\|_u p_I). 
\]
As $w^m\leq \|w^m\|_u q_{I'}$ for all $m$ large, we have that $U_{p_I}w^m\leq \|w^m\|_uU_{p_I}q_{I'}$ for all $m$ large.  Thus, 
\begin{eqnarray*}
-U_{p_I}w^m -\sum_{i\in I} \alpha_i p_i +\|w^m\|_u p_I & \geq & -\|w^m\|_uU_{p_I}q_{I'} -\sum_{i\in I} \alpha_i p_i +\|w^m\|_u p_I\\
   & = & \|w^m\|_u (p_I -U_{p_I}q_{I'}) - \sum_{i\in I} \alpha_i p_i
\end{eqnarray*}
for all $m$ large. 

We know from Lemma \ref{pneq} that  $p_I-U_{p_I}q_{I'}>0$.  As $p_I-U_{p_I}q_{I'}\in V(p_I)$ we also have  that $p_I-U_{p_I}q_{I'} = \sum_{j=1}^s \gamma_j r_j$, where $\gamma_j>0$ for all $j$ and the $r_j$'s are orthogonal idempotents in $V(p_I)$. 
It now follows that for all $m$ large, 
\begin{eqnarray*}
\Lambda_{V(p_I)}(-U_{p_I}w^m -\sum_{i\in I} \alpha_i p_i +\|w^m\|_up_I) & \geq & (\|w^m\|_u \sum_{j=1}^s\gamma_j r_j - \sum_{i\in I}\alpha_i p_i|r_1)(p_I|r_1)^{-1}\\
  & =&( \|w^m\|_u\gamma_1 - ( \sum_{i\in I}\alpha_i p_i|r_1))(p_I|r_1)^{-1}.
\end{eqnarray*}
The right-hand side goes to $\infty$ as $m\to \infty$, and hence (\ref{infty}) holds. 

For the case $p_J\neq q_{J'}$ a similar argument can be used. 
\end{proof} 
We now prove Theorem \ref{main2}. 
\begin{proof}[Proof of Theorem \ref{main2}] Parts (i) and (ii) follow directly from Propositions \ref{prop:4.9} and \ref{prop:4.11}. Clearly the map $\rho\colon\mathcal{P}_h\to V(p_I,p_J)/\mathbb{R}p_{IJ}$ given by $\rho(h') = [b]$, where 
\[
h'(x)  = \max\left\{\Lambda_{V(q_{I'})}\left(-U_{q_{I'}}x - \sum_{i\in I'} \beta_iq_i\right),  \Lambda_{V(q_{J'})}\left(U_{q_{J'}}x - \sum_{j\in J'} \beta_j q_j\right)\right\},
\]
and $b =  \sum_{i\in I'} \beta_iq_i + \sum_{j\in J'} \beta_j q_j\in V(p_I,p_J)$ with $\min_{i\in I\cup J} \beta_i=0$, is a bijection. Indeed, for each $[c]\in V(p_I,p_J)/\mathbb{R}p_{IJ}$, there is a unique $c'\in [c]$ with $\min \sigma_{V(p_I,p_J)}(c') =0$.  So, by Proposition \ref{prop:4.9}, $\rho$ is an isometry from $(\mathcal{P}_h,\delta)$ onto $ (V(p_I,p_J)/\mathbb{R}p_{IJ},\|\cdot \|_\mathrm{var} )$. 
\end{proof}

\subsection{The homeomorphism onto the dual unit ball} 
In this subsection we show  Theorem \ref{hom}. To start we prove a basic lemma that will be useful in the sequel. 
\begin{lemma}\label{lem:A} 
If $q\leq p$ are idempotents in $V$ and $z\in V(p)$, then $\Lambda_{V(q)}(U_qz)\leq \Lambda_{V(p)}(z)$. 
\end{lemma}
\begin{proof}
If $\lambda = \Lambda_{V(p)}(z)$, then $0\leq \lambda p-z$, so that $0\leq \lambda U_qp -U_qz$. 
As $q= U_q q \leq U_qp\leq U_qu =q^2=q$, we find that $0\leq\lambda U_qp -U_qz =\lambda q -U_qz$, and hence  $\Lambda_{V(q)}(U_qz)\leq \lambda$. 
\end{proof}

We will show that $\phi$ given by (\ref{phi1}) and (\ref{phi2}) is a continuous bijection from $\overline{V}^h$ onto $B_1^*$. As $\overline{V}^h$  is compact and $B_1^*$ is Hausdorff, we can then conclude that $\phi$ is a homeomorphism. We begin by showing that $\phi$ maps $V$ into the interior of $B_1^*$. 
\begin{lemma}\label{homlem1} 
For each $x\in V$ we have that $\phi(x) \in\mathrm{int}\,B^*_1$. 
\end{lemma}
\begin{proof}
For $x\in V$ there exists $y\in V$ with $\|y\|_u= 1$ such that 
\[
\|\phi(x)\|_u^* = \sup_{w\in V\colon  \|w\|_u \leq 1} |(w|\phi(x))| = (y|\phi(x)),
\]
where $(v|w) =\mathrm{tr}(v\bullet w)$.  
So, if $x$ has spectral decomposition $x=\sum_{i=1}^r \lambda_ip_i$, then we can consider the Peirce decomposition of $y$, 
\[
y = \sum_{i=1}^r \mu_i p_i+\sum_{i<j} y_{ij},
\]
to find that
\[
\|\phi(x)\|_u^* = (\phi(x)|y) = \frac{1}{\sum_{i=1}^r e^{\lambda_i} +e^{-\lambda_i}} (\sum_{i=1}^r  (e^{\lambda_i} -e^{-\lambda_i})p_i| y) 
\leq \frac{ \sum_{i=1}^r (e^{\lambda_i} -e^{-\lambda_i})|\mu_i| }{\sum_{i=1}^r e^{\lambda_i} +e^{-\lambda_i}}<1,\]
as $\mu_i  =(y|p_i)\leq (u|p_i) =1$ and $\mu_i  =(y|p_i)\geq (-u|p_i) =-1$.
\end{proof}
\begin{lemma}\label{lem:inj1} The map $\phi$ is injective on $V$.
\end{lemma}
\begin{proof}
Suppose that $x,y\in V$ with $x=\sum_{i=1}^r \sigma_i p_i$ and $y=\sum_{i=1}^r \tau_i q_i$, where $\sigma_1\leq \ldots\leq \sigma_r$ and $\tau_1\leq \ldots\leq \tau_r$, satisfy $\phi(x) =\phi(y)$. Then $\phi(x) =\sum_{i=1}^r \alpha_i p_i  =\sum_{i=1}^r \beta_i q_i=\phi(y)$. where 
\[
\alpha_j =\frac{e^{\sigma_j}-e^{-\sigma_j}}{\sum_{i=1}^r e^{\sigma_i}+e^{-\sigma_i}}
\mbox{\quad and \quad}
\beta_j =\frac{e^{\tau_j}-e^{-\tau_j}}{\sum_{i=1}^r e^{\tau_i}+e^{-\tau_i}}
\]
for all $j$. As $\alpha_1 \leq \ldots\leq \alpha_r$ and $\beta_1\leq \ldots\leq \beta_r$, it follows from the spectral theorem (version 2) \cite[Theorem III.1.2]{FK} that $\alpha_j = \beta_j$ for all $j$. Lemma \ref{lem:calc} now implies that $\sigma = (\sigma_1,\ldots,\sigma_r) = (\tau_1,\ldots,\tau_r) =\tau$, as 
\[
(\alpha_1,\ldots,\alpha_r)  = \nabla\log \mu(\sigma)\mbox{\quad and \quad} 
(\beta_1,\ldots,\beta_r)  = \nabla\log \mu(\tau)
\] 
Note that $\alpha_i =\alpha_j$ if and only if $\sigma_i=\sigma_j$, and, $\beta_i =\beta_j$ if and only if $\tau_i=\tau_j$, as $\nabla \log \mu(x)$ is injective. 
It now follows from the spectral theorem (version 1) \cite[Theorem III.1.1]{FK} that $x=y$. 
\end{proof}
\begin{lemma}\label{lem:ontoint} The map $\phi$ maps $V$ onto $\mathrm{int}\,B_1^*$. 
\end{lemma}
\begin{proof}
As $\phi$ is continuous on $V$ and $\phi(V)\subseteq \mathrm{int}\,B_1^*$ it follows from Brouwer's domain invariance theorem that  $\phi(V)$ is open in $\mathrm{int}\,B_1^*$. Suppose, for the sake of contradiction, that
$\phi(V)\neq \mathrm{int}\,B_1^*$. Then $ \partial \phi(V)\cap \mathrm{int}\, B_1^*$ is nonempty, as otherwise $\phi(V)$ is closed and open, which would imply that $\mathrm{int}\, B_1^*$ is the disjoint union of  two nonempty open sets contradicting the connectedness of $\mathrm{int}\, B^*_1$. So we can find a $z\in \partial \phi(V)\cap \mathrm{int}\,B_1^*$. Let $(y^n)$ in $V$ be such that $\phi(y^n)\to z$ and write $y^n = \sum_{i=1}^r \lambda_i^np_i^n$. As $\phi$ is continuous on $V$, we may assume that $r^n =\|y^n\|_u\to\infty$. Furthermore, after taking a subsequence, we may  assume that $(y^n)$ satisfies the conditions in Proposition \ref{prop:3.1}. So, using the notation as in Proposition \ref{prop:3.1}, we get that 
\[
\phi(y^n)  =   \frac{\sum_{i=1}^r(e^{\lambda^n_i}-e^{-\lambda^n_i})p_i^n}{\sum_{i=1}^r e^{\lambda^n_i}+e^{-\lambda^n_i}} 
  =    \frac{\sum_{i=1}^r(e^{-r^n+\lambda^n_i}-e^{-r^n-\lambda^n_i})p_i^n}{\sum_{i=1}^r e^{-r^n+\lambda^n_i}+e^{-r^n-\lambda^n_i}}. \]
The right-hand side converges to 
\[ 
 \frac{1}{\sum_{i\in I} e^{-\alpha_i}+\sum_{j\in J} e^{-\alpha_j}} \left( \sum_{i\in I}e^{-\alpha_i}p_i-\sum_{j\in J}e^{-\alpha_j}p_j\right) =z.
\]
But this implies that $z\in \partial B_1^*$, which is impossible. Indeed, if we let $p_I =\sum_{i\in I} p_i$ and $p_J=\sum_{j\in J} p_j$, then  $1 \geq \|z\|_u^*\geq ( z|p_I-p_J) = 1$, as $ -u\leq p_I -p_J\leq u$.
\end{proof}
For simplicity we denote the (closed) boundary faces of $B_1^*$ by  
\[
F_{p,q} = \mathrm{conv}\,( (U_p(V)\cap S(V))\cup (U_q(V)\cap -S(V)))
\]
where $p$ and $q$ are orthogonal idempotents in $V$. 
\begin{lemma}\label{homlem2} If $h$ is a horofunction given by (\ref{horofunction1}), then $\phi$ maps $\mathcal{P}_h$ into  $ \mathrm{relint} \,F_{p_I,p_J}$.  
\end{lemma}
\begin{proof}  Clearly, $\phi(h)\in F_{p_I,p_J}$ if $h$ is given by (\ref{horofunction1}). So, $\phi$ maps $\mathcal{P}_h$ into $F_{p_I,p_J}$ by Theorem \ref{main2}(i). To show that $\phi$ maps $\mathcal{P}_h$ into $ \mathrm{relint} \,F_{p_I,p_J}$, it suffices to show that $\phi(h)\in  \mathrm{relint} \,F_{p_I,p_J}$.

To do this we first consider  $w= (|I|+|J|)^{-1}(p_I-p_J)\in F_{p_I,q_J}$ and show that $w\in \mathrm{relint} F_{p_i,q_J}$. 
Let $c\in F_{p_I,p_J}$ be arbitrary. Note  that we can write $c= \sum_{i\in I'} \lambda_i q_i -\sum_{j\in J'} \lambda_jq_j$, where 
$\sum_{i\in I'} q_i =p_I$, $\sum_{j\in J'} q_j =p_J$, and $\sum_{i\in I'} \lambda_i + \sum_{j\in J'} \lambda_j =1$ with $0\leq \lambda_i,\lambda_j\leq 1$ for all $i$ and $j$. We see that $w+\epsilon(w-c) = (1+\epsilon)w -\epsilon c\in F_{p_I,p_j}$ for all $\epsilon>0$ small, and hence $w\in\mathrm{relint}\, F_{p_i,p_j}$. 

Clearly, $\phi(h)\in F_{p_I,p_J} = \mathrm{conv}\,( (U_{p_I}(V)\cap S(V))\cup (U_{p_J}(V)\cap -S(V)))$. To complete the proof we argue by contradiction. So suppose that $\phi(h)\not\in \mathrm{relint} F_{p_I,p_J}$. Then $\phi(h)$ is in the (relative) boundary of $F_{p_I,p_J}$, and hence 
\[
z_\epsilon = (1+\epsilon) \phi(h) -\epsilon w\not\in F_{p_I,p_J}
\]
for all $\epsilon>0$, as $w\in\mathrm{relint} F_{p_I,p_J}$ and $F_{p_I,p_J}$ is convex.

However, for each $i\in I$ we have that the coefficient of $p_i$ in $z_\epsilon$,
\[
\frac{(1+\epsilon)e^{-\alpha_i} }{\sum_{i\in I} e^{-\alpha_i}+\sum_{j\in J}e^{-\alpha_j}} - \frac{\epsilon}{|I|+|J|},
\]
is strictly positive for all $\epsilon>0$ sufficiently small. Likewise,  for each $j\in J$ we have that the coefficient of $-p_j$ in $z_\epsilon$, 
\[
\frac{(1+\epsilon)e^{-\alpha_j} }{\sum_{i\in I} e^{-\alpha_i}+\sum_{j\in J}e^{-\alpha_j}} - \frac{\epsilon}{|I|+|J|},
\]
is strictly positive for all $\epsilon>0$ sufficiently small. This implies that $z_\epsilon\in F_{p_I,p_J}$ for all $\epsilon>0$ small, which is impossible. 
This completes the proof.
\end{proof}
Using  the previous results we now show that $\phi$ is injective on $\overline{V}^h$. 
\begin{corollary}\label{injective}
The map $\phi\colon \overline{V}^h\to B_1^*$ is injective.
\end{corollary}
\begin{proof}
We already saw in Lemmas \ref{homlem1} and  \ref{lem:inj1} that $\phi$ maps $V$ into $\mathrm{int}\, B_1^*$ and is injective on $V$. So by the previous lemma, it suffices to show that if $\phi(h)=\phi(h')$ for horofunctions $h\sim h'$, then $h =h'$. Let $h$ be given by (\ref{horofunction1}) and suppose that $h'$ is given by 
\[
h'(x)  = \max\left\{\Lambda_{V(q_{I'})}\left(-U_{q_{I'}}x - \sum_{i\in I'} \beta_iq_i\right),  \Lambda_{V(q_{J'})}\left(U_{q_{J'}}x - \sum_{j\in J'} \beta_j q_j\right)\right\}.
\]
Then 
\[
  \frac{\sum_{i\in I} e^{-\alpha_i}p_i -\sum_{j\in J} e^{-\alpha_j}p_j}{\sum_{i\in I} e^{-\alpha_i} +\sum_{j\in J} e^{-\alpha_j}}
  =   \frac{\sum_{i\in I'} e^{-\beta_i}q_i -\sum_{j\in J'} e^{-\beta_j}q_j}{\sum_{i\in I'} e^{-\beta_i} +\sum_{j\in J'} e^{-\beta_j}}.
\]

As $\min_k \alpha_k =0 =\min_k\beta_k$, it follows from the spectral theorem  \cite[Theorem III.1.2]{FK} that 
\[
 \frac{1}{\sum_{i\in I} e^{-\alpha_i} +\sum_{j\in J} e^{-\alpha_j}} =\|\phi(h)\|_u=\|\phi(h')\|_u =  \frac{1}{\sum_{i\in I'} e^{-\beta_i} +\sum_{j\in J'}e^{-\beta_j}}, 
\]
so that 
\[
\sum_{i\in I} e^{-\alpha_i}p_i -\sum_{j\in J} e^{-\alpha_j}p_i = \sum_{i\in I'} e^{-\beta_i}q_i -\sum_{j\in J'} e^{-\beta_j}q_j.
\]

As each $x\in V$ can be written in a unique way as $x= x^+-x^-$, where $x^+$ and $x^-$ are orthogonal element $x^+$ and $x^-$ in $V_+$, see \cite[Proposition 1.28]{AS1}, we find that  $\sum_{i\in I} e^{-\alpha_i}p_i = \sum_{i\in I'} e^{-\beta_i}q_i$ and $\sum_{j\in J} e^{-\alpha_j}p_i = \sum_{j\in J'} e^{-\beta_j}q_j$. This implies that 
\[
\sum_{i\in I} \alpha_ip_i  =-\log (\sum_{i\in I} e^{-\alpha_i}p_i +(u-p_I))=-\log \sum_{i\in I'} (e^{-\beta_i}q_i+(u-q_{I'}))= \sum_{i\in I'} \beta_iq_i
\]
and 
\[
\sum_{j\in J} \alpha_jp_i = -\log( \sum_{j\in J} e^{-\alpha_j}p_i +(u-p_J))=-\log(  \sum_{j\in J'} e^{-\beta_j}q_j +(u-q_{J'}))=  \sum_{j\in J'}\beta_jq_j,
\]
and hence $h=h'$.
\end{proof}
The next result shows that $\phi$ is continuous on $\partial\overline{V}^h$. 
\begin{theorem}\label{continuity} The map $\phi\colon \overline{V}^h\to B_1^*$ is continuous. 
\end{theorem}
\begin{proof}
Clearly $\phi$ is continuous on $V$. Suppose  $(y^n)$ is a sequence in $V$ with $h_{y^n}\to h\in\partial\overline{V}^h$.  We wish to show that $\phi(y^n)\to\phi(h)$. Let $(\phi(y^{n_k}))$ be a subsequence. We will show that it has a subsequence which converges to $\phi(h)$. 

As $h$ is a horofunction, we know that $r^n=\|y^{n_k}\|_u\to\infty$ by Lemma \ref{Rieffel}. For each $k$ there exists a Jordan frame $q_1^{n_k},\ldots,q_r^{n_k}$ in $V$ and $\lambda_1^{n_k},\ldots,\lambda_r^{n_k}\in\mathbb{R}$ such that 
\[
y^{n_k} = \sum_{i=1}^r \lambda_i^{n_k}q_i^{n_k}.
\]
By taking a subsequence we may assume that there exists $I_+\subseteq \{1,\ldots,r\}$ and $1\leq s\leq r$ such that 
$r^{n_k} =\|y^{n_k}\|_u =|\lambda^{n_k}_s|$, $\lambda_i^{n_k}>0$ if and only if  $i\in I_+$, for all $k$. 

For each $k$, let $\beta^{n_k}_i = r^{n_k} -\lambda_i^{n_k}$ for $i \in I_+$, and $\beta^{n_k}_i = r^{n_k} +\lambda_i^{n_k}$ for $i \not\in I_+$. Note that $\beta^{n_k}\geq 0$ for all $i$ and $k$, and $\beta_s^{n_k}=0$ for all $k$. By taking a further subsequence we may assume that $\beta^{n_k}_i\to\beta_i\in [0,\infty]$  and $q^{n_k}_i\to q_i$ for all $i$. 
Let $I'=\{i\in I_+\colon \beta_i<\infty\}$ and $J'=\{j\not\in I_+\colon \beta_j<\infty\}$. Note that $s\in I'\cup J'$ and we can apply Proposition \ref{prop:3.1} to conclude that $h_{y^{n_k}}\to h'\in\partial\overline{V}^h$, where 
\[
h'(x) = \max\{\Lambda_{V(q_{I'})}( -U_{q_{I'}}x - \sum_{i\in I'}\beta_iq_i),\Lambda_{V(q_{J'})}( U_{q_{J'}}x - \sum_{j\in J'}\beta_jq_j)\}. 
\]
As $h_{y^{n_k}}\to h$, we know that $h=h'$ and hence $\delta(h,h')=0$. This implies that $p_I=q_{I'}$ and $p_J=q_{J'}$ by Theorem \ref{main2}. Moreover, 
\[
\sum_{i\in I}\alpha_ip_i +\sum_{j\in J}\alpha_jp_j = \sum_{i\in I'}\beta_iq_i+\sum_{j\in J'}\beta_jq_j.
\]
It follows that 
\[
\sum_{i\in I}\alpha_ip_i = U_{p_I}(\sum_{i\in I}\alpha_ip_i +\sum_{j\in J}\alpha_jp_j) = U_{q_{I'}}(\sum_{i\in I'}\beta_iq_i+\sum_{j\in J'}\beta_jq_j)= \sum_{i\in I'}\beta_iq_i
\]
and 
\[
\sum_{j\in J}\alpha_jp_j = U_{p_J}(\sum_{i\in I}\alpha_ip_i +\sum_{j\in J}\alpha_jp_j) = U_{q_{J'}}(\sum_{i\in I'}\beta_iq_i+\sum_{j\in J'}\beta_jq_j)= \sum_{j\in J'}\beta_jq_j,
\]
so that $\sum_{i\in I}e^{\alpha_i} p_i = \sum_{i\in I'}e^{\beta_i}q_i$ and $\sum_{j\in J}e^{\alpha_j} p_j = \sum_{j\in J'}e^{\beta_j}q_j$. We conclude that 
\begin{eqnarray*}
\lim_{k\to\infty} \phi(y^{n_k}) & =& \lim_{k\to\infty}  \frac{\sum_{i=1}^r (e^{-r^{n_k}+\lambda_i^{n_k}} -e^{-r^{n_k}-\lambda_i^{n_k}})q_i^{n_k}}{\sum_{i=1}^r (e^{-r^{n_k}+\lambda_i^{n_k}} +e^{-r^{n_k}-\lambda_i^{n_k}})}\\ 
& = &  
\frac{\sum_{i\in I'} e^{-\beta_i}q_i - \sum_{j\in J'}e^{-\beta_j}q_j}{\sum_{i\in I'} e^{-\beta_i} + \sum_{j\in J'}e^{-\beta_j}} \\
& =& \frac{\sum_{i\in I} e^{-\alpha_i}p_i - \sum_{j\in J}e^{-\alpha_j}p_j}{\sum_{i\in I} e^{-\alpha_i} + \sum_{j\in J}e^{-\alpha_j}} =\phi(h).
\end{eqnarray*}

From Lemmas \ref{homlem1} and \ref{homlem2} we know that $\phi$ maps $V$ into $\mathrm{int}\, B_1^*$ and $\partial \overline{V}^h$ into $\partial B_1^*$. So to complete the proof it remains to show that if $(h_n)$ in $\partial \overline{V}^h$ converges to $h\in \partial \overline{V}^h$, then $\phi(h_n)\to\phi(h)$.  Suppose  $h$ is given by (\ref{horofunction1}) and  for each $n$ the horofunction $h_n$ is given by 
\begin{equation}\label{h_n}
h_n(x)  = \max\left\{\Lambda_{V(q^n_{I_n})}\left(-U_{q^n_{I_n}}x - \sum_{i\in I_n} \beta^n_iq^n_i\right),  \Lambda_{V(q^n_{J_n})}\left(U_{q^n_{J_n}}x - \sum_{j\in J_n} 
\beta^n_j q^n_j\right)\right\} \mbox{\quad for $x\in V$,} 
\end{equation}
where $I_n,J_n\subseteq\{1,\ldots,r\}$ are disjoint, $I_n\cup J_n$ is nonempty, and  $\min\{\beta^n_k\colon k\in I_n\cup J_n\} =0$.

To prove the assertion we show that each subsequence of  $(\phi(h_n))$ has a convergent subsequence with limit $\phi(h)$. Let $(\phi(h_{n_k}))$ be a subsequence. 
By taking a subsequences we may assume that  
\begin{enumerate}[(1)]
\item There exist $I_0,J_0\subseteq\{1,\ldots,r\}$  disjoint with $I_0\cup J_0$ nonempty such that $I_{n_k} = I_0$ and $J_{n_k} = J_0$ for all $k$.  
\item $\beta^{n_k}_i\to\beta_i\in [0,\infty]$  and $q^{n_k}_i\to q_i$ for all $i\in I_0\cup J_0$. 
\item There exists $i^*\in I_0\cup J_0$ such that $\beta^{n_k}_{i^*}=0$ for all $k$. 
\end{enumerate}
Let $I' = \{i\in I_0\colon \beta_i<\infty\}$ and $J'=\{j\in J_0\colon \beta_j<\infty\}$, and note that $i^*\in I'\cup J'$. 
 
 Using Lemma \ref{lem:C}  we now show that $h_{n_k}\to h'$, where 
 \begin{equation}\label{h'}
 h'(x)=\max\left\{\Lambda_{V(q_{I'})}\left(-U_{q_{I'}}x - \sum_{i\in I'} \beta_iq_i\right),  \Lambda_{V(q_{J'})}\left(U_{q_{J'}}x - \sum_{j\in J'} \beta_j q_j\right)\right\}.
\end{equation}
Note that if  $I'$ is nonempty, then by Lemma \ref{lem:C} we have that 
\[
\lim_{k\to\infty} \Lambda_{V(q^{n_k}_{I_0})}\left(-U_{q^{n_k}_{I_0}}x - \sum_{i\in I_0} \beta^{n_k}_iq^{n_k}_i\right)= 
\Lambda_{V(q_{I'})}\left(-U_{q_{I'}}x - \sum_{i\in I'} \beta_iq_i\right), 
\] 
as $ U_{q^{n_k}_{I_0}}x \to U_{q_{I_0}}x$ by Lemma \ref{conv} and $U_{q_{I'}}(U_{q_{I_0}}x) = U_{q_{I'}}x$ by \cite[Proposition 2.26]{AS1}. Likewise if  $J'$ is nonempty, we have that 
\[
\lim_{k\to\infty} \Lambda_{V(q^{n_k}_{J_0})}\left(U_{q^{n_k}_{J_0}}x - \sum_{j\in J_0} \beta^{n_k}_jq^{n_k}_j\right)= 
\Lambda_{V(q_{J'})}\left(U_{q_{J'}}x - \sum_{j\in J'} \beta_jq_j\right).
\] 
Thus, if $I'$ and $J'$ are  both nonempty (\ref{h'}) holds. 

Now suppose that $I'$ is empty, so $J'$ is nonempty. As $-x\leq \|x\|_u u$, we get that 
\[
-U_{q_{I_0}^{n_k}}x\leq \|x\|_u U_{q_{I_0}^{n_k}}u =  \|x\|_u q_{I_0}^{n_k}. 
\]
This implies that $-U_{q_{I_0}^{n_k}}x - \sum_{i\in I_0}\beta_i^{n_k}q_i^{n_k}\leq \sum_{i\in I_0} (\|x\|_u -\beta_i^{n_k})q_i^{n_k}$, and hence 
\[
\Lambda_{V(q^{n_k}_{I_0})}\left(-U_{q^{n_k}_{I_0}}x - \sum_{i\in I_0} \beta^{n_k}_iq^{n_k}_i\right)\leq \max_{i\in I_0}(\|x\|_u -\beta_i^{n_k})\to -\infty.
\]
On the other hand, $h_{n_k}(x)\geq -\|x\|_u$ for all $k$. Thus, for all $k$ sufficiently large, we have that 
\[
h_{n_k}(x)  =   \Lambda_{V(q^{n_k}_{J_0})}\left(U_{q^{n_k}_{J_0}}x - \sum_{j\in J_0} \beta_j^{n_k}q_j^{n_k}\right),
\] 
which implies  that (\ref{h'}) holds if $I'$ is empty. In the same way it can be shown that (\ref{h'}) holds if $J'$ is empty. 

As $h_n\to h$, we know that $h'=h$, so $\delta(h,h')=0$.  It follows from  Theorem \ref{main2} that 
$p_I=q_{I'}$, $p_J=q_{J'}$, and $\sum_{i\in I}\alpha_ip_i +\sum_{j\in J}\alpha_jp_j = \sum_{i\in I'}\beta_iq_i+\sum_{j\in J'}\beta_jq_j$. This implies that   
\[
\sum_{i\in I}\alpha_ip_i = \sum_{i\in I'}\beta_iq_i\mbox{\quad 
and
\quad }
\sum_{j\in J}\alpha_jp_j = \sum_{j\in J'}\beta_jq_j,
\]
so that $\sum_{i\in I}e^{\alpha_i} p_i = \sum_{i\in I'}e^{\beta_i}q_i$ and $\sum_{j\in J}e^{\alpha_j} p_j = \sum_{j\in J'}e^{\beta_j}q_j$. Thus,  
\begin{eqnarray*}
\lim_{k\to\infty} \phi(h_{n_k}) & =& \lim_{k\to\infty}  \frac{\sum_{i\in I_0} e^{-\beta_i^{n_k}}q^{n_k}_i -\sum_{j\in J_0} e^{-\beta_j^{n_k}}q^{n_k}_j}{\sum_{i\in I_0} e^{-\beta_i^{n_k}}+\sum_{j\in J_0} e^{-\beta_j^{n_k}}}\\ 
& = &  \frac{\sum_{i\in I'} e^{-\beta_i}q_i - \sum_{j\in J'}e^{-\beta_j}q_j}{\sum_{i\in I'} e^{-\beta_i} + \sum_{j\in J'}e^{-\beta_j}} \\
& =& \frac{\sum_{i\in I} e^{-\alpha_i}p_i - \sum_{j\in J}e^{-\alpha_j}p_j}{\sum_{i\in I} e^{-\alpha_i} + \sum_{j\in J}e^{-\alpha_j}} =\phi(h),
\end{eqnarray*}
which completes the proof.
\end{proof}
\begin{theorem}\label{surjective} The map $\phi\colon \overline{V}^h\to B_1^*$ is onto.
\end{theorem}
\begin{proof}
From Lemma \ref{lem:ontoint} we know that $\phi(V) = \mathrm{int}\, B_1^*$.  Let $z\in\partial B_1^*$. As $B_1^*$ is the disjoint union of the relative interiors of its faces, see \cite[Theorem 18.2]{Rock}, we know that there exist orthogonal idempotents $p_I$ and $p_J$ such that $z\in\mathrm{relint} F_{p_I,p_J}$. So we can write \[
z =\sum_{i\in I} \lambda_ip_i-\sum_{j\in J}\lambda_jp_j,
\]
where $p_I =\sum_{i\in I}p_i$, $q_J =\sum_{j\in J}q_j$, $0<\lambda_k\leq 1$ for all $k\in I\cup J$, and $\sum_{k\in I\cup J}\lambda_k =1$. 

Define $\mu_k = -\log \lambda_k$ for $k\in I\cup J$.  So, $\mu_k\geq 0$ and let $\mu^*=\min\{\mu_k\colon k\in I\cup J\}$. Set $\alpha_k = \mu_k-\mu^*\geq 0$ and note that $\min\{\alpha_k\colon k\in I\cup J\}=0$.  

Then $h$, given by 
\[
h(x)  = \max\left\{\Lambda_{V(p_I)}\left(-U_{p_I}x - \sum_{i\in I} \alpha_ip_i\right),  \Lambda_{V(p_J)}\left(U_{p_J}x - \sum_{j\in J} \alpha_j p_j\right)\right\} 
\]
for $x\in V$, is a horofunction by Proposition \ref{prop:2}. Moreover, 
\begin{eqnarray*}
\phi(h) & = &  \frac{1}{\sum_{i\in I} e^{-\alpha_i} +\sum_{j\in J} e^{-\alpha_j}}\left(\sum_{i\in I} e^{-\alpha_i}p_i -\sum_{j\in J} e^{-\alpha_j}p_j\right)\\
  & = &  \frac{1}{\sum_{i\in I} e^{-\mu_i} +\sum_{j\in J} e^{-\mu_j}}\left(\sum_{i\in I} e^{-\mu_i}p_i -\sum_{j\in J} e^{-\mu_j}p_j\right)\\
  & = &  \frac{1}{\sum_{i\in I} \lambda_i +\sum_{j\in J} \lambda_j}\left(\sum_{i\in I} \lambda_ip_i -\sum_{j\in J} \lambda_jp_j\right),\\
\end{eqnarray*}
and hence $\phi(h)=z$, which  completes the proof.
\end{proof}
The proof of Theorem \ref{hom} is now straightforward. 
\begin{proof}[Proof of Theorem \ref{hom}] 
It follows from Theorems  \ref{continuity} and \ref{surjective} and Corollary \ref{injective} that $\phi\colon\overline{V}^h\to B_1^*$ is a continuous bijection. As $\overline{V}^h$ is compact and $B_1^*$ is Hausdorff, we conclude that $\phi$ is a homeomorphism. It follows from Lemma \ref{homlem2} that $\phi$ maps parts onto the relative interior of a boundary face of $B_1^*$.
\end{proof}

\begin{remark}
It is interesting to note that a similar idea can be used to show that the horofunction compactification of a finite dimensional normed space $(V,\|\cdot\|)$ with a smooth, strictly convex, norm is homeomorphic to the closed dual unit ball. Indeed, in that case the horofunctions are given by $h\colon z\mapsto -x^*(z)$, where $x^*\in V^*$ has norm 1, see for example \cite[Lemma 5.3]{Gu2}.  Moreover for $(y^n)$ in $V$ we have that $h_{y^n}\to h$ if and only if $y^n/\|y^n\|\to x$ and $\|y^n\|\to\infty$. 

In this case we define a map $\psi\colon \overline{V}^h\to B_1^*$ as follows. For $x\in V$ with $x\neq 0$, let 
\[
\psi(x)  = -\left(\frac{e^{\|x\|} - e^{-\|x\|}}{e^{\|x\|}+e^{-\|x\|}}\right) x^*,
\]
where $x^* \in  V^*$ is the unique functional with $x^*(x) =\|x\|$ and $\|x^*\| =1$, and let $\psi(0)=0$. 
For $h\in\partial \overline{V}^h$ with $h\colon z\mapsto -x^*(z)$ let 
\[
\psi(h) = -x^*.
\]
It is straightforward to check that $\psi$ is a bijection from $\overline{V}^h$  onto $B_1^*$, and $\psi$ is continuous on $\mathrm{int}\,B_1^*$. To show continuity on $\partial \overline{V}^h$, we assume, by way of contradiction, that $(h_n)$ is a sequence of horofunctions with $h_n\to h$ and $h_n (z) = -x^*_n(z)$ for all $z\in V$, and  there exists a neighbourhood $U$ of $\psi(h)$ in $B_1^*$ such that $\psi(h_n)\not\in U$ for all $n$. Then for each $z^*\in\partial B_1^*$ with $z^*\not\in U$ we have that $z^*(x)<1$. So, by compactness,  $\delta = \max\{ 1 - z^*(x)\colon z^*\in\partial B_1^*\setminus U\}>0$. It now follows that 
\[
h_n(x) - h(x) = -x^*_n(x) +x^*(x) = 1- x^*_n(x)\geq \delta>0 
\]
for all $n$, which contradicts $h_n\to h$. This shows  that $\psi$ is continuous bijection, and hence a homeomorphism, as $\overline{V}^h$ is compact and $B_1^*$ is Hausdorff.

More generally, one can consider product spaces $V=\prod_{i=1}^r V_i$ with norm  $\|x\|_V = \max_{i=1}^r \|v_i\|_i$, where each  $(V_i,\|\cdot\|_i)$ is a finite dimensional normed space with a smooth, strictly convex, norm. In that case we have by  \cite[Theorem 2.10]{L1} that the horofunctions of $V$ are given by 
\begin{equation}\label{eq:5.1}
h(v) = \max_{j\in J} (h_{\xi_j^*}(v_j)-\alpha_j),
\end{equation}
where $J\subseteq\{1,\ldots,r\}$ nonempty, $\min_{j\in J}\alpha_j =0$, $\xi_j^*\in V_j^*$ with $\|\xi_j^*\| =1$, and $h_{\xi_j^*}(v_j) = -\xi_j^*(v_j)$. 
One can use similar ideas as the ones in Section 3 to  show that the horofunction compactification is homeomorphic to the closed unit dual ball of $V$.  Indeed, one can define a  map $\phi_V\colon \overline{V}^h\to B_1^*$ as follows. For $v\in V$ let  
\[
\phi_V(v) = \frac{1}{\sum_{i=1}^r e^{\|v_i\|_i} + e^{-\|v_i\|_i}} \left(\sum_{i=1}^r (e^{\|v_i\|_i} - e^{-\|v_i\|_i})p(v^*_i)\right),
\]
and $\phi_V(0)=0$. Here $p(v^*_i) = (0,\ldots,0,v^*_i,0,\ldots,0)$ and $v^*_i$ is the unique functional such that  $v^*_i(v_i) =\|v_i\|_i$ and $\|v^*_i\|_i=1$, if $v_i\neq 0$, and we set $p(v_i^*)=0$, if $v_i=0$.  
For a horofunction $h$ given by (\ref{eq:5.1}) we define 
\[
\phi_V(h) = \frac{1}{\sum_{j\in J} e^{-\alpha_j}}\left( \sum_{j\in J} e^{-\alpha_j}p(\xi_j^*)\right).
\]
Following the same line of reasoning as in Section 3 one can prove that $\phi_V$ is a homeomorphism. 
\end{remark}
The connection between the global topology of the horofunction compactification and the dual unit ball  seems hard to establish for general finite dimensional normed spaces, and might not even hold. In the settings discussed in this paper all horofunctions are Busemann points, but there are normed spaces with horofunctions that are not Busemann, see \cite{Wa2}. It could well be the case that the horofunction compactification of these spaces is not homeomorphic to the closed unit dual ball, but no counter example is known at present.

\section{Hilbert geometries} 
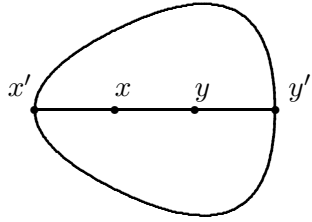
\begin{figure}[h]
\begin{center}
\begin{picture}(30, 70)  
\thicklines
  \closecurve(0,0, 60,30, 0,60)
  \put(-40,35){$x'$}
  \put(-30,30){\circle*{3}}
\put(0,35){$x$}
  \put(0,30){\circle*{3}}
\put(30,35){$y$}
  \put(30,30){\circle*{3}}
\put(65,35){$y'$}
  \put(60,30){\circle*{3}}
  \put(-30,30){\line(1,0){90}}
   \end{picture}
   \caption{Hilbert distance}\label{fig:1}
  \end{center}  
\end{figure}

In this section we study global topology and geometry of the horofunction compactification  of certain Hilbert geometries. 
Recall that the Hilbert distance  is defined as follows.  Let $A$ be a real finite dimensional affine space. Consider a bounded, open, convex set  $\Omega\subseteq A$. For $x,y\in \Omega$, let $\ell_{xy}$ denote the straight-line through $x$ and $y$ in $A$, and denote the points of intersection of $\ell_{xy}$ and $\partial\Omega$ by $x'$ and $y'$, where $x$ is between $x'$ and $y$, and $y$ is between $x$ and $y'$, as in Figure \ref{fig:1}. 

On $\Omega$ the  {\em Hilbert distance}  is defined by 
\begin{equation}\label{eq:1.1}
\rho_H(x,y) = \log \Big{(}\frac{|x'-y|}{|x'-x|} \frac{|y'-x|}{|y'-y|}\Big{)}
\end{equation}
for all $x\neq y$ in $\Omega$, and $\rho_H(x,x) =0$ for all $x\in \Omega$.
The metric space $(\Omega, \rho_H)$ is called the  {\em Hilbert geometry} on $\Omega$. 

These metric spaces generalise Klein's model of hyperbolic space and have a Finsler structure, see \cite{Nu,PT}.  
In our analysis we will work with Birkhoff's version of the Hilbert metric, which is defined on a cone in an order-unit space in terms of its partial ordering. This provides a convenient way to work with the Hilbert distance and its Finsler structure. In the next subsection we will recall the basic concepts involved in our analysis. Throughout we will follow the terminology used in \cite[Chapter 2]{LNBook}, which contains a detailed discussion of  Hilbert geometries  and some their applications. 
We refer the reader to \cite{PT} for a comprehensive account of the theory of Hilbert geometries. 

\subsection{Preliminaries and Finsler structure} 
Let $(V,V_+,u)$  be a finite dimensional order-unit space. So $V_+$ is a closed cone in $V$ with $u\in\mathrm{int}\,V_+$. Recall that the cone $V_+$ induces a partial ordering on $V$ by $x\leq y$  if $y-x\in V_+$, see Section 4.1. For $x\in V$ and $y\in V_+$, we say that $y$ {\em dominates} $x$ if there exist $\alpha,\beta\in\mathbb{R}$ such that $\alpha y\leq x\leq \beta y$. In that case, we write  
\[
M(x/y) =\inf\{\beta\in\mathbb{R}\colon x\leq \beta y\}
\]
and
\[
m(x/y) = \sup\{\alpha\in\mathbb{R}\colon \alpha y\leq x\}.
\]

By the Hahn-Banach theorem, $x\leq y$ if and only if $\psi(x)\leq\psi (y)$ for all $\psi\in V^*_+=\{\phi\in V^*\colon \phi\mbox{ positive}\}$, which is equivalent to $\psi(x)\leq\psi(y)$ for all $\psi\in S(V)$. Using this fact, it easy to verify that for each $x\in V$ and $y\in\mathrm{int}\,V_+$ we have  
\[
M(x/y) = \sup_{\psi\in S(V)}\frac{\psi(x)}{\psi(y)}\mbox{\quad and\quad} m(x/y) =\inf_{\psi\in S(V)}\frac{\psi(x)}{\psi(y)}.
\]

We also note that if $A\in \mathrm{GL}(V)$ is a linear automorphism of $V_+$, i.e., $A(V_+)= V_+$, then $x\leq \beta y$ if, and only if, $Ax\leq \beta Ay$. It follows that $M(Ax/Ay)=M(x/y)$ and $m(x/y)= m(Ax/Ay)$. 

If $w\in\mathrm{int}\, V_+$, then $w$ dominates each $x\in V$, and we define 
\[
|x|_w = M(x/w)- m(x/w).
\]
One can verify that $|\cdot|_w$ is a semi-norm on $V$, see \cite[Lemma A.1.1]{LNBook}, and a genuine norm on the quotient space $V/\mathbb{R}w$, as $|x|_w=0$ if and only if $x=\lambda w$ for some $\lambda\in\mathbb{R}$. 

Clearly, if $x,y\in V$ are such that $y=0$ and $y$ dominates $x$,  then $x=0$, as $V_+$ is a cone. On the other hand, if $y\in V_+\setminus\{0\}$, and $y$ dominates $x$, then $M(x/y)\geq m(x/y)$. The domination relation yields an equivalence relation on $V_+$ by $x\sim y$ if $y$ dominates $x$ and $x$ dominates $y$. The equivalence classes are called the {\em parts} of $V_+$. As $V_+$ is closed, one can check that  $\{0\}$ and $\mathrm{int}\,V_+$ are  parts of $V_+$. The parts of a finite dimensional cone are closely related to its faces. Indeed, if $V_+$ is the cone of a finite dimensional order-unit space, then it can be shown that the parts correspond to the relative interiors of the faces of $V_+$, see \cite[Lemma 1.2.2]{LNBook}. Recall that a {\em face} of a convex set $S\subseteq V$ is a subset  $F$ of $S$ with the property that if $x,y\in S$ and $\lambda x+ (1-\lambda)y\in F$ for some $0<\lambda <1$, then $x,y\in F$. 

It is easy to verify that if $x,y\in V_+\setminus\{0\}$, then $x\sim y$ if, and only if, there exist $0<\alpha\leq \beta$ such that $\alpha y\leq  x\leq \beta y$. Furthermore, if $x\sim y$, then 
\begin{equation}\label{eq:2.1}
m(x/y) = \sup\{\alpha>0\colon y\leq \alpha^{-1}x\}= M(y/x)^{-1}.
\end{equation}
{\em Birkhoff's version of the  Hilbert distance}  on $V_+$ is defined as follows:
\begin{equation}\label{eq:2.2}
d_H(x,y) = \log\Big{(}\frac{M(x/y)}{m(x/y)}\Big{)} =  \log M(x/y) +\log M(y/x) 
\end{equation}
for all $x\sim y$ with $y\neq 0$, $d_H(0,0)=0$, and $d_H(x,y)=\infty$ otherwise.  

Note that $d_H(\lambda x,\mu y) = d_H(x,y)$ for all $x,y\in V_+$ and $\lambda,\mu>0$, so $d_H$ is not a distance on $V_+$. It is, however, a distance between pairs of rays in each part of $V_+$.  In particular, if $\phi\colon V\to\mathbb{R}$ is a linear functional such that $\phi(x)>0$ for all $x\in V_+\setminus\{0\}$, then $d_H$ is a distance on 
\[
\Omega_V=\{x\in \mathrm{int}\, V_+\colon \phi(x)=1\},
\]
which is a (relatively) open, bounded, convex set, see \cite[Lemma 1.2.4]{LNBook}. Moreover, the following holds, see \cite[Proposition 2.1.1 and Theorem 2.1.2]{LNBook}. 
\begin{theorem} $(\Omega_V,d_H)$ is a metric space and $d_H=\rho_H$ on $\Omega_V$. 
\end{theorem}
It is worth noting that any Hilbert geometry can be realised as $(\Omega_V,d_H)$ for some order-unit space $V$ and strictly positive linear functional $\phi$.

A Hilbert geometry  $(\Omega_V,d_H)$ has a Finsler structure, see \cite{Nu}. Indeed, if one  defines the length of a piecewise $C^1$-smooth path $\gamma\colon [0,1]\to \Omega_V$  by 
\[
L(\gamma) = \int^1_0 |\gamma'(t)|_{\gamma(t)}\mathrm{d}t,
\]
then 
\begin{equation}\label{HilbertFinsler}
d_H(x,y) = \inf_\gamma L(\gamma),
\end{equation}
where the infimum is taken over all piecewise $C^1$-smooth paths in $\Omega_V$ with $\gamma(0)=x$ and $\gamma(1)=y$. 

So for Hilbert geometries Problem \ref{keyproblem} can be formulated more explicitly as follows. 
\begin{problem}\label{HilbertProblem} Let $(V,V_+,u)$  be a finite dimensional  order-unit space and $\phi\colon V\to\mathbb{R}$ be a linear functional with $\phi(x)>0$ for all $x\in V_+\setminus\{0\}$ and $\phi(u)=1$. For which Hilbert geometries $(\Omega_V,d_H)$ does there exists a homeomorphism from the horofunction compactification $ \overline{\Omega}_V^h$ with  basepoint $u$ onto the closed dual unit ball $B_1^*$ of $|\cdot|_u$ on $V/\mathbb{R}u$, which maps each part of the horofunction boundary onto the relative interior of a boundary face of $B_1^*$?
\end{problem}
It should be noted that in the case of Hilbert geometries the unit ball $\{x\in V/\mathbb{R}w\colon |x|_w\leq 1\}$ in the tangent space at $w\in\Omega_V$ may have a different facial structure for different $w$.  This phenomenon appears frequently in the  case where $\Omega_V$ is a polytope. 

This problem, however, does not arise in the spaces we analyse here. Indeed, we will consider order-unit spaces $(V,V_+,u)$, where $V$ is a Euclidean Jordan algebra of rank $r$, $V_+$ is the cone of squares, and $u$ is the algebraic unit.  So $\mathrm{int}\, V_+$ is a symmetric cone and $\mathrm{Isom}(\Omega_V)$ acts transitively on $\Omega_V$. Throughout we will take $\phi\colon V\to\mathbb{R}$ with $\phi(x) = \frac{1}{r}\mathrm{tr}(x)$, which is a state and
\[
\Omega_V=\{x\in\mathrm{int}\, V_+\colon \phi(x)=1\}=\{x\in \mathrm{int}\, V_+\colon \mathrm{tr}(x) = r\}.
\]
In that case we call $(\Omega_V,d_H)$ {a \em symmetric} Hilbert geometry.  A prime example is  
\[
\Omega_V =\{ A\in\mathrm{Herm}_n(\mathbb{C})\colon \mathrm{tr}(A) = n\mbox{ and }A\mbox{ positive definite}\}.
\]

In a symmetric Hilbert geometry  the distance can be expressed in terms of the spectrum. Indeed, we know that for $x\in V$ invertible, the quadratic representation $U_x\colon V\to V$ is a linear automorphism of $V_+$, see \cite[Proposition III.2.2]{FK}. Moreover, $U_x^{-1} = U_{x^{-1}}$ and $U_{x^{-1/2}}x=u$.  Furthermore, for $x\in V$ we have that 
\[
M(x/u) =\inf\{\lambda\colon x\leq \lambda u\} =  \max\sigma(x)\mbox{\quad and\quad} 
m(x/u) =\sup\{\lambda\colon  \lambda u\leq x\} =  \min\sigma(x),
\] 
so that $|x|_u =  \max\sigma(x) - \min\sigma(x)$. 
Also for $x,y\in\mathrm{int}\, V_+$ we have that 
\[
\log M(x/y) = \log M(U_{y^{-1/2}}x/u) = \log \max\sigma(U_{y^{-1/2}}x) =  \max\sigma(\log U_{y^{-1/2}}x)
\]
and 
\[
\log M(y/x) = \log m(x/y)^{-1}= -\log m(U_{y^{-1/2}}x/u) = -\min\sigma(\log U_{y^{-1/2}}x).
\]
It follows  that 
\[
d_H(x,y) =\log M(x/y)+\log M(y/x) =  |\log U_{y^{-1/2}}x|_u=\mathrm{diam\,}\sigma(\log U_{y^-1/2}x)\mbox{\quad for all }x,y\in\mathrm{int}\, V_+.
\]
Moreover, for each $w\in \Omega_V$ we have that 
\[
|x|_w = M(x/w)- m(x/w) = M(U_{w^{-1/2}}x/u)-m(U_{w^{-1/2}}x/u) = |U_{w^{-1/2}}x|_u\mbox{\quad for all }x\in V,
\]
which shows that the facial structure of the unit ball in each tangent space is identical, as $U_{w^{-1/2}}$ is an invertible linear map. 

By using the Jordan algebra structure there is a direct way to show that a symmetric Hilbert geometry has a Finsler structure. 
\begin{proposition}
If $(\Omega_V,d_H)$ is a symmetric Hilbert geometry, then for each $x,y\in \Omega_V$ we have that $d_H(x,y)=\inf L(\gamma)$, where the infimum is taken over all piecewise $C^1$-smooth paths $\gamma\colon [0,1]\to \Omega_V$ with $\gamma(0) =x$ and $\gamma(1)=y$, and 
\[
L(\gamma)=\int_0^1 |\gamma'(t)|_{\gamma(t)}\mathrm{d}t.
\]
\end{proposition}
\begin{proof}
Let $\gamma\colon [0,1]\to \Omega_V$ be a piecewise $C^1$-path with $\gamma(0)=x$ and $\gamma(1)=y$. We have  
\begin{eqnarray*}
d_H(x,y) & = &\log M(y/x)  -\log m(y/x)\\
 & = &  \max_{\psi\in S(V)} \log \frac{\psi(y)}{\psi(x)} - \min_{\psi\in S(V)} \log \frac{\psi(y)}{\psi(x)}\\
& = & \max_{\psi\in S(V)}\int_0^1 \frac{\mathrm{d}}{\mathrm{d}t}\log \psi(\gamma(t)) \mathrm{d}t
- \min_{\psi\in S(V)}\int_0^1 \frac{\mathrm{d}}{\mathrm{d}t}\log \psi(\gamma(t)) \mathrm{d}t\\
& = & \max_{\psi\in S(V)}\int_0^1 \frac{ \psi(\gamma'(t))}{ \psi(\gamma(t))} \mathrm{d}t
-\min_{\psi\in S(V)}\int_0^1 \frac{ \psi(\gamma'(t))}{ \psi(\gamma(t))} \mathrm{d}t\\
&\leq & \int_0^1 \max_{\psi\in S(V)}\frac{ \psi(\gamma'(t))}{ \psi(\gamma(t))} \mathrm{d}t
-\int_0^1 \min_{\psi\in S(V)} \frac{\psi(\gamma'(t))}{\psi(\gamma(t))} \mathrm{d}t\\
& = & \int_0^1 M(\gamma'(t)/\gamma(t))-m(\gamma'(t)/\gamma(t))\mathrm{d}t\\
& = & \int_0^1 |\gamma'(t)|_{\gamma(t)}\mathrm{d}t.
\end{eqnarray*}

Now let $x,y\in \Omega_V$ and consider the $C^1$-smooth path $\sigma$ in $C^\circ$ given by, 
\[
\sigma(t) = U_{x^{1/2}}(U_{x^{-1/2}}y)^t\mbox{\quad for  }0\leq t\leq 1.
\] 
Note that $\sigma(0)=U_{x^{1/2}}u= x$ and $\sigma(1) = y$. Define 
\[
\mu(t) =\frac{\sigma(t)}{\phi(\sigma(t))}\mbox{\quad for all }0\leq t\leq 1.
\] 
So, $\mu$ is a $C^1$-smooth path connecting $x$ and $y$ in $\Omega_V$. 
A direct calculation gives 
\[
\mu'(t) = \frac{\sigma'(t)}{\phi(\sigma(t))} - 
\frac{\phi(\sigma'(t))}{\phi(\sigma(t))^2}\sigma(t) 
\mbox{\quad for  } 0\leq t\leq 1.\]

We also have that $U_{\mu(t)^{-1/2}} =\phi(\sigma(t))U_{\sigma(t)^{-1/2}}$ for $0\leq t\leq 1$, 
which implies 
\begin{equation}\label{eq:4.4.1} 
U_{\mu(t)^{-1/2}}\mu'(t) = U_{\sigma(t)^{-1/2}}\sigma'(t) - \frac{\phi(\sigma'(t))}{\phi(\sigma(t))} u.
\end{equation}
Furthermore
\[
\sigma'(t) = U_{x^{1/2}}((U_{x^{-1/2}}y)^t\log(U_{x^{-1/2}}y))\mbox{\quad for  $0\leq t\leq 1$.
}\]

Write $z=U_{x^{-1/2}}y$ and let $z=\sum_{i=1}^r \lambda_ip_i$ be the spectral decomposition of $z$. As $z\in\mathrm{int}\, V_+$, we have that  $z^t = \sum_{i=1}^r\lambda_i^t p_i$ and $\log z =\sum_{i=1}^r (\log \lambda_i)p_i$, and hence  
\[
z^t\log z = \sum_{i=1}^r (\lambda_i^t\log\lambda_i)p_i.
\mbox{\quad and\quad }
U_{z^{-t/2}}(z^t\log z) =\log z.
\]

From  (\ref{eq:4.4.1}) we get that 
\begin{eqnarray*}
M(\mu'(t)/\mu(t)) - m(\mu'(t)/\mu(t)) &=& M(U_{\mu(t)^{-1/2}}\mu'(t)/u) - m(U_{\mu(t)^{-1/2}}\mu'(t)/u)\\
& = &M(U_{\sigma(t)^{-1/2}}\sigma'(t)/u) - m(U_{\sigma(t)^{-1/2}}\sigma'(t)/u).
\end{eqnarray*}
It follows that 
\begin{eqnarray*}
 M(\mu'(t)/\mu(t)) - m(\mu'(t)/\mu(t))   & = &M(\sigma'(t)/\sigma(t)) - m(\sigma'(t)/\sigma(t))\\
 & =&  M(U_{x^{-1/2}}\sigma'(t)/U_{x^{-1/2}}\sigma(t)) - m(U_{x^{-1/2}}\sigma'(t)/U_{x^{-1/2}}\sigma(t))  \\
&=& M(z^t\log z/z^t) -  m(z^t\log z/z^t)  \\
&=&  M(\log z/u) -  m(\log z/u)  \\
&= &\log M(U_{x^{-1/2}}y/u) -  \log m(U_{x^{-1/2}}y/u)  \\
 &=& \log M(y/x)-\log m(y/x).\\
\end{eqnarray*}

We conclude that 
\[
L(\mu) =\int_0^1\log  M(y/x) -\log m(y/x)\mathrm{d}t = d_H(x,y),
\] 
which completes the proof.
\end{proof}

\subsection{Horofunctions of symmetric Hilbert geometries} 
The main objective is to confirm Problem \ref{HilbertProblem} for symmetric Hilbert geometries. To describe the homeomorphism, we recall the description of the horofunction compactification of symmetric Hilbert geometries given in \cite[Theorem 5.6]{LLNW}. 
\begin{theorem}\label{Hilberthoro}
The horofunctions of a symmetric Hilbert geometry $(\Omega_V,d_H)$ are precisely the functions $h\colon \Omega_V\to\mathbb{R}$ of the form 
\begin{equation}\label{Hilberthf}
h(x) = \log M(y/x)+\log M(z/x^{-1})\mbox{\quad for } x\in\Omega_V,
\end{equation} 
where $y,z\in\partial V_+$ are such  that $\|y\|_u=\|z\|_u =1$ and $(y|z)=0$.
\end{theorem}
It follows from the proof of \cite[Theorem 5.6]{LLNW} that all horofunctions are in fact Busemann points. Indeed, if $y$ and $z$ have spectral decompositions 
\[
y = \sum_{i \in I} \lambda_ip_i\mbox{\quad and\quad }z =\sum_{j\in J}\mu_jp_j,
\]
where $I,J\subset\{1,\ldots,r\}$ are nonempty and disjoint, and $p_1,\ldots,p_r$ is a Jordan frame, then the sequence $(y_n)\in \mathrm{int}\, V_+$ given by, 
\[
y_n = \sum_{i \in I } \lambda_ip_i+\sum_{j\in J}\frac{1}{n^2\mu_j}p_j + \sum_{k\not\in I\cup J} \frac{1}{n}p_k,
\]
has the property that $y_n\to y$, $y_n^{-1}/\|y_n^{-1}\|_u\to z$ and $h_{y_n}\to h$, where $h$ is as in (\ref{Hilberthf}).  Note that if we let $v_n = y_n/\phi(y_n)\in\Omega_V$, then $h_{v_n}(z)=h_{y_n}(z)$ for all $z\in \Omega_V$, so $h_{v_n}\to h$. 

Also note that  for $n,m\geq 1$, 
\[
U_{y_n^{-1/2}}y_m = \sum_{i \in I } p_i+\sum_{j\in J}\frac{n^2}{m^2}p_j + \sum_{k\not\in I\cup J} \frac{n}{m}p_k.
\]
This implies that for each $n\geq m\geq 1$,  
\[
M(y_m/y_n) = M(U_{y_n^{-1/2}}y_m/u) = \|U_{y_n^{-1/2}}y_m\|_u = n^2/m^2,
\]
so that $\log M(y_m/y_n) = 2\log n -2\log m$.  Moreover, $\log M(y_n/y_m) =\log 1 =0$ for all $n\geq m\geq 1$. 
It follows that 
\[
d_H(v_n,v_m)+d_H(v_m,v_1)= d_H(y_n,y_m)+d_H(y_m,y_1) = d_H(y_n,y_1) = d_H(v_n,v_1)
\]
for all $n\geq m\geq 1$. Thus, $(v_n)$ is an almost geodesic sequence in $\Omega_V$, and hence each horofunction in 
$\overline\Omega_V^h$ is a Busemann point. 

To identify the parts and describe the detour distance we need the following general lemma. 
\begin{lemma}\label{lem:5.1}
Let $(V,V_+,u)$ be a finite dimensional order-unit space. If $v\in\partial V_+\setminus\{0\}$ and $w_n\in\mathrm{int}\, V_+$ with $w_{n+1}\leq w_n$ for all $n\geq 1$ and $w_n\to w\in\partial V_+\setminus\{0\}$, then 
\[
\lim_{n\to\infty} M(v/w_n) = \left[\begin{array}{ll} M(v/w) & \mbox{ if $w$ dominates $v$}\\ \infty & \mbox{ otherwise.}\end{array}\right. 
\]
\end{lemma}
\begin{proof}
Set $\lambda_n = M(v/w_n)$ for $n\geq 1$. Then for $n\geq m\geq 1$ we have that $0\leq \lambda_nw_n -v \leq \lambda_nw_m-v$. This implies that  $\lambda_m\leq \lambda_n$ for all $m\leq n$, and hence $(\lambda_n)$ is monotonically increasing. 

Now suppose that $\lambda = M(v/w)<\infty$, i.e., $w$ dominates $v$.  Then $0\leq \lambda w-v\leq \lambda w_n -v$, and hence $\lambda_n\leq \lambda$ for all $n$.  This implies that $\lambda_n\to \lambda^*\leq \lambda<\infty$. 
As $0\leq \lambda_n w_n -v$ for all $n$ and $V_+$ is closed, we know that $\lim_{n\to\infty} \lambda_nw_n -v = \lambda^* w-v\in V_+$. So $\lambda^*\geq \lambda$, and hence $\lambda^*= \lambda$.  We conclude that if $w$ dominates $v$, then $\lim_{n\to\infty} M(v/w_n) = M(v/w)$. 

On the other hand, if $w$ does not dominate $v$, then 
\begin{equation}\label{eq:notinV_+} 
\lambda w-v\not\in V_+\mbox{\quad  for all }\lambda\geq 0.
\end{equation}
Assume, by way of contradiction, that $(\lambda_n)$ is bounded. Then $\lambda_n\to \lambda^*<\infty$, since $(\lambda_n)$ is increasing,  and $\lambda_nw_n -v \to \lambda^* w -v\in V_+$, as $V_+$ is closed. This contradicts (\ref{eq:notinV_+}), and hence $\lambda_n = M(v/w_n)\to\infty$, if $w$ does not dominate $v$.
\end{proof}

Before we identify the parts in $\partial \overline{\Omega}^h_V$ and the detour distance, it is useful to recall the following fact: 
\[
M(x/y) = M(y^{-1}/x^{-1}) \mbox{\quad for all }x,y\in\mathrm{int}\,V_+,
\]
if $\mathrm{int}\, V_+$ is a symmetric cone, see \cite[Section 2.4]{LRW}.

\begin{proposition}\label{d_Hparts}
Let $(\Omega_V,d_H)$ be a symmetric Hilbert geometry and $h,h'\in\partial\overline{\Omega}^h_V$ with 
\[
h(x) = \log M(y/x)+\log M(z/x^{-1})
\]
and
\[
h'(x) = \log M(y'/x)+\log M(z'/x^{-1})
\]
for $x\in \Omega_V$. The following assertions hold: 
\begin{enumerate}[(i)] 
\item $h$ and $h'$ are in the same part if and only if $y\sim y'$ and $z\sim z'$. 
\item If $h$ and $h'$ are in the same part, then $\delta(h,h') = d_H(y,y') +d_H(z,z')$. 
\end{enumerate}
\end{proposition}
\begin{proof}
Consider the spectral decompositions: $y = \sum_{i\in I}\lambda_ip_i$, $z = \sum_{j\in J}\mu_jp_j$,  $y' = \sum_{i\in I'}\alpha_iq_i$, and $z' = \sum_{j\in J'}\beta_jq_j$. 
Set 
\[
y_n = \sum_{i\in I}\lambda_ip_i + \sum_{j\in J}\frac{1}{n^2\mu_j}p_j + \sum_{k\not\in I\cup J} \frac{1}{n} p_k
\]
and 
\[
w_n = \sum_{i\in I'}\alpha_iq_i + \sum_{j\in J'}\frac{1}{n^2\beta_j}q_j + \sum_{k\not\in I'\cup J'} \frac{1}{n} q_k.
\]
Then $h_{y_n}\to h$ and $h_{w_n}\to h'$ by the proof of \cite[Theorem 5.6]{LLNW}. 

For all $n\geq 1$ large we have that $ \|w_n\|_u =\|y'\|_u=1$, so that 
\[
d_H(w_n,u)  = \log M(w_n/u) +\log M(u/w_n) = \log\|w_n\|_u+\log M(w_n^{-1}/u) =  \log\|w_n^{-1}\|_u.
\]
 Now set $v_n  = w^{-1}_n/\|w_n^{-1}\|_u$ and note that by (\ref{detourcost}), 
 \begin{eqnarray*}
 H(h',h) & = & \lim_{n\to\infty} d_H(w_n,u) +h(w_n)\\
   & = & \lim_{n\to\infty}  \log\|w_n^{-1}\|_u + \log M(y/w_n) +\log M(z/w_n^{-1})\\
   & = & \lim_{n\to\infty} \log M(y/w_n) +\log M(z/v_n^{-1}).
 \end{eqnarray*}
 
Clearly $w_{n+1}\leq w_n$ and $w_n\to y'$. Also 
\[
w_n^{-1} = \sum_{i\in I'} \alpha_i^{-1}q_i +\sum_{j\in J'} n^2\beta_j q_j +\sum_{k\not\in I'\cup J'} nq_k.
\]
So, for all $n\geq 1$ large, we have that $\|w_n^{-1}\|_u =n^2$, as $\max_{j\in J}\beta_j = \|z'\|_u =1$. It follows that 
\[
v_n = \sum_{i\in I'} \frac{1}{n^2\alpha_i}q_i +\sum_{j\in J'}\beta_j q_j +\sum_{k\not\in I'\cup J'} \frac{1}{n}q_k
\]
for all $n\geq 1$ large. So, $v_{n+1}\leq v_n$ for all $n\geq 1$ large and $v_n\to z'$. 
It now follows from Lemma \ref{lem:5.1} that $H(h',h)=\infty$ if $y'$ does not dominate $y$, or, $z'$ does not dominate $z$.  Moreover, if $y'$ dominates $y$, and, $z'$ dominates $z$, then $H(h',h) = \log M(y/y') +\log M(z/z')$. 
 
Interchanging the roles between $h$ and $h'$ we find that $H(h,h') = \infty$  if $y$ does not dominate $y'$, or, $z$ does not dominate $z'$, and $H(h,h') = \log M(y'/y) +\log M(z'/z)$, otherwise. 
Thus, $\delta(h,h') = d_H(y,y')+d_H(z,z')$ if and only if $y\sim y'$ and $z\sim z'$, and $\delta(h,h')=\infty$ otherwise.  
  \end{proof}
The characterisation of the parts and the detour distance  is a more explicit description of the general one one given in \cite[Theorem 4.9]{LW} in the case of symmetric Hilbert geometries. 
\subsection{The homeomorphism} 
Let us now define a map $\phi_H\colon \overline{\Omega}_V^h\to B^*_1$, where $B^*_1$ is the unit ball of the dual norm of $|\cdot|_u$ on $V/\mathbb{R}u$.  For $x\in \Omega_V$ let 
\[
\phi_H(x) = \frac{x}{\mathrm{tr}(x)}-\frac{x^{-1}}{\mathrm{tr}(x^{-1})},
\]
and for $h\in\partial\overline{\Omega}_V^h$ given by (\ref{Hilberthf}) let 
\[
\phi_H(h) = \frac{y}{\mathrm{tr}(y)}-\frac{z}{\mathrm{tr}(z)}.
\]
We note that $\phi_H(h)$ is well-defined by Proposition \ref{d_Hparts}. 

We will prove the following theorem in the sequel.
\begin{theorem}\label{Hilberthom}
If $(\Omega_V,d_H)$ is a symmetric Hilbert geometry, then the map $\phi_H\colon  \overline{\Omega}_V^h\to B^*_1$ is a homeomorphism which maps each part of $\partial\overline{\Omega}_V^h$ onto the relative interior of  a boundary face of $B^*_1$. 
\end{theorem}

We first analyse the dual unit ball $B_1^*$ of $|\cdot|_u$ and its facial structure. The following fact, see also \cite[Section 2.2]{LRW}, will be useful. 
\begin{lemma}\label{lem:5.2} Given an order-unit space $(V,V_+,u)$, the norm $|\cdot|_u$ on $V/\mathbb{R}u$ coincides with the quotient norm of $2\|\cdot\|_u$ on $V/\mathbb{R}u$. 
\end{lemma}
\begin{proof}
Denote the quotient norm of $2\|\cdot\|_u$ on $V/\mathbb{R}u$ by $\|\cdot\|_q$. Then 
\begin{eqnarray*}
\|\overline{x}\|_q & = & 2\inf_{\mu\in\mathbb{R}}\|x-\mu u\|_u\\
  & = & 2\inf_{\mu\in\mathbb{R}} \max_{\phi\in S(V)} |\phi(x) -\mu| \\
  & = & 2\inf_{\mu\in\mathbb{R}} \max\{\max_{\phi\in S(V)} (\phi(x)) -\mu,  \max_{\phi\in S(V)}( -\phi(x))+\mu\}\\
& = & \max_{\phi\in S(V)}(\phi(x)) + \max_{\phi\in S(V)}( -\phi(x))\\
 & = & |\overline{x}|_u
\end{eqnarray*}
for all $\overline{x}\in V/\mathbb{R}u$, as $\inf_{\mu\in\mathbb{R}} \max\{a -\mu,b+\mu\} = (a+b)/2$ for all $a,b\in\mathbb{R}$.  
\end{proof}

Recall that in a Euclidean Jordan algebra $V$ each $x$ has a unique orthogonal decomposition $x= x^+-x^-$, where $x^+$ and $x^-$ are orthogonal elements in $V_+$, see \cite[Proposition 1.28]{AS1}.    Let 
\[
\mathbb{R}u^\perp =\{x\in V\colon (u|x) = 0\} = \{x\in V\colon \mathrm{tr}(x^+) = \mathrm{tr}(x^-)\}.
\]
It follows from Lemma \ref{lem:5.2} that 
\[
(V/\mathbb{R}u,|\cdot|_u)^* = (\mathbb{R}u^\perp,\frac{1}{2}\|\cdot\|_u^*).
\]
So the dual unit ball $B_1^*$ in $\mathbb{R}u^\perp$ is given by 
\[
B_1^* =2\mathrm{conv}(S(V)\cup -S(V))\cap \mathbb{R}u^\perp,
\]
see \cite[Theorem 1.19]{AS0}, and its (closed) boundary faces are precisely the nonempty sets of the form, 
\[
A_{p,q} = 2\mathrm{conv}\, ( (U_p(V)\cap S(V))\cup(U_q(V)\cap -S(V)))\cap \mathbb{R}u^\perp,
\]
where $p$ and $q$ are orthogonal idempotents, see \cite[Theorem 4.4]{ER}. 

To prove Theorem \ref{Hilberthom} we collect a number of preliminary results. 
\begin{lemma}\label{lem:5.3} For each $x\in\Omega_V$ we have that $\phi_H(x)\in\mathrm{int}\, B_1^*$, and  for each $h\in\partial\overline{\Omega}_V^h$ we have that $\phi_H(h)\in\partial B_1^*$. 
\end{lemma}
\begin{proof}
Let $x = \sum_{i=1}^r \lambda_i p_i\in\Omega_V$, so $\lambda_i>0$ for all $i$.  Note that $(u|\phi_H(x)) = 1 -1 =0$ and hence $\phi_H(x) \in \mathbb{R}u^\perp$. Given $-u\leq z\leq u$, we have the Peirce decomposition of $z$ with respect to the frame $p_1,\ldots,p_r$, 
\[
z = \sum_{i=1}^r \sigma_i p_i +\sum_{i<j}z_{ij}
\]
with $-1 =-(u|p_i)\leq \sigma_i = (z|p_i) \leq (u|p_i) = 1$.   As this is an orthogonal decomposition we have that 
\[
(z|\phi_H(x))  =  \sum_{i=1}^r\sigma_i\left( \frac{\lambda_i}{\sum_{j=1}^r \lambda_j} -  \frac{\lambda^{-1}_i}{\sum_{j=1}^r \lambda^{-1}_j}\right)
    < \sum_{i=1}^r\left( \frac{\lambda_i}{\sum_{j=1}^r \lambda_j}\right) + \sum_{i=1}^r\left( \frac{\lambda^{-1}_i}{\sum_{j=1}^r \lambda^{-1}_j}\right)
    =  2.
\]
This implies that $\frac{1}{2}\|\phi_H(x)\|_u^* = \frac{1}{2}\sup_{-u\leq z\leq u} (z|\phi_H(x))<1$, and hence $\phi_H(x)\in\mathrm{int}\, B_1^*$.

To prove the second assertion let $h$ be a horofunction given by $h(x) =\log M(y/x)+\log M(z/x^{-1})$, where $\|y\|_u=\|z\|_u =1$ and $(y|z)=0$. Write $y=\sum_{i\in I}\alpha_iq_i$ and $z=\sum_{j\in J}\beta_jq_j$. If we now let $q_I = \sum_{i\in I}q_i$ and $q_J = \sum_{j\in J}q_j$, then $-u\leq q_I-q_J\leq u$ and 
\[
\|\phi_H(h)\|_u^*\geq \frac{1}{2}(q_I-q_J|\phi_H(h)) = (1+1)/2=1.
\]
Moreover, for each $-u\leq w\leq u$ we have that 
\[
|(w|\phi_H(h))| \leq |(w|y/\mathrm{tr}(y))|+|(w|z/\mathrm{tr}(z))|\leq 
(u|y/\mathrm{tr}(y))+(u|z/\mathrm{tr}(z)) =2.
\]
Combining the inequalities shows that $\phi_H(h)\in \partial B_1^*$. 
\end{proof}

To prove injectivity of $\phi_H$ on $\Omega_V$ we need the following lemma, which is similar to Lemma \ref{lem:calc}. 
\begin{lemma}\label{lem:5.4}
Let $\mu_i\colon \mathbb{R}^r\to \mathbb{R}$, for $i=1,2$, be given by 
\[
\mu_1(x) = \sum_{i=1}^r e^{x_i}\mbox{\quad and \quad }\mu_2(x) = \sum_{i=1}^r e^{-x_i}\mbox{\quad for $x\in\mathbb{R}^r$, }
\]
 and let  $g\colon x\mapsto \log \mu_1(x)+\log\mu_2(x)$. If $x,y\in\mathbb{R}^r$ are such that $y\neq x+c(1,\ldots,1)$ for all $c\in\mathbb{R}$, then $\nabla g(x)\neq \nabla g(y)$. 
\end{lemma}
\begin{proof}
For $0<t<1$, $p =1/t$ and $q=1/(1-t)$ we have, by H\"{o}lder's inequality, that
\[
\mu_1(tx+(1-t)y) = \sum_i e^{tx_i}e^{(1-t)y_i} \leq \left(\sum_i e^{x_i}\right)^{1/p}\left(\sum_i e^{y_i}\right)^{1/q} = \mu_1(x)^t\mu_1(y)^{1-t},
\]
and we have equality if and only if there exists a $C_1>0$ such that $e^{y_i} =(e^{(1-t)y_i})^q = C_1(e^{tx_i})^p=C_1e^{x_i}$ for all $i$, which is equivalent to $y_i=x_i+c_1$ for all $i$. 

Likewise, 
\[
\mu_2(tx+(1-t)y) = \mu_2(x)^t\mu_2(y)^{1-t}
\]
and we have equality if and only if  $y_i=x_i+c_2$ for all $i$. 

It follows that $g\colon x\mapsto \log \mu_1(x)+\log\mu_2(x)$ satisfies
\[
g(tx+(1-t)y) \leq tg(x)+(1-t)g(y)
\]
for $x,y\in\mathbb{R}^r$ and $0<t<1$. Moreover, we have equality if, and only if, there exists $c\in\mathbb{R}$ such that $y_i=x_i+c$ for all $i$.   This implies that if  $x,y\in\mathbb{R}^r$ are such that $y\neq x+c(1,\ldots,1)$ for all $c\in\mathbb{R}$, then 
$g(x)-g(y)>\nabla g(y)\cdot (x-y)$ and $g(y)-g(x)>\nabla g(x)\cdot (y-x)$. So, 
\[
0>(\nabla g(y)-\nabla g(x))\cdot (x-y),
\]
and hence $\nabla g(y)\neq \nabla g(x)$. 
\end{proof}
\begin{lemma}\label{lem:5.5} The map $\phi_H$ is injective on $\Omega_V$.
\end{lemma}
\begin{proof}
Suppose that $\phi_H(x)=\phi_H(y)$, where $x=\sum_{i=1}^r \lambda_i p_i$ and $y=\sum_{i=1}^r \mu_i q_i$ in $\Omega_V$.  Note that $0<\lambda_i,\mu_i$ for all $i$ and $(x|u)= \mathrm{tr}(x) =r = \mathrm{tr}(y)= (y|u)$.
After possibly relabelling we can write 
\[
\phi_H(x) = \sum_{i=1}^r \left( \frac{\lambda_i}{\sum_{j=1}^r \lambda_j} - \frac{\lambda^{-1}_i}{\sum_{j=1}^r \lambda^{-1}_j} \right)p_i =\sum_{i=1}^r \alpha_i p_i
\]
and 
\[
\phi_H(y) = \sum_{i=1}^r \left( \frac{\mu_i}{\sum_{j=1}^r \mu_j} - \frac{\mu^{-1}_i}{\sum_{j=1}^r \mu^{-1}_j} \right)q_i =\sum_{i=1}^r \beta_i q_i,
\]
where  $\alpha_1\leq \ldots\leq \alpha_r$ and $\beta_1\leq \ldots\leq \beta_r$. By the spectral theorem (version 2) \cite{FK} we conclude that $\alpha_i =\beta_i$ for all $i$. 

Consider the injective map $\mathrm{Log}\colon\mathrm{int}\,\mathbb{R}_+^r\to \mathbb{R}^r$ given by $\mathrm{Log}(\gamma)= (\log \gamma_1,\ldots,\log \gamma_r)$. Let $\Delta = \{\gamma\in \mathrm{int}\,\mathbb{R}^r_+\colon \sum_{i=1}^r \gamma_i=r\}$.  The map $ (\nabla g)\circ \mathrm{Log}$ is injective on $\Delta$ by Lemma \ref{lem:5.4} and 
\[
\nabla g(\mathrm{Log}(\gamma)) = \left(\frac{\gamma_1}{\sum_{i=1}^r \gamma_i} - \frac{\gamma^{-1}_1}{\sum_{i=1}^r \gamma^{-1}_i}, \ldots, \frac{\gamma_r}{\sum_{i=1}^r \gamma_i} - \frac{\gamma^{-1}_r}{\sum_{i=1}^r \gamma^{-1}_i}\right).
\]

Writing $\lambda =(\lambda_1,\ldots,\lambda_r)$ and $\mu = (\mu_1,\ldots,\mu_r)$, we have that $\lambda,\mu\in \Delta$ and  
\[
\nabla g(\mathrm{Log}(\lambda)) =(\alpha_1,\ldots,\alpha_r) = (\beta_1,\ldots,\beta_r) = \nabla g(\mathrm{Log}(\mu)),  
\]
so that $\lambda = \mu$.  

As $(\nabla g)\circ\mathrm{Log}$ is injective on $\Delta$, we also know that $\alpha_k=\alpha_{k+1}$ if and only if $\lambda_k =\lambda_{k+1}$. Likewise,  $\beta_k=\beta_{k+1}$ if and only if $\mu_k =\mu_{k+1}$. From the spectral theorem (version 1) \cite{FK} we now conclude that $x=y$.
\end{proof}

In the next couple of lemmas we show that $\phi_H$ is onto. 
\begin{lemma}\label{lem:5.6}
The map $\phi_H$ maps $\Omega_V$ onto $\mathrm{int}\, B_1^*$. 
\end{lemma}
\begin{proof}
Note that $\Omega_V$ is an open set of the affine space $\{x\in V\colon \mathrm{tr}(x) =r\}$ which has dimension $\dim V-1$. Also $B_1^*\subset \mathbb{R}u^\perp$ has dimension $\dim V - 1$. As $\phi_H$ is a continuous injection from $\Omega_V$ into $\mathrm{int}\, B_1^*$ by Lemmas \ref{lem:5.3} and \ref{lem:5.5}, we know that $\phi_H(\Omega_V)$ is a open subset of $\mathrm{int}\, B_1^*$ by Brouwer's invariance of domain theorem.  We now argue by contradiction. So, suppose that  $\phi_V(\Omega_V)\neq \mathrm{int}\, B_1^*$.  Then there exists a $w\in \partial \phi_H(\Omega_V)\cap \mathrm{int}\, B_1^*$, as otherwise $\phi_H(\Omega_V)$ is closed and open, which would imply that $\mathrm{int}\,B_1^*$ is the disjoint union of two nonempty open sets contradicting the connectedness of $\mathrm{int}\,B_1^*$. So let $w\in\partial \phi_H(\Omega_V)\cap \mathrm{int}\, B_1^*$ and   let $(v_n)$ in $\Omega_V$ be such that $\phi_H(v_n)\to w$.  

As $\phi_H$ is continuous on $\Omega_V$, we may assume that $d_H(v_n,u)\to\infty$. After taking a subsequence, we may also assume that $v_n\to v\in\partial \Omega_V$. Now let $y_n = v_n/\|v_n\|_u$ and set $y = v/\|v\|_u$. Furthermore let $z_n = y_n^{-1}/\|y_n^{-1}\|_u$ After taking subsequences we may assume that $z_n\to z\in\partial V_+$ and $y_n\to y\in\partial V_+$, so $\|y\|_u=\|z\|_u=1$. 
As $y_n\bullet z_n = u/\|y_n^{-1}\|_u \to 0$, we find that $y\bullet z = 0$, which implies that  $(y|z) = 0$. 

Using the spectral decomposition we write $y_n =\sum_{i=1}^r \lambda_i^np_i^n$ and $ y =\sum_{i\in I}\lambda_ip_i$ where $\lambda_i>0$ for all $i\in I$.  Likewise we let $z_n = \sum_{i=1}^r \mu_i^np_i^n$ and $z =\sum_{j\in J}\mu_jp_j$ with $\mu_j>0$ for all $j\in J$. Note that $\mu_i^n = (\lambda_i^n)^{-1}/\|y^{-1}_n\|_u$. 

Then 
\[
\phi_h(v_n) = \frac{\sum_{i=1}^r \lambda_i^np_i^n}{\sum_{k=1}^r \lambda_k^n} -  
\frac{\sum_{i=1}^r(\lambda_i^n)^{-1}p_i^n}{\sum_{k=1}^r (\lambda_k^n)^{-1}} = 
\frac{\sum_{i=1}^r \lambda_i^np_i^n}{\sum_{k=1}^r \lambda_k^n} -  
\frac{\sum_{i=1}^r \mu_i^np_i^n}{\sum_{k=1}^r \mu^n_k}\to \frac{\sum_{i\in I} \lambda_ip_i}{\sum_{k\in I} \lambda_k} -  
\frac{\sum_{j\in J} \mu_jp_j}{\sum_{k\in J} \mu_j} = w.
\]

Now let $w^* = \sum_{i\in I}p_i - \sum_{j\in J}p_j$ and note that $-u\leq w^*\leq u$, as $(y|z)=0$. We find that 
\[
\frac{1}{2}\|w\|^*_u \geq \frac{1}{2}(w|w^*) = (1+1)/2 =1,
\]
and hence $w\in\partial B_1^*$, which is a contradiction.
\end{proof}

\begin{lemma}\label{lem:5.7} The map $\phi_H$ maps $\partial\overline{\Omega}_V^h$ onto $\partial B_1^*$. 
\end{lemma}
\begin{proof}
We know from Lemma \ref{lem:5.3} that $\phi_H$ maps $\partial\overline{\Omega}_V^h$ into $\partial B_1^*$. To prove that it is onto let $w\in\partial B_1^*$. Then there exists a face, say 
\[
A_{p,q} = 2\mathrm{conv}\, ( (U_p(V)\cap S(V))\cup(U_q(V)\cap -S(V)))\cap \mathbb{R}u^\perp
\]
where $p$ and $q$ are orthogonal idempotents, such that $w$ is in the relative interior of $A_{p,q}$, as $B_1^*$ is the disjoint union of the relative interiors of its faces \cite[Theorem 18.2]{Rock}.  
So, 
\[
w= \sum_{i\in I} \alpha_i p_i -\sum_{j\in J}\beta_j q_j,
\]
where $\alpha_i>0$ for all $i\in I$, $\beta_j>0$ for all $j\in J$, and $\sum_{i\in I}\alpha_i +\sum_{j\in J} \beta_j =2$. Moreover, $\sum_{i\in I}p_i = p$ and $\sum_{j\in J}q_j =q$.  

As $w\in \mathbb{R}u^\perp$, we have that $0=(u|w) = \sum_{i\in I}\alpha_i -\sum_{j\in J} \beta_j$, and hence 
$\sum_{i\in I}\alpha_i =\sum_{j\in J} \beta_j =1$. 

Put $\alpha^* =\max_{i\in I} \alpha_i$ and $\beta^*=\max_{j\in J} \beta_j$. Furthermore, for $i\in I$ set $\lambda_i = \alpha_i/\alpha^*$ and for $j\in J$ set $\mu_j= \beta_j/\beta^*$. Then 
\[
w = \left( \frac{\sum_{i\in I} \alpha_i p_i}{\sum_{k\in I} \alpha_k}\right) - \left( \frac{\sum_{j\in J} \beta_jq_j}{\sum_{k\in J} \beta_k}\right) = 
\left( \frac{\sum_{i\in I} \lambda_i p_i}{\sum_{k\in I} \lambda_k}\right) - \left( \frac{\sum_{j\in J} \mu_jq_j}{\sum_{k\in J} \mu_k}\right).  
\]
Note that $0<\lambda_i\leq 1$ for all $i\in I$ and $\max_{i\in I}\lambda_i =1$. Likewise, $0<\mu_j\leq 1$ for all $j\in J$ and $\max_{j\in J}\beta_j =1$.  

Now let $y = \sum_{i\in I}\lambda_i p_i$ and $z=\sum_{j\in J} \mu_jq_j$. Then $\|y\|_ u =\|z\|_u =1$ and $(y|z) = 0$. Furthermore, if we let $h\colon \Omega_V\to \mathbb{R}$ be given by 
\[
h(x) = \log M(y/x) +\log M(z/x^{-1})
\]
for $x\in \Omega_V$, then $h$ is a horofunction by Theorem \ref{Hilberthoro} and 
\[
\phi_H(h) = \left( \frac{\sum_{i\in I} \lambda_i p_i}{\sum_{k\in I} \lambda_k}\right) - \left( \frac{\sum_{j\in J} \mu_jq_j}{\sum_{k\in J} \mu_k}\right) =w.
\]
This completes the proof. 
\end{proof}
We already saw in Lemma \ref{lem:5.5} that $\phi_H$ is injective on $\Omega_V$. The next lemma shows that $\phi_H$ is injective on $\overline{\Omega}_V^h$. 
\begin{lemma}
The map $\phi_H\colon \overline{\Omega}_V^h\to B_1^*$ is injective.
\end{lemma}
\begin{proof}
We know from Lemmas \ref{lem:5.5}, \ref{lem:5.6} and \ref{lem:5.7} that $\phi_H$ is injective on $\Omega_V$, $\phi_H$ maps $\Omega_V$ onto $\mathrm{int}\, B_1^*$, and $\phi_H(\partial \overline{\Omega}_V^h)\subseteq \partial B_1^*$. So to prove that $\phi_H$ is injective on $\overline {\Omega}_V^h$, it remains to show that if $h,h'\in\partial \overline {\Omega}_V^h$ and $\phi_H(h)=\phi_H(h')$, then $h=h'$.  

Suppose  
\[
h(x) =\log M(y/x) +\log M(z/x^{-1})\mbox{\quad and\quad }h'(x)= \log M(y'/x)+\log M(z'/x^{-1})
\]
for all $x\in\Omega_V$. Then 
\[
\phi_H(h) = \frac{y}{\mathrm{tr}(y)} - \frac{z}{\mathrm{tr}(z)} = \frac{y'}{\mathrm{tr}(y')} - \frac{z'}{\mathrm{tr}(z')} =\phi_H(h').
\]
Using the fact that the orthogonal decomposition of an element in $V$ is unique, see \cite[Proposition 1.26]{AS1}, we conclude that 
\[
\frac{y}{\mathrm{tr}(y)} = \frac{y'}{\mathrm{tr}(y')} \mbox{\quad and\quad } \frac{z}{\mathrm{tr}(z)} = \frac{z'}{\mathrm{tr}(z')}.
\]
As $\|y\|_u=\|y'\|_u=1$, we get that $\mathrm{tr}(y)=\mathrm{tr}(y')$, and hence $y=y'$. Likewise, $\|z\|_u=\|z'\|_u=1$ implies that $z=z'$. Thus, $h=h'$, which completes the proof.
\end{proof}
\subsection{Proof of Theorem \ref{Hilberthom}}
Before we prove Theorem \ref{Hilberthom}, we recall some terminology from Jordan theory. 
For $x,z\in V$ we let $[x,z]=\{y\in V\colon x\leq y\leq z\}$, which is called an {\em order-interval}. Given $y\in V_+$  we let 
\[
\mathrm{face}(y) = \{x\in V_+\colon x\leq \lambda y\mbox{ for some }\lambda\geq 0\}.
\]
In a Euclidean Jordan algebra $V$ every idempotent $p$ satisfies  
\[
\mathrm{face}(p) \cap [0,u] =[0,p],
\]
see \cite[Lemma 1.39]{AS1}. Also note that $y\sim y'$ if and only if $\mathrm{face}(y) =\mathrm{face}(y')$. 

\begin{proof}[Proof of Theorem \ref{Hilberthom}] 
We know from the results in the previous subsection that $\phi_H\colon\overline{\Omega}_V^h\to B_1^*$ is a bijection, which is continuous on $\Omega_V$. 

To prove continuity of $\phi_H$ on the whole of $\overline{\Omega}_V^h$ we first show that if $(v_n)$ in $\Omega_V$ is such that $h_{v_n}\to h\in\partial \overline{\Omega}_V^h$, then $\phi_H(v_n)\to \phi_H(h)$.  Let $h(x) = \log M(y/x)+\log M(z/x^{-1})$ for $x\in\Omega_V$, where $\|y\|_u=\|z\|_u =1$ and $(y|z)=0$.  For $n\geq 1$ let $y_n =v_n/\|v_n\|_u$ and note that $\phi_H(v_n) = \phi_H(y_n)$ for all $n$.  Let $w_k = \phi_H(v_{n_k})$, $k\geq 1$ be a subsequence of $(\phi_H(v_n))$.  We need to show that $(w_k)$ has a subsequence that converges to $\phi_H(h)$.

As $h$ is a horofunction and $(\Omega_V,d_H)$ is a proper metric space, we have that $d_H(v_n,u)=d_H(y_n,u)\to\infty$  by Lemma \ref{Rieffel}. 
It follows that $(y_{n_k})$ has a subsequence $(y_{k_m})$ with $y_{k_m}\to y'\in\partial V_+$ and $z_{k_m} = y^{-1}_{k_m}/\|y^{-1}_{k_m}\|_u  \to z'\in V_+$.  Note that as $y\in\partial V_+$, we have that $\|y_{k_m}^{-1}\|_u\to\infty$. This implies that 
\[
y'\bullet z' =\lim_{m\to\infty} y_{k_m}\bullet \frac{y_{k_m}^{-1}}{\|y_{k_m}^{-1}\|_u} = \lim_{m\to\infty} \frac{u}{\|y_{k_m}^{-1}\|_u}=0,
\]
which implies that $(y'|z')=0$ by \cite[III, Exercise 3.3]{FK}, and hence $z'\in\partial V_+$. 
 Moreover, for each $x\in\Omega_V$, 
\begin{eqnarray*}
\lim_{m\to\infty} h_{y_{k_m}}(x) & = &\lim_{m\to\infty} \log M(y_{k_m}/x) +\log M(x/y_{k_m})  - \log M(y_{k_m}/u)- \log M(u/y_{k_m})\\
	& = & \lim_{m\to\infty}\log M(y_{k_m}/x) +\log M(y^{-1}_{k_m}/x^{-1})  - \log \|y_{k_m}\|_u- \log M(y^{-1}_{k_m}/u)\\
	& = & \lim_{m\to\infty}\log M(y_{k_m}/x) +\log M(y^{-1}_{k_m}/x^{-1})  - \log \|y^{-1}_{k_m}\|_u\\
& = & \lim_{m\to\infty}\log M(y_{k_m}/x) +\log M(z_{k_m}/x^{-1}) \\
& = & \log M(y'/x)+\log M(z/x^{-1}).
\end{eqnarray*}
So, if we let $h'(x)=  \log M(y'/x)+\log M(z/x^{-1})$, then $h'$ is a horofunction by Theorem \ref{Hilberthoro} and $h_{y_{k_m}}\to h'$.  As $h=h'$, we know that $\delta(h,h') = d_H(y,y') +d_H(z,z') =0$, and hence 
$y=y'$ and $z=z'$. It follows that 
\[
\phi_H(v_{k_m})=\phi_H(y_{k_m}) = \frac{y_{k_m}}{\mathrm{tr}(y_{k_m})} - \frac{y^{-1}_{k_m}}{\mathrm{tr}(y^{-1}_{k_m})}\to \frac{y}{\mathrm{tr}(y)} -\frac{z}{\mathrm{tr}(z)}= \phi_H(h). 
\]

Recall that $\phi_H$ maps $\Omega_V$ into $\mathrm{int}\, B_1^*$  and $\phi_H$ maps $\partial \overline{\Omega}_V^h$ into $\partial B_1^*$ by Lemma \ref{lem:5.3}. So, to prove the continuity of $\phi_H$ it remains to show that if $(h_n)$ is a sequence in $\partial\overline{\Omega}_V^h$ converging to $h\in \partial\overline{\Omega}_V^h$, then $\phi_H(h_n)\to\phi_H(h)$.

Let $(\phi_H(h_{n_k}))$ be a subsequence of $(\phi_H(h_n))$. We need to show that it has a subsequence $(\phi_H(h_{k_m}))$ converging to $\phi_H(h)$. 
We know there exists $v_m,w_m\in\partial V_+$ with $\|v_m\|_u=\|w_m\|_u =1$ and $(v_m|w_m) =0$ such that 
\[
h_{k_m}(x) = \log M(v_m/x)+\log M(w_m/x^{-1})
\]
for $x\in\Omega_V$. By taking a further subsequence we may assume that $v_m\to v\in\partial V_+$ and $w_m\to w\in\partial V_+$. Then $\|v\|_u=\|w\|_u=1$ and $(v|w)=0$. Moreover,
\[
\log M(v_m/x)\to\log M(v/x)\mbox{\quad and\quad }\log M(w_m/x^{-1})\to \log M(w/x^{-1}), 
\]
for each $x\in\Omega_V$, as $y\mapsto M(y/x)$ is a continuous map on $V$, see \cite[Lemma 2.2]{LLNW}. Thus, $h_{k_m}\to h^*\in\partial\overline{\Omega}_V^h$, where 
\[
h^*(x) = \log M(v/x)+\log M(w/x^{-1}),
\]
by Theorem \ref{Hilberthoro}. As $h_n\to h$, we have that $h=h^*$.  This implies that $y=v$ and $z=w$, as otherwise $\delta(h,h^*)\neq 0$ by Proposition \ref{d_Hparts}. Thus, $v_m\to y$ and $w_m\to z$, and hence 
\[
\phi_H(h_{k_m}) = \frac{v_m}{\mathrm{tr}(v_m)} -\frac{w_m}{\mathrm{tr}(w_m)}\to 
\frac{y}{\mathrm{tr}(y)} -\frac{z}{\mathrm{tr}(z)}=\phi_H(h).
\]
This completes the proof of the continuity of $\phi_H$. 

Thus, $\phi_H$ is a continuous bijection from $\overline{\Omega}_V^h$ onto $B_1^*$. As 
$\overline{\Omega}_V^h$ is compact and  $B_1^*$ is Hausdorff, we conclude that $\phi_H$ is a homeomorphism. 

To complete the proof of the theorem it remains to show that $\phi_H$ maps each part onto the relative interior of a boundary face of $B_1^*$. Let $h(x) =\log M(y/x)+\log M(z/x^{-1})$ be a horofunction, where $y =\sum_{i\in I}\lambda_i p_i$ and $z=\sum_{j\in J}\mu_j p_j$ with $\lambda_i,\mu_j>0$  for all $i\in I$ and $j\in J$. 
Let $p_I=\sum_{i\in I} p_i$ and $p_J=\sum_{j\in J}p_j$. As $\phi_H$ is surjective, it suffices to show that $\phi_H$ maps $\mathcal{P}_h$ into the relative interior of 
\[
A_{p_I,p_J} = 2\mathrm{conv}\, ( (U_{p_I}(V)\cap S(V))\cup(U_{p_J}(V)\cap -S(V)))\cap \mathbb{R}u^\perp.
\]
So, let $h'\in\mathcal{P}_h$ where $h'(x) = \log M(y'/x)+\log M(z'/x^{-1})$ for $x\in\Omega_V$. Then $p_I\sim y\sim y'$ and $p_J\sim z\sim z'$. Using the spectral decomposition write $y'=\sum_{i\in I'} \alpha_iq_i$ and $z'=\sum_{j\in J'} \beta_jq_j$, where $\alpha_i>0$ for all $i\in I'$ and $\beta_j>0$ for all $j\in J'$. Now let $q_{I'} = \sum_{i\in I'} q_i$ and $q_{J'} = \sum_{j\in J'} q_j$.  It follows that $p_I\sim q_{I'}$ and $p_J \sim q_{J'}$.  So, 
$\mathrm{face}(p_I)=\mathrm{face}(q_{I'})$ and $\mathrm{face}(p_J) =\mathrm{face}(q_{J'})$.  As $\mathrm{face}(p_I) \cap [0,u] = [0,p_I]$ and 
$\mathrm{face}(q_{I'}) \cap [0,u] = [0,q_{I'}]$ by \cite[Lemma 1.39]{AS1}, we conclude that $p_I =q_{I'}$. In the same way we get that $p_J = q_{J'}$. As $\alpha_i>0$ for all $i\in I'$ and $\beta_j>0$ for all $j\in J'$, we have that 
\[
\phi_H(h')= \frac{y'}{\mathrm{tr}(y')} - \frac{z'}{\mathrm{tr}(z')} 
\]
is in the relative interior of $A_{q_{I'},q_{J'}} = A_{p_I,p_J}$.
\end{proof}

\section{Final remarks}  

It would be interesting to find a general class of simply connected smooth manifolds $M$ with a Finsler distance for which Problem \ref{keyproblem} has a positive solution. A common feature of the spaces considered in this paper is the property that the facial structure of the unit ball $\{x\in T_bM\colon F(b,v)\leq 1\}$ is the same for all $b\in M$. In particular, one could consider spaces where the $d_F$-isometry group of $M$ acts transitively on $M$. This is the case for all  normed spaces and the symmetric Hilbert geometries.  A second feature of the spaces considered here is that all horofunctions arise as limits of geodesics.  This property might be a useful  further assumption to make. 

Even if both these properties hold in a finite dimensional  normed space or a Hilbert geometry, then  it is not clear how one can define a homeomorphism for these spaces, despite the fact that we know all  horofunctions by Walsh \cite{Wa2,Wa1}. What made things work in our settings was the Jordan algebra structure and its associated spectral theory, which allowed us to give a more explicit description of the horofunctions and the parts of the horofunction boundary that gave a clear link with the dual norm.   
 
It is also worth noting that if both  $M$  and the normed space $(T_bM, \|\cdot\|_b)$ at the basepoint $b$  have a positive solution to Problem \ref{keyproblem}, then there exists a homeomorphism between the horofunction compactifications of these spaces that maps parts onto parts. It would be interesting to know if this connection exists more generally.  More specifically, one can ask the following general question. 
\begin{problem} Suppose $M$ is a simply connected smooth manifold with a Finsler distance, such that  the restriction of $F$ to the tangent space $T_bM$ at $b$ is a norm. When does there exist a homeomorphism between the horofunction compactification of $M$ with basepoint $b$ and the horofunction compactification of the normed space $(T_bM, \|\cdot\|_b)$, which maps parts onto parts?
\end{problem}
A solution to this problem would allow one to study the horofunction compactifications of these manifolds by analysing the  horofunction compactifications of finite dimensional normed spaces, which might be easier.

\end{document}